\documentclass[11pt]{article}
\usepackage{latexsym}
\usepackage{amsfonts,amssymb,amsmath}
\newtheorem{thm}{Theorem}[section]

\newtheorem{lem}[thm]{Lemma}
\newtheorem{pro}[thm]{Proposition}
\newtheorem{defn}[thm]{Definition}

\newtheorem{equat}[thm]{(\hspace{-.07in}}
\title{The ${\cal R}$- and ${\cal L}$-orders of the 
       Thompson-Higman monoid $M_{k,1}$  and their complexity }

\author{ Jean-Camille Birget  }

\date{\today}
\begin{document}
\maketitle

\begin{abstract}
We study the monoid generalization $M_{k,1}$ of the Thompson-Higman 
groups, and we characterize the ${\cal R}$- and the ${\cal L}$-preorder of
$M_{k,1}$. Although $M_{k,1}$ has only one non-zero $\cal J$-class and $k-1$
non-zero $\cal D$-classes, the ${\cal R}$- and the ${\cal L}$-preorder
are complicated; in particular, $<_{\cal R}$ is dense (even within an
$\cal L$-class), and $<_{\cal L}$ is dense (even within an $\cal R$-class).

We study the computational complexity of the ${\cal R}$- and the 
${\cal L}$-preorder.  When inputs are given by words over a finite 
generating set of $M_{k,1}$, the ${\cal R}$- and the ${\cal L}$-preorder 
decision problems are in {\sf P}.  The main result of the paper is that 
over a ``circuit-like'' generating set, the ${\cal R}$-preorder 
decision problem of $M_{k,1}$ is $\Pi_2^{\sf P}$-complete, whereas the
${\cal L}$-preorder decision problem is {\sf coNP}-complete.
We also prove related results about circuits: For combinational circuits,
the surjectiveness problem is $\Pi_2^{\sf P}$-complete, whereas the
injectiveness problem is {\sf coNP}-complete.
\end{abstract}


\section{Introduction}

The groups of Richard J.\ Thompson are well known and have been extensively
studied since their introduction in the 1960s \cite{Th0, McKTh, Th}. 
They were generalized by Graham Higman \cite{Hig74}. A classical survey is 
\cite{CFP}. Here we will follow the notation of \cite{BiThomps} (which is
similar to \cite{Scott}). 

The Thompson-Higman group $G_{k,1}$ (for any integer $k \geq 2$) is defined 
by taking all maximally extended right-ideal isomorphisms between essential 
right ideals of a free monoid $A^*$ (where $A$ is an alphabet with $k$
elements). This can be directly generalized to a monoid $M_{k,1}$ consisting 
of all right-ideal {\em homomorphisms} between right ideals of the free 
monoid $A^*$ \cite{BiThomMon}. Subsection 1.1 gives more detailed definitions.

The motivations for studying $G_{k,1}$ (and some of its subgroups) have 
been, first, its remarkable collection of properties (e.g., $G_{k,1}$ is 
finitely presented, it is simple, it contains all finite groups and all 
countable free groups), and second, its magical appearance in many contexts 
(see e.g.\ the many references in \cite{BiThomps, BiThomMon}).
The monoid $M_{k,1}$ also has remarkable properties: it is finitely 
generated, congruence-simple, ${\cal J}$-0-simple, it has exactly $k-1$ 
${\cal D}$-classes, and its word problem is in {\sf P} \cite{BiThomMon}. 
Moreover, there are strong connections between $M_{k,1}$ and circuit 
complexity \cite{BiFact, BiDistor}. Indeed, although $M_{k,1}$ is finitely
generated, it is also interesting to use generating sets of the form 
$\Gamma \cup \tau$ where $\Gamma$ is any finite generating set of $M_{k,1}$
and $\tau$ consists of all position transpositions in words in $A^*$ (the
exact definition of $\tau$ is given at the beginning of Section 4).
Then every combinational circuit $C$ (of size $|C|$) can be represented by 
a word of $M_{k,1}$ over $\Gamma \cup \tau$ of length $O(|C|)$; 
conversely, if a function $A^m \to A^n$ can be described by a word $w$ 
over $\Gamma \cup \tau$ of length $|w|$ then this function has a circuit 
of size $O(|w|^2)$ (Prop.\ 2.4 and Theorem 2.9 in \cite{BiDistor}).  
We call generating sets of $M_{k,1}$ of the form $\Gamma \cup \tau$  
{\it circuit-like} generating sets.

The monoid $M_{k,1}$ (and the Thompson-Higman group $G_{k,1}$ too) has a 
hybrid nature with respect to finiteness and infinity: On the one hand,
each element of $M_{k,1}$ has a finite description (by a finite function 
between finite prefix codes). However, $M_{k,1}$ is countably infinite, and 
its elements have a partial action on $A^*$; moreover, $M_{k,1}$ acts
faithfully on the Cantor space $A^{\omega}$, which is uncountable. 

The subject of this paper is to characterize the ${\cal R}$- and the
${\cal L}$-preorder of $M_{k,1}$. 
The results reflect the hybrid nature of $M_{k,1}$: 
The characterization of the ${\cal R}$- and ${\cal L}$-preorders is quite 
simple in terms of the action on the uncountable space $A^{\omega}$, but 
looks a little bit more complicated when 
formulated in terms of the pre-action on the countable set $A^*$. 
The ${\cal R}$- and ${\cal L}$-preorders are dense, i.e., whenever 
$\psi <_{\cal R} \varphi$ in $M_{k,1}$ then 
$\psi <_{\cal R} \chi <_{\cal R} \varphi$ for some $\chi \in M_{k,1}$.
We also study the computational complexity of the ${\cal R}$- and the
${\cal L}$-preorder of $M_{k,1}$. We prove that over a finite generating 
set the ${\cal R}$- and the ${\cal L}$-preorder decision problems of 
$M_{k,1}$ are in {\sf P}. The main result of the paper is that over a 
``circuit-like'' generating set, however, the ${\cal R}$-preorder decision 
problem of $M_{k,1}$ is $\Pi_2^{\sf P}$-complete, whereas the 
${\cal L}$-preorder decision problem is {\sf coNP}-complete. 
We also prove related results about circuits: For combinational circuits, 
the surjectiveness problem is $\Pi_2^{\sf P}$-complete, and the 
injectiveness problem is {\sf coNP}-complete.
We also prove that when $\psi <_{\cal R} \varphi$ over a circuit-like 
generating set then the length of right-multipliers cannot be polynomially
bounded, unless the polynomial hierarchy collapses. Also, for surjective
elements of $M_{k,1}$ and for surjective combinational circuits, right 
inverses cannot have polynomially bounded size, unless the polynomial 
hierarchy collapses.
 
This paper is a continuation of \cite{BiThomMon}, where the monoid
generalizations $M_{k,1}$ of the Thompson-Higman groups $G_{k,1}$ was 
introduced. The present paper focuses on the $\cal R$- and $\cal L$-orders 
of $M_{k,1}$, and the complexity of problems associated with the $\cal R$- 
and $\cal L$-orders of $M_{k,1}$.  
The $\cal J$-order and the $\cal D$-relation of $M_{k,1}$, and their
complexity, will be studied in \cite{BiDrelJord}. 


\subsection{Definition of the Thompson-Higman groups and monoids}

We first review some of the material from \cite{BiThomMon}.
Let $A$ be a finite alphabet of cardinality $k$, and let $A^*$ denote the 
free monoid over $A$ (consisting of all finite sequences of elements of 
$A$); elements of $A^*$ are called words. The {\it empty word} is denoted 
by $\varepsilon$. The {\it length} of $w \in A^*$ is denoted by $|w|$. 
For $u,v \in A^*$, the {\it concatenation} is denoted by $uv$ or by 
$u \cdot v$; more generally, the concatenation of $B, C \subseteq A^*$ is
$BC = \{uv : u \in B, v \in C\}$. 

We say that two sets $X$ and $Y$ {\em intersect} iff
$X \cap Y \neq \varnothing$.
A {\it right ideal} of $A^*$ is a subset 
$R \subseteq A^*$ such that $RA^* \subseteq R$. 
A right ideal $R$ is said to be {\it essential} iff $R$ intersects every 
right ideal of $A^*$.
For right ideals $R' \subseteq R \subseteq A^*$ we say that $R'$ is
{\em essential in} $R$ (or that $R'$ is an essentially equal right subideal 
of $R$) iff  $R'$ intersects all right subideals of $R$.

We say that a right ideal is generated by a  set $C$ iff $R$ is the
intersection of all right ideals that contain $C$; equivalently, $R = CA^*$.  
We say that $u \in A^*$ is a {\it prefix} of $v \in A^*$ iff $uz = v$ for 
some $z \in A^*$, and we write $u \ {\sf pref} \ v$. 
A {\it prefix code} is a subset $C \subseteq A^*$ such that no element of 
$C$ is a prefix of another element of $C$. 
A prefix code is {\it maximal} iff it is not a strict subset of another 
prefix code.
It is easy to prove that a right ideal $R$ has a unique minimal (under
inclusion) generating set, and that this minimal generating set is a prefix 
code; this prefix code is maximal iff $R$ is an essential right ideal.

It is often helpful to picture $A^*$ as the infinite $k$-ary tree, 
with vertex set $A^*$ and edge set $\{(v,va): v \in A^*, a \in A\}$. This 
is the right Cayley graph of the monoid $A^*$ with generating set $A$, and 
we will call it {\it the tree of} $A^*$. We turn this into a rooted tree 
by choosing the empty word $\varepsilon$ as the root.  
A {\it path} in this tree will always be taken to be directed away from 
the root.
In general, a rooted tree is called a {\it $k$-ary tree} iff every vertex 
has $\leq k$ children.
A $k$-ary tree is called {\it saturated} iff every non-leaf vertex has
exactly $k$ children.

A word $v$ is a prefix of a word $w$ iff $v$ is is an ancestor of $w$ in 
the tree of $A^*$.
A set $P$ is a prefix code iff no two elements of $P$ are on the same path.
A set $R$ is a right ideal iff any path that starts in $R$ has all its
vertices in $R$.
A finitely generated right ideal $R$ is essential iff every infinite path
of the tree reaches $R$ (and then stays in it from there on).
Similarly, for two finitely generated right ideals $R' \subset R$, $R'$ is
essential in $R$ iff any infinite path starting in $R$ also intersects
$R'$ (and then stays in $R'$ from there on). 
A finite prefix code $P$ is maximal iff any infinite path starting at the 
root intersects $P$.

A {\em right ideal homomorphism} of $A^*$ is a function
$\varphi: R_1 \to A^*$ with domain $R_1$ such that $R_1$ is a right ideal 
of $A^*$, and such that for all $x_1 \in R_1$ and all $w \in A^*$:
 \ $\varphi(x_1w) = \varphi(x_1) \ w$.
It is easy to prove that the image set of a right ideal homomorphism is 
a right ideal, which is finitely generated (as a right ideal) if the domain
$R_1$ is finitely generated.
Let $f:A^* \to A^*$ be any partial function; then ${\sf Dom}(f)$ denotes the 
domain and ${\sf Im}(f)$ denotes the image (range) of $f$.

For any set $X \subseteq A^*$ we denote the {\em partial identity} of $X$ 
by ${\sf id}_X$, i.e., ${\sf id}_X(x) = x$ for all $x \in X$, and 
${\sf id}_X(x)$ is undefined when
$x \not\in X$. For any sets $X,Y \subseteq A^*$ we have
 \ ${\sf id}_X \circ {\sf id}_Y = {\sf id}_Y \circ {\sf id}_X = $
${\sf id}_{X \cap Y}$ .

Note that we write the action of functions on the left of the argument 
and thus, functions are read and composed from right to left.

\medskip

Let $\varphi: R_1 \to R_2$ be a right ideal homomorphism with 
$R_1 = {\sf Dom}(\varphi)$ and $R_2 = {\sf Im}(\varphi)$. Then $\varphi$ 
can be described by a total surjective function $P_1 \to S_2$, where $P_1$ 
is the prefix code (not necessarily maximal) that generates $R_1$ as a right 
ideal, and $S_2$ is a set (not necessarily a prefix code) that generates 
$R_2$ as a right ideal. This function $P_1 \to S_2$, which is just the 
restriction of $\varphi$ to $P_1$, is called the {\em table} 
of $\varphi$. The prefix code $P_1$ is called the {\em domain code} of 
$\varphi$ and we write $P_1 = {\rm domC}(\varphi)$.  When $S_2$ is a prefix 
code we call $S_2$ the {\em image code} of $\varphi$ and we write 
$S_2 = {\rm imC}(\varphi)$.

We say that right ideal homomorphism $\Phi: R'_1 \to A^*$ is an 
{\em essentially equal restriction} of a right ideal homomorphism
$\varphi: R_1 \to A^*$ (or, equivalently, that $\varphi$ is an 
{\em essentially equal extension} of $\Phi$) iff $R'_1$ is essential in 
$R_1$, and for all $x'_1 \in R'_1$: \ $\varphi(x'_1) = \Phi(x'_1)$. 
In the earlier papers \cite{BiThomps, BiCoNP, BiFact, BiThomMon, BiDistor}
we used the terms ``essential restriction (extension)'' instead of 
``essentially equal restriction (extension)''; the new terminology is more
precise and will help prevent mix-ups with other concepts that will be
introduced later in this paper. 

The following is a crucial fact:
Every homomorphism $\varphi$ between finitely generated right ideals of
$A^*$ has a {\em unique} maximal essentially equal extension, denoted 
${\sf max}(\varphi)$  (Prop.\ 1.2(2) in \cite{BiThomMon}). 
Another interesting fact (remark after Prop.\ 1.2 in \cite{BiThomMon}): 
Every right ideal homomorphism $\varphi$ has an essentially equal 
restriction $\varphi'$ whose table $P' \to Q'$ is such that both $P'$ and 
$Q'$ are prefix codes.

We can now define the {\it Higman-Thompson monoid} $M_{k,1}$: As a set, 
$M_{k,1}$ consists of all homomorphisms (between finitely generated right 
ideals of $A^*$) that have been maximally essentially equally extended. 
In other words, as a set, 

\smallskip

 \ \ \ $M_{k,1} = \{ {\sf max}(\varphi) : \varphi$ is a homomorphism 
between finitely generated right ideals of $A^* \}$. 

\smallskip

\noindent
The multiplication is composition followed by maximal essentially equal 
extension. This multiplication is associative (Prop.\ 1.4 in 
\cite{BiThomMon}).  
Similarly, we define the uncountable monoid ${\cal M}_{k,1}$ exactly like
$M_{k,1}$ except that we allow all homomorphisms between (not necessarily 
finitely generated) right ideals of $A^*$. 

An important motivation for considering $M_{k,1}$ is that its elements are
closely related to combinational circuits, as mentioned before.
Another possible motivation is the connection with
C$^*$-algebras; indeed, $M_{k,1}$ is a submonoid of the multiplicative part 
of the Cuntz algebra ${\cal O}_k$ \cite{BiThomMon}. 

In \cite{BiThomMon} it was proved that $M_{k,1}$ is finitely generated, and 
that its word problem is decidable in deterministic polynomial time; moreover, 
$M_{k,1}$ is congruence-simple, i.e., the only congruences on $M_{k,1}$ are 
the trivial congruence (in which all of $M_{k,1}$ is one class), and the 
discrete congruence (in which each element constitutes a class). 
It was also proved that the Higman-Thompson group $G_{k,1}$ is the group of
units of $M_{k,1}$ (i.e., the set of invertible elements).

To describe the structure of a semigroup, the {\it Green relations} play an
important role; they determine the left, right, and two-sided ideal structure
of the semigroup. We consider the following Green relations: 
 \ $\leq_{\cal J}$ (the $\cal J$-preorder), 
$\leq_{\cal L}$ (the $\cal L$-preorder),  
$\leq_{\cal R}$ (the $\cal R$-preorder),  and $\equiv_{\cal D}$ (the
$\cal D$-equivalence relation). From $\leq_{\cal J}$, $\leq_{\cal L}$, and
$\leq_{\cal R}$ one also defines $\equiv_{\cal J}$ (meaning ``$\leq_{\cal J}$
{\sc and} $\geq_{\cal J}$''), and similarly, $\equiv_{\cal L}$, and
$\equiv_{\cal R}$.
In this paper we will not consider the Green relations $\leq_{\cal H}$ and
$\equiv_{\cal H}$.   See e.g.\ \cite{CliffPres, Grillet} for more 
information on the Green relations.
For any $u,v \in M$ (where $M$ is a monoid) we have, by definition, 
$u \leq_{\cal J} v$ iff every ideal of $M$ containing $v$ also contains $u$;
equivalently, $u \leq_{\cal J} v$ iff there exist $x,y \in M$ such that 
$u = xvy$. 
Similarly, $u \leq_{\cal L} v$ iff any left ideal of $M$ containing $v$ 
also contains $u$; equivalently, there exists $x \in M$ such that $u = xv$; 
the definition of $\leq_{\cal R}$ is similar.  
By definition, $u \equiv_{\cal D} v$ iff there exists $s \in M$ such that
$u \equiv_{\cal R} s \equiv_{\cal L} v$; this is equivalent to the existence
of $t \in M$ such that $u \equiv_{\cal L} t \equiv_{\cal R} v$.

In \cite{BiThomMon} it was proved that $M_{k,1}$ has only one non-zero 
$\equiv_{\cal J}$-class, and that $M_{k,1}$ has exactly $k-1$ non-zero 
$\equiv_{\cal D}$-classes.
These results mean that $M_{k,1}$ has almost no structure as far as 
congruences, $\leq_{\cal J}$, $\equiv_{\cal J}$ and $\equiv_{\cal D}$ are 
concerned.
However, we will see in this paper that $\leq_{\cal L}$, $\leq_{\cal R}$
have a complicated structure.


\subsection{ Cantor space and topological aspects of the Thompson-Higman 
monoids }

For the study of $\leq_{\cal L}$ and $\leq_{\cal R}$ in $M_{k,1}$ it will 
be useful to consider the action of $M_{k,1}$ on infinite words, i.e, on 
the Cantor space $A^{\omega}$; this is somewhat similar to Thompson's 
\cite{Th} original definition of the group $G_{2,1}$.
We also call the elements of $A^{\omega}$ the {\bf ends} of the tree of 
$A^*$, or the ends of $A^*$. For an end $z = (z_n: n \geq 1)$ (with 
$z_n \in A$ for all $n$), we call
$\{ z_1 \ldots z_n : n \geq 1\} \cup \{ \varepsilon\}$ the set of prefixes 
of $z$. For notational convenience, we will define $\omega$ to consist of
the positive integers (starting at 1).

Let $R$ be a finitely generated right ideal of $A^*$ and let 
$z = (z_n: n \in \omega) \in A^{\omega}$.
We say that {\em the end $z$ belongs to the right ideal} $R$ iff 
$z_1 \ldots z_m \in R$ for some $m \geq 1$; 
since $R$ is a right ideal, this is equivalent to the existence of 
$m \geq 0$ such that for all $n \geq m$: \ $z_1 \ldots z_n \in R$. 
The set of all ends in a finitely generated right
ideal $R \subseteq A^*$ is denoted by ${\sf ends}(R)$.

The set of ends $A^{\omega}$ can be given the {\bf Cantor space topology}, 
defined by taking $\{ v A^{\omega} : v \in A^{\omega}\}$ as a base of open 
sets. Then a set is {\it open} iff it is of the form $PA^{\omega}$, 
where $P$ is any subset of $A^*$; moreover, $P$ can be taken to be a prefix 
code. Hence, for every right ideal $R$ the set ${\sf ends}(R)$ is open.
Every subset $S$ of $A^{\omega}$ can be written as $S = PA^{\omega} \cup E$, 
with $E \cap PA^{\omega} = \varnothing$, where $P \subset A^*$ is
countable, and $PA^{\omega}$ is the interior of $S$; moreover, here $P$ can 
be chosen to be a prefix code; 
the set $E \subset A^{\omega}$ is the (countable) set of isolated elements 
of $S$. By definition, $w \in S \subseteq A^{\omega}$ is an isolated point 
of $S$ if $w$ has a neighborhood $N_w$ such that $S \cap N_w = \{w\}$.
A set $S \subseteq A^{\omega}$ is {\it closed} (or, equivalently, 
compact) iff $S$ can be written as $S = PA^{\omega} \cup E$, with $P$ and 
$E$ as above, with the additional condition that $P$ and $E$ are finite.
Finally, a set $S \subseteq A^{\omega}$ is {\it clopen} iff 
$S = PA^{\omega}$ for some finite set $P \subset A^*$; moreover, here $P$ 
can be taken to be a finite prefix code (see e.g.\ \cite{PerrinPin}).
Hence for any right ideal $R$ we have: \ $R$ is finitely generated 
iff ${\sf ends}(R)$ is clopen in the Cantor space $A^{\omega}$.
For a right ideal $R$, ${\sf ends}(R)$ is always open; hence, 
since $A^{\omega}$ is compact, we also have for any right ideal $R$:
 \ {\it $R$ is finitely generated \ iff \ ${\sf ends}(R)$ is compact.}

This yields a topological characterization of the Higman-Thompson monoid 
$M_{k,1}$ within ${\cal M}_{k,1}$. For every $\varphi \in {\cal M}_{k,1}$ 
the following are equivalent: \ (1) \ $\varphi \in M_{k,1}$, 
 \ (2) \ ${\sf ends}({\sf Dom}(\varphi))$ is compact, 
 \ (3) \ ${\sf ends}({\sf Im}(\varphi))$ is compact.

For any set $S \subseteq A^{\omega}$ we denote the {\it closure} of $S$ by 
${\sf clos}(S)$. By definition of closure we have: \ $w \in {\sf clos}(S)$ 
iff for every prefix $p$ of $w$, $p A^{\omega}$ intersects $S$. Indeed, 
$\{ p A^{\omega} : $ is a prefix of $w\}$ is a neighborhood base of $w$. 
In particular we also have for any right ideal $R \subseteq A^*$ and any 
$w \in A^{\omega}$: \ $w \in {\sf clos}({\sf ends}(R))$ \ iff 
 \ there exists a prefix $p$ of $w$ such that for all $s \in A^*$,
 \ $ps \in R$.

Any element of $\varphi \in {\cal M}_{k,1}$, represented by a homomorphism 
between right ideals of $A^*$, can be extended to a partial function 
$\Phi$ on $A^* \cup A^{\omega}$; indeed, if $\varphi(u) = v$ for 
$u,v \in A^*$ then we define $\Phi(uw) = vw$ for all $w \in A^{\omega}$. 
The extension $\Phi$ can now be restricted to a partial function $\Phi'$ on
$A^{\omega}$; $\Phi'$ is well defined as a partial function, and $\Phi'$
uniquely determines $\Phi$ and $\varphi$. Indeed, if $\varphi$ is a right 
ideal homomorphism of the form $PA^* \to QA^*$, where $P, Q \subset A^*$ are 
prefix codes, then $\Phi'$ is a total function $PA^{\omega} \to QA^{\omega}$;
and vice-versa.  Hence, ${\cal M}_{k,1}$ acts faithfully on $A^{\omega}$, and
when the elements of ${\cal M}_{k,1}$ are represented by partial functions on 
$A^{\omega}$, the multiplication is just composition (without any further 
maximal extension).
In this action on $A^{\omega}$, the elements of ${\cal M}_{k,1}$ are 
{\it continuous} partial functions, with open domains and ranges.
 
%
%
%

\subsection{Three descriptions of the Thompson-Higman groups and monoids}

The Thompson-Higman groups $G_{k,1}$ and monoids $M_{k,1}$ can be defined
in equivalent ways, by finite, countable, or uncountable structures.
Each description has advantages and drawbacks.

In the {\bf finite description}, every element of $G_{k,1}$ or $M_{k,1}$
is given by a function between finite prefix codes over the alphabet $A$. 
Such a function can be represented concretely by a finite table. However, 
the representation is not unique: ``extensions'' or ``restrictions'' of a 
table do not change the element of $M_{k,1}$ being represented. In a 
{\it restriction}, a table entry $(x, y)$ is replaced by the $k$ entries 
$(xa_1, ya_1), \ldots, (xa_k, ya_k)$, where $A = \{a_1, \ldots, a_k\}$; 
for $G_{k,1}$; an {\it extension} is the inverse of a restriction.  The 
notion of restriction of tables looks a little contrived, but it has a 
natural interpretation in the countable description below, and in the 
uncountable description, restrictions and extensions are not used at all.
For background, see Thompson's Def.\ 1.1 in \cite{Th}, Higman's pages 24-25 
in \cite{Hig74}, and Lemma 2.2 in \cite{BiThomps}; the functions between the
leaves of finite trees, in \cite{CFP}, are very similar to the description
by tables; for $M_{k,1}$, see  Prop.\ 1.2 in \cite{BiThomMon}. 
Nevertheless, every element of $M_{k,1}$ has a unique maximally extended
table. Multiplication is somewhat complicated in the finite description: 
first the tables have to be restricted so as to become composable (i.e., 
so that the image code of the first equals the domain code of the second), 
then they are composed, and finally the resulting table has to be maximally 
extended (if the unique representation is desired). Associativity is not 
obvious. The $\cal L$- and $\cal R$-preorders of $M_{k,1}$ are complicated 
to characterize in the finite description. 

In the {\bf countable description}, every element of $G_{k,1}$ or $M_{k,1}$
is given by right ideal homomorphisms between finitely generated right ideals
of $A^*$.  Again, the representation is not unique: extensions or restrictions
of a right ideal homomorphism do not change the element of $M_{k,1}$ being
represented. However, the concept of extension or restriction is a special 
case of the usual one (for partial functions), and not ad hoc (as it was 
for tables). 
Every element of $M_{k,1}$ can be represented by a unique maximally extended
right ideal homomorphism.
Multiplication is a little simpler than for tables: it is just the usual 
function composition, followed by maximal extension (if the unique 
representation is desired).
Associativity is not obvious. The characterization of the $\cal L$- and 
$\cal R$-preorders of $M_{k,1}$ is manageable but somewhat complicated.

In the {\bf uncountable description}, every element of $G_{k,1}$ or $M_{k,1}$
is given by a permutation, respectively a partial function on the Cantor
space $A^{\omega}$. 
Multiplication is simply composition of permutations or of partial functions.
However, to define which permutations or partial functions on $A^{\omega}$
belong to $G_{k,1}$ or $M_{k,1}$, the countable description needs to be 
referred to (at least indirectly). 
The characterization of the $\cal L$- and $\cal R$-preorders of $M_{k,1}$ is 
easy to state.

\subsection{Overview}

In Section 2 the $\cal R$-order of $M_{k,1}$ is characterized (in terms of 
image sets),
and in Section 3 the $\cal L$-order is characterized (in terms of right
congruences). 
In Section 4 the $\cal R$- and $\cal L$-orders are embedded into the 
idempotent order, and various density properties of the $\cal R$- and 
$\cal L$-orders are shown.
In Section 5 it is proved that the $\leq_{\cal R}$-decision problem of
$M_{k,1}$ is in {\sf P} when a finite generating set is used, and that it 
is $\Pi_2^{\sf P}$-complete with circuit-like generating sets; the
surjectiveness problem for combinational circuits is also proved to be 
$\Pi_2^{\sf P}$-complete.  
In Section 6 it is proved that the $\leq_{\cal L}$-decision problem of
$M_{k,1}$ is in {\sf P} when a finite generating set is used, and that it 
is {\sf coNP}-complete with circuit-like generating sets; the 
injectiveness problem for combinational circuits is also proved to be 
{\sf coNP}-complete.


\section{The ${\cal R}$-order of $M_{k,1}$  }

The monoid $M_{k,1}$ has some similarities with the monoid ${\sf PF}_X$ 
of all partial functions on a set $X$. In both cases, the ${\cal R}$-order
between functions is related to the inclusion order of the image sets. 
However, for $M_{k,1}$ the notion of ``inclusion'' has several different 
characterizations, and depends on whether we use the finite, countable, or
uncountable description of $M_{k,1}$.
We show that the ${\cal R}$-order of $M_{k,1}$ is dense. 

\medskip

\noindent {\bf Notation.} For a finite set $S \subset A^*$, the length of
a longest word in $S$ will be denoted by $\ell(S)$.

Recall that a tree is called $k$-ary iff every vertex has $\leq k$ 
children, and at least one vertex has exactly $k$ children. 
A $k$-ary tree is said to be {\it saturated} iff every non-leaf 
vertex has exactly $k$ children.

\subsection{Characterization of the ${\cal R}$-order of $M_{k,1}$ }

\begin{thm} \label{R_order} {\bf (${\cal R}$-order of $M_{k,1}$).} \ \
For all $\psi, \varphi \in M_{k,1}$ the following are equivalent:

\smallskip

\noindent {\bf (1)} \ \ \ $\psi(.) \ \leq_{\cal R} \ \varphi(.)$; 

\smallskip

\noindent {\bf (2)} \ \ \  
${\sf ends}({\sf Im}(\psi)) \ \subseteq \ {\sf ends}({\sf Im}(\varphi))$;

\smallskip

\noindent {\bf (3)} \ \ \ every right ideal of $A^*$ that intersects 
${\sf Im}(\psi)$ also intersects ${\sf Im}(\varphi)$;

\smallskip

\noindent {\bf (4)} \ \ \ every monogenic right ideal of $A^*$ that 
intersects ${\sf Im}(\psi)$ also intersects ${\sf Im}(\varphi)$;

\smallskip

\noindent {\bf (5)} \ \ \  every path in the tree of $A^*$, starting at 
$\varepsilon$, of length 
 \, ${\sf max}\{ \ell({\sf imC}(\varphi)), \ell({\sf imC}(\psi))\}$, 
 that 

 \ \ \ \ intersects ${\sf imC}(\psi)$, also intersects ${\sf imC}(\varphi)$;

\smallskip

\noindent {\bf (6)} \ \ \ for every $y \in {\sf imC}(\psi)$ we have: 

 \ \ \ \ $\bullet$ \ either ${\sf imC}(\varphi)$ contains a prefix of $y$,

 \ \ \ \ $\bullet$ \ or the subtree of the tree of $A^*$ with root $y$ and 
leaf-set \, $yA^* \cap {\sf imC}(\varphi)$ \, is saturated. 
\end{thm}
The characterizations above reflect the various representations of
$M_{k,1}$, (2) the uncountable representation, (3) and (4) the countable 
infinite representation, (5) and (6) the finite representation.
We first prove intermediate results, and then the Theorem.

\begin{lem} \label{equiv_R_order} \hspace{-.13in} {\bf .}

\noindent For any two right ideals
$R_1, R_2 \subseteq A^*$ the following are equivalent:

\smallskip

{\bf (a)} \ \ \   
 ${\sf clos}({\sf ends}(R_1)) \ \subseteq \ {\sf clos}({\sf ends}(R_2))$.

\smallskip

{\bf (b)} \ \ \ Every right ideal of $A^*$ that intersects $R_1$
also intersects $R_2$.

\smallskip

{\bf (c)} \ \ \ Every monogenic right ideal of $A^*$ that
intersects $R_1$ also intersects $R_2$.

\smallskip

\noindent If $R_2$ is finitely generated then {\bf (a)} is equivalent to 
the following:

\smallskip

{\bf (d)} \ \ \ ${\sf ends}(R_1) \ \subseteq \ {\sf ends}(R_2)$.

\smallskip

\noindent If $R_1$ and $R_2$ are generated by the prefix code $P_1$, 
respectively $P_2$, then {\bf (d)} is equivalent to the following
(where we allow ${\sf max}\{ \ell(P_1), \ell(P_2)\}$ to be infinite):

\smallskip

{\bf (e)} \ \ \ Every path in the tree of $A^*$, starting at
$\varepsilon$, of length ${\sf max}\{ \ell(P_1), \ell(P_2)\}$, that intersects

\hspace{.4in} $P_1$, also intersects $P_2$.

\smallskip

\noindent If $R_1$ and, $R_2$ are generated by the prefix code $P_1$,
respectively $P_2$, and if $P_2$ is finite,  
then {\bf (d)} is equivalent to the following:

\smallskip

{\bf (f)} \ \ \ For every $y \in P_1$ we have:

 \hspace{.4in} $\bullet$ \ either $P_2$ contains a prefix of $y$,

 \hspace{.4in} $\bullet$ \ or the subtree of the tree of $A^*$ with root $y$ 
and leaf-set $yA^* \cap P_2$  is saturated. 
\end{lem}
{\bf Proof.} 
 [(a) $\Rightarrow$ (b)] \ Let $R$ be any right ideal that intersects 
$R_1$, and let $w \in R_1 \cap R$. Then there is an end in 
${\sf ends}(R_1 \cap R)$ with $w$ as a prefix (in fact, all 
ends with prefix $w$ are in ${\sf ends}(R_1 \cap R)$). By (a), 
this end is also in ${\sf clos}({\sf ends}(R_2))$. 
Hence, some finite word of the form $wx$ belongs to $R_2$. Since $w \in R$
and since $R$ is a right ideal, $wx \in R$. So, $wx \in R \cap R_2$, hence
$R$ and $R_2$ intersect.   

\smallskip

\noindent [(b) $\Rightarrow$ (c)] \ This is trivial.

\smallskip

\noindent
[(c) $\Rightarrow$ (a)] \ Let $z = z_1 \ldots z_m \ldots z_n \ldots $  
$ \ \in {\sf clos}({\sf ends}(R_1))$. Hence, $z_1 \ldots z_m \in R_1$ 
for some $m \geq 0$. Hence, for all $n > m:$ \ 
$z_1 \ldots z_m \ldots z_n \in R_1$, so 
$z_1 \ldots z_m \ldots z_n \, A^*$ intersects $R_1$.
By (c) this implies that $z_1 \ldots z_m \ldots z_n \, A^*$ intersects 
$R_2$, for all $n > m$. So for each $n > m$ there exists $u_n \in A^*$  
such that $z_1 \ldots z_m \ldots z_n \, u_n \in \ R_2$.

Let $P_2$ be a finite generating set of $R_2$, i.e., $R_2 = P_2A^*$.
So for each $n > m$ there exists $u_n \in A^*$ and there exist 
$p_n \in P_2$ and $v_n \in A^*$ such that 
$z_1 \ldots z_m \ldots z_n \, u_n = p_n v_n$.
When $n$ is longer than the longest word in $P_2$ (such an $n$ exists since 
$P_2$ is finite), $p_n$ will be a prefix of $z_1 \ldots z_m \ldots z_n$. 
Then some element of $P_2$ ($\subset R_2$) is a prefix of the end $z$.
It follows that the end $z$ is in $R_2$.

\smallskip

\noindent
[(a) $\Leftrightarrow$ (d), when $R_2$ is finitely generated] 
  \ Trivially, (d) always implies (a).
Suppose now that (a) holds. When $R_2$ is finitely generated, 
${\sf ends}(R_2)$ is compact, so 
${\sf ends}(R_2) = {\sf clos}({\sf ends}(R_2))$. Then 
${\sf ends}(R_1) \subseteq {\sf clos}({\sf ends}(R_1)) \subseteq $
${\sf clos}({\sf ends}(R_2)) = {\sf ends}(R_2)$, so (d) follows. 

\smallskip

\noindent
[(d) $\Leftrightarrow$ (e)] \ 
This is fairly obvious.

\smallskip

\noindent
[(e) $\Leftarrow$ (f), when $R_2$ is finitely generated] \
Let $p$ be a long-enough path, intersecting $P_1$ at some element $y$.
By (f), if some element $x$ of $P_2$ is a prefix of $y$ then $x$ is also 
on the path $p$, hence $p$ also intersects $P_2$ (at $x$).
If, on the other hand, the subtree of $A^*$ with root $y$ and leaf-set 
$yA^* \cap P_2$ is saturated, then every way to continue the path $p$ 
beyond $y$ will lead to an intersection with $P_2$. 

\smallskip

\noindent
[(e) $\Rightarrow$ (f), when $R_2$ is finitely generated] \
Consider any $y \in P_1$ and consider any long path $p$ through $y$.
By (e), $p$ intersects $P_2$ at some point, say $x \in P_2$. If $x$ is
shorter than $y$ then $x$ is between $y$ and the root $\varepsilon$, hence
$x$ is a prefix of $y$, and no other path through $y$ intersects $P_2$ 
(indeed, there is only one path from the root to, and that path intersects
$P_2$).

If $x$ is longer than $y$ then the intersection of the paths through $y$ and
$P_2$ is $yA^* \cap P_2$. If the tree with root $y$ and leaves $yA^* \cap P_2$
were not saturated then some path through $y$ (away from the root)
could ``escape'' without intersecting $P_2$.  
 \ \ \ $\Box$

\bigskip

\noindent {\bf Remark R\ref{equiv_R_order}.} \  Part (d) of  
Lemma \ref{equiv_R_order} is not true if we omit the assumption that $R_2$
is {\em finitely generated}. For example, let \ $A = \{a,b\}$, 
 \ $R_1 = a \, \{a,b\}^*$, \ $R_2 = \{a^mb: m \geq 0\} \ \{a,b\}^*$ \   
($ = a^*b \, \{a,b\}^*$). Then every right ideal that intersects $R_1$ also 
intersects $R_2$. However, the end \ $a^{\omega} = (a^n : n \geq 0)$ \ is 
in $R_1$ but not in $R_2$.

\begin{defn} \label{essent_inclusion} \ For any two finitely generated
right ideals $R_1, R_2 \subseteq A^*$ we say that $R_1$ is 
{\bf end-included} in $R_2$ iff $R_1, R_2$ satisfy the equivalent 
properties of Lemma \ref{equiv_R_order}.

We denote this by \ $R_1 \subseteq_{\sf end} R_2$.
 
Two finitely generated right ideals $R_1, R_2 \subseteq A^*$ are called 
{\bf essentially equal} (or {\bf end-equal}) iff 
$R_1 \subseteq_{\sf end} R_2$ and $R_2 \subseteq_{\sf end} R_1$. \ We 
denote this by $R_1 =_{\sf ess} R_2$. 
\end{defn}

The relation $\subseteq_{\sf end}$ is a pre-order (reflexive and transitive)
on the set of finitely generated right ideals of $A^*$. Moreover,
$\subseteq_{\sf end}$ is a {\em lattice pre-order}, i.e., the order on the 
set of $=_{\sf ess}$-equivalence classes is a lattice order.
We have $R_1 \subseteq_{\sf end} R_2$ iff
$R_1 \ =_{\sf ess} \ R_1 \cap R_2$ iff
$R_1 \cup R_2 \ =_{\sf ess} \ R_2$. 
Clearly, $R_1 \subseteq_{\sf end} R_2$ does not imply $R_1 \subseteq R_2$,
and essential equality does not imply equality; e.g., 
all finitely generated {\em essential} right ideals of $A^*$  (as defined in
Subsection 1.1) are essentially equal.

\begin{lem} \label{righIdealRewrite} \           
Let $R_1 = P_1A^*$ and $R_2 = P_2A^*$ be right ideals, where $P_1$ and 
$P_2$ are finite prefix codes. Then $R_1 =_{\sf ess} R_2$ iff $P_2$
can be transformed into $P_1$ by a finite sequence of replacements steps
of the following form:

\smallskip

\noindent {\bf (r1)} \ \ \ For a finite prefix code $C$ and for $c \in C$,
replace $C$ by \ $(C - \{c\}) \cup cA$.
 
\smallskip

\noindent {\bf (r2)} \ \ \ For a finite prefix code $C'$ such that 
$cA \subseteq C'$ for some word $c$, replace $C'$ by \ $(C' - cA) \cup \{c\}$.  
\end{lem}
{\bf Proof.} This is straightforward. \ \ \ $\Box$

\begin{lem} \label{idealInter} {\rm (Lemma 3.3 of \cite{BiThomps}).} \
Let $P, Q, R \subseteq A^*$ be such that $PA^* \cap QA^* = RA^*$, and
$R$ is a prefix code. Then $R \subseteq P \cup Q$.

As a consequence, the intersection of two finitely generated right
ideals is finitely generated.  \ \ \ $\Box$
\end{lem}


\begin{defn} \label{riHom} \hspace{-.13in} {\bf .} 
 \ By ${\sf riHom}(A^*)$ we denote the set of right ideal homomorphisms 
between finitely generated right ideals of $A^*$. This is a monoid under
function composition.
\end{defn}
Recall (Definition 1.3 in \cite{BiThomMon}) that for 
$\psi, \psi_0 \in {\sf riHom}(A^*)$ we say that $\psi_0$ is an 
{\em essentially equal restriction} of $\psi$  iff $\psi_0$ is a 
restriction of $\psi$ and
${\sf Dom}(\psi_0) =_{\sf ess} {\sf Dom}(\psi)$.

\begin{lem} \label{ess_restr} \ Assume $\psi, \psi_0 \in {\sf riHom}(A^*)$
are such that  $\psi_0$ is an essentially equal restriction of $\psi$.
Then  \ ${\sf Im}(\psi) =_{\sf ess} {\sf Im}(\psi_0)$ \ and \
${\sf Dom}(\psi) =_{\sf ess} {\sf Dom}(\psi_0)$.
\end{lem}
{\bf Proof.} This follows directly from the definition of essentially
equal restriction.
 \ \ \ $\Box$

\begin{pro} \label{ess_equ_Im_Dom} \ If $\psi, \varphi \in {\sf riHom}(A^*)$
represent the same element of $M_{k,1}$ then
 \ ${\sf Im}(\psi) =_{\sf ess} {\sf Im}(\varphi)$ \ and 
 \ ${\sf Dom}(\psi) =_{\sf ess} {\sf Dom}(\varphi)$.
\end{pro} 
{\bf Proof.} $\psi$ and $\varphi$ represent the same element of $M_{k,1}$ 
iff ${\sf max}(\psi) = {\sf max}(\varphi)$. By Lemma \ref{ess_restr},
 \ ${\sf Im}(\psi) =_{\sf ess} {\sf Im}({\sf max}(\psi)) \ = \ $
${\sf Im}({\sf max}(\varphi)) =_{\sf ess} {\sf Im}(\varphi)$, and similarly
for ${\sf Dom}$.  
 \ \ \ $\Box$

\begin{lem} \label{ess_incl_morph} \ 
Let $\varphi, \psi \in {\sf riHom}(A^*)$ and
assume that \ ${\sf Im}(\psi) \subseteq_{\sf end} {\sf Im}(\varphi)$. 
Then there exist $\varphi_0, \psi_0 \in {\sf riHom}(A^*)$ such that 

\smallskip

$\bullet$ \ $\varphi_0$ is an essentially equal restriction of $\varphi$, 
 and $\psi_0$ is an essentially equal restriction of $\psi$, 

\smallskip

$\bullet$ \ ${\sf imC}(\psi_0)$ and ${\sf imC}(\varphi_0)$ are prefix codes,

\smallskip

$\bullet$ \ we have the inclusion (in the ordinary sense): 
 \ \ ${\sf imC}(\psi_0) \ \subseteq \ {\sf imC}(\varphi_0).$
\end{lem}
{\bf Proof.} By essentially restricting $\varphi$ and $\psi$, if necessary,
we can assume that ${\sf imC}(\varphi)$ and ${\sf imC}(\psi)$ are prefix
codes. Let $P = {\sf imC}(\varphi)$ and $Q = {\sf imC}(\psi)$, so
${\sf Im}(\psi) = QA^*$ and ${\sf Im}(\varphi) = PA^*$. 
By Lemma \ref{idealInter} (Lemma 3.3 of \cite{BiThomps}), there exists a 
prefix code $Q_0$ such that \ $PA^* \cap QA^* = Q_0A^*$ \ (hence, 
obviously, $Q_0A^* \subseteq QA^*$), and \ $Q_0 \subseteq P \cup Q$.
Moreover, since $\subseteq_{\sf end}$ is a lattice pre-order,  
${\sf Im}(\psi) \subseteq_{\sf end} {\sf Im}(\varphi)$ implies that 
${\sf Im}(\psi) \cap {\sf Im}(\varphi) =_{\sf ess} {\sf Im}(\varphi)$, so 
$Q_0 A^* =_{\sf ess} Q A^*$.

We now restrict $\psi$ (whose image is ${\sf Im}(\psi) = QA^*$) so as to 
obtain $\psi_0$ with ${\sf Im}(\psi_0) = Q_0A^*$ $(\subseteq QA^*)$. 
Since $QA^* =_{\sf ess} Q_0A^* \subseteq QA^*$, $\psi_0$ is an essentially 
equal restriction of $\psi$. 

Next, we partition $P$ into
$P_1 = \{ p \in P : p$ is a prefix of an element of $Q_0 \}$, and
$P_2 = P - P_1$. Then we define $P_0 = Q_0 \cup P_2$.  
Since $Q_0A^* \subseteq PA^*$, every element $p$ of $P$ is either a
prefix of an element of $Q_0$ or prefix-incomparable with all of $Q_0A^*$.
It follows that $P_0A^* \subseteq PA^*$ and $P_0A^* =_{\sf ess} PA^*$. 

Finally, we restrict $\varphi$ (whose image is ${\sf Im}(\varphi) = PA^*$)
so as to obtain $\varphi_0$ with ${\sf Im}(\varphi_0) = P_0A^*$ 
$(\subseteq PA^*)$. Since $P_0A^* =_{\sf ess} PA^*$, $\varphi_0$ is an
essentially equal restriction of $\varphi$. 

Now we have, \ \ ${\sf imC}(\psi_0) = Q_0 \subseteq Q_0 \cup P_2 = P_0 = $
 ${\sf imC}(\varphi_0)$.  \ \ \ \ \ $\Box$

\bigskip


\noindent {\bf Proof of Theorem \ref{R_order}.}  \\  
 Because of Lemma \ref{equiv_R_order} we only need to prove that 
$\psi \leq_{\cal R} \varphi$ iff 
${\sf Im}(\psi) \subseteq_{\sf end} {\sf Im}(\varphi)$. 

\smallskip

\noindent [$\Rightarrow$] 
 \ Suppose $\psi(.) = \max(\varphi \circ \alpha(.))$,
for some $\alpha \in M_{k,1}$. It is easy to see that 
${\sf Im}(\varphi \circ \alpha(.)) \subseteq {\sf Im}(\varphi(.))$   
It follows now from Lemma \ref{ess_restr} that 
${\sf Im}(\varphi \circ \alpha(.)) \subseteq_{\sf end} $
${\sf Im}(\max(\varphi \circ \alpha(.)))$.

\smallskip

\noindent
[$\Leftarrow$] \ If ${\sf Im}(\psi) \subseteq_{\sf end} {\sf Im}(\varphi)$
we can apply Lemma \ref{ess_incl_morph} and choose representations 
$\psi_0, \varphi_0$ of $\psi$, respectively $\varphi$, such that 
 \ ${\sf imC}(\psi_0) \subseteq {\sf imC}(\varphi_0)$. 
We will now define $\alpha \in {\sf riHom}(A^*)$ so that 
 \ $\varphi_0 \circ \alpha(.) = \psi_0(.) : $

\medskip

$\bullet$ We pick \ ${\sf domC}(\alpha) = {\sf domC}(\psi_0)$; 

\medskip

$\bullet$ for every $y \in {\sf imC}(\psi_0)$ we choose an element 
 \ $\overline{y} \in \varphi_0^{-1}(y)$;

\medskip

$\bullet$ for each $x \in {\sf domC}(\psi_0)$ we define 
 \ $\alpha(x) = \overline{\psi_0(x)}$ \ 
 $\big( \, \in \varphi_0^{-1}(\psi_0(x)) \big)$.

\medskip

\noindent Now for every $x \in {\sf domC}(\psi_0)$ we have: \ \  
$\varphi_0 \circ \alpha(x) \ = \ $
$\varphi_0 \big( \, \overline{\psi_0(x)} \, \big)$
$ \ \in \ \varphi_0 \big(\varphi_0^{-1}(\psi_0(x))\big) \ = \ \psi_0(x)$. 
So, $\varphi_0 \circ \alpha(.) = \psi_0(.)$.
 \ \ \ $\Box$


\section{The ${\cal L}$-order of $M_{k,1}$  }

Just as for the ${\cal R}$-order, the ${\cal L}$-order of the monoid 
$M_{k,1}$ has some similarities with the ${\cal L}$-order of the monoid 
${\sf PF}_X$ of all partial functions on a set $X$. In both cases, the 
${\cal L}$-order between functions is related to the refinement order of 
the partitions on the domains of the functions. However, for $M_{k,1}$ we 
need more complicated notions of partition and of refinement.

In the following subsections we first define right congruences in $A^*$ and 
essential equality of right congruences. We associate a right congruence 
with every element of $M_{k,1}$. 
We define the refinement order of right congruences, and 
finally we use that to characterize the ${\cal L}$-order of $M_{k,1}$. 


\subsection{Right congruences, prefix code partitions, 
  and essential equivalence}

A {\em right congruence on} a right ideal $R \subseteq A^*$ is 
an equivalence relation $\simeq$ on $R$ such that for all $x,y \in R$ and 
all $w \in A^*:$ \ $x \simeq y$ implies $xw \simeq yw$; moreover, 
$x \simeq y$ is undefined if $x$ or $y$ are not both in $R$. 
The right ideal $R$ is called the {\em domain} of $\simeq$ and we denote it
by ${\sf Dom}(\simeq)$.  
We will only consider the case when $R$ is {\it finitely generated} as a 
right ideal. 
The equivalence class containing an element $x$ will be denoted by $[x]$.

Let $P \subset A^*$ be a finite prefix code, and let $=_P$ be an equivalence 
relation on $P$. Then $=_P$ determines a right congruence $\simeq_P$ on the 
right ideal $PA^*$, as follows: 
If $p_1, p_2 \in P$ and $p_1 =_P p_2$ then $p_1w \simeq_P p_2w$ for all 
$w \in A^*$; and if $p_1x, p_2y \in PA^*$ are such that 
$p_1 \neq_P p_2$ or $x \neq y$, then $p_1x \not\simeq_P p_2y$.
Thus, if the set of equivalence classes of $=_P$ is \ $\{ P_1, \ldots, P_n\}$,
then the set of congruence classes of $\simeq_P$ is
 \ $\{ P_j w \ : \ w \in A^*, \ 1 \leq j \leq n \}$.
Hence $\simeq_P$ is the {\it coarsest} right congruence on $PA^*$ 
that agrees with $=_P$ on $P$. We also say that $\simeq_P$ is the 
{\em right congruence generated by the equivalence relation} $=_P$.

\begin{defn} \label{prefix_congruence} \ Let $P \subset A^*$ be a finite 
prefix code, let $\simeq$ be a right congruence on the right ideal $PA^*$,
and let $=_P$ be the restriction of $\simeq$ to $P$.
We call $\simeq$ a {\bf prefix code congruence} \ iff \ $\simeq$ is equal 
to the right congruence generated by its restriction $=_P$. In that case
we call $P$ the {\bf domain code} of $\simeq$, and denote it by
${\sf domC}(\simeq)$.
\end{defn}
Not every right congruence on $PA^*$ is a prefix code congruence.
For example, $\simeq$ is not determined by its restriction to $P$ if 
$pu \simeq puv$ (for some $p \in P$ and some $u,v \in A^*$ with $v$ 
non-empty), or if $p_1x \simeq p_2y$ (for some $p_1 \not\simeq p_2 \in P$ 
and some $x, y \in A^*$).  
In this paper we will only consider prefix code congruences.

\medskip

A prefix code congruence $\simeq$ can be extended to the {\em ends} of
${\sf Dom}(\simeq)$.
For $w_1, w_2 \in {\sf ends}({\sf Dom}(\simeq))$, we say that \ 
$w_1 \simeq w_2$ \ iff \ there exist \ $p_1, p_2 \in {\sf Dom}(\simeq)$ 
 \ and $v \in A^{\omega}$ such that 

\smallskip

$w_1 = p_1v$, \ $w_2 = p_2v$, \ \ and \ \ $p_1 \simeq p_2$ .   

\smallskip

\noindent Hence, the set of right congruence classes of $\simeq$ in
${\sf ends}({\sf Dom}(\simeq))$ is 
 \  $\{ [p] \, v \ : \ p \in {\sf Dom}(\simeq), \ v \in A^{\omega} \}$.

\smallskip

\noindent Notational Remark: Although $\simeq$ can be extended to a 
partition of ${\sf ends}({\sf Dom}(\simeq))$, our notation 
${\sf Dom}(\simeq)$ will continue to refer to the right ideal of finite 
words on which $\simeq$ is defined; i.e., in our notation we still have 
${\sf Dom}(\simeq) \subseteq A^*$.

\medskip

For a prefix code congruence $\simeq$, the prefix code ${\sf domC}(\simeq)$ 
is finite, by definition. It follows that for a prefix code congruence we 
have: {\em Every $\simeq$-class in ${\sf ends}({\sf Dom}(\simeq))$ or in 
${\sf Dom}(\simeq)$ is finite, with cardinality uniformly bounded from above 
by \ $|{\sf domC}(\simeq)|$. }

\begin{defn} \label{congr_ext_restr} \    
Let $\simeq_1$ and $\simeq_2$ be two prefix code congruences. We say that
$\simeq_1$ is an {\bf essentially equal extension} of $\simeq_2$ (and that 
$\simeq_2$ is an {\bf essentially equal restriction} of $\simeq_1$) iff the 
following three conditions hold: 

\smallskip

\noindent
{\bf (1)} \ \ \ ${\sf Dom}(\simeq_2) \ =_{\sf ess} \ {\sf Dom}(\simeq_1)$,
 \ and 

\smallskip

\noindent
{\bf (2)} \ \ \ ${\sf Dom}(\simeq_2) \ \subseteq \ {\sf Dom}(\simeq_1)$,  
  \ and

\smallskip

\noindent
{\bf (3)} \ \ \ $\simeq_2$ agrees with $\simeq_1$ on ${\sf Dom}(\simeq_2)$;
 \ i.e., for all $x, y \in {\sf Dom}(\simeq_2): \ $
$ x \simeq_2 y \Leftrightarrow x \simeq_1 y$.
\end{defn}
Conditions (1) and (2) are equivalent to saying that every element of 
${\sf Dom}(\simeq_1)$ is the prefix of some element of 
${\sf Dom}(\simeq_2)$, and that element of ${\sf Dom}(\simeq_2)$ has a 
prefix in ${\sf Dom}(\simeq_1)$.

Conditions (2) and (3) are equivalent to saying that every $\simeq_2$-class 
is also a $\simeq_1$-class; (2) and (3) are also equivalent to saying 
that \ ${\sf Dom}(\simeq_2) \ \subseteq \ {\sf Dom}(\simeq_1)$ \ and for 
every $x \in {\sf Dom}(\simeq_2), \ [x]_1 = [x]_2$. Here, $[x]_1$ and 
$[x]_2$ denote the equivalence class of $x$ for $\simeq_1$, respectively
$\simeq_2$. 

\medskip

Essentially equal restrictions and extensions of prefix code congruences 
can be determined by a replacement (or rewriting) process, based on the 
following {\bf replacement rules} (where $A = \{a_1, \ldots, a_k\}$):

\begin{equat}\hspace{-.06in}{\bf )} \label{congr_replacement} 
\hspace{.4in}
Replace the class $C$ in the domain code $P$ by the set of 
classes \ $\{C a_1, \ldots, C a_k\}$.
\end{equat}
\begin{equat}\hspace{-.08in}{\bf )} \label{congr_InvReplacement}
\hspace{.4in} Replace the set of classes 
 \ $\{C a_1, \ldots, C a_k\}$ \ of the domain code $P$ by the new class $C$.
\end{equat}
When rule (\ref{congr_replacement}) is applied, the domain code $P$ is
replaced by \, $(P - C) \cup CA$; as a result, an essentially equal 
restriction of the prefix code congruence is obtained.
Similarly, when rule (\ref{congr_InvReplacement}) is applied the domain
code $P$ is replaced by \, $(P - CA) \cup C$; as a result, an essentially 
equal extension of the prefix code congruence is obtained. 
The replacement steps (\ref{congr_replacement}) and
(\ref{congr_InvReplacement}) can be iterated. It turns out that all 
essentially equal restrictions and extensions can be obtained in the above 
way:

\begin{pro} \label{congr_ext_restr_Rewriting} \   
Let $\simeq_1$ and $\simeq_2$ be two prefix code congruences. Then 
$\simeq_2$ is an essentially equal restriction of $\simeq_1$ \ iff 
 \ $\simeq_2$ can be obtained from $\simeq_1$ by a finite sequence of 
replacements of the form {\rm (\ref{congr_replacement})}. And $\simeq_1$ is an 
essentially equal extension of $\simeq_2$ \ iff \ $\simeq_1$ can be obtained 
from $\simeq_2$ by a finite sequence of replacements of the form
{\rm (\ref{congr_InvReplacement})}.
\end{pro}
{\bf Proof.} The proof is similar to the proof of Prop.\ 1.4 in
\cite{BiThomMon} (as well as the proof of Lemma 2.2 in \cite{BiThomps},
going back to Thompson).
The direction $[\Leftarrow]$ is easy to see. 

Conversely, suppose that $\simeq_2$ is an essentially equal restriction of 
$\simeq_1$. Let $P_1 = {\sf domC}(\simeq_1)$ and 
$P_2 = {\sf domC}(\simeq_2)$.
Since ${\sf Dom}(\simeq_1)$ and ${\sf Dom}(\simeq_2)$
are essentially equal and ${\sf Dom}(\simeq_2) \subseteq {\sf Dom}(\simeq_1)$, 
every path (in the tree of $A^*$) starting in $P_1$ reaches $P_2$.  Hence, 
the set difference \ ${\sf Dom}(\simeq_1) - {\sf Dom}(\simeq_2)$ \ is finite.
Also, since $\simeq_2$ agrees with $\simeq_1$ on ${\sf Dom}(\simeq_2)$ it 
follows that \ ${\sf Dom}(\simeq_1) - {\sf Dom}(\simeq_2)$ \ consists of 
$\simeq_1$-equivalence classes. 

Let $\pi$ be a $\simeq_1$-equivalence class that lies in 
 \ ${\sf domC}(\simeq_1) - {\sf domC}(\simeq_2)$; the latter set is not empty,
otherwise $\simeq_1$ and $\simeq_2$ would be equal. Removing $\pi$ from 
$\simeq_1$ yields a new prefix code congruence $\simeq_1'$ with domain 
 \ ${\sf Dom}(\simeq_1') = {\sf Dom}(\simeq_1) - \pi$. 
In the tree of $A^*$, the children of the elements of $\pi$ are
 \ $\bigcup_{i=i}^k \pi a_i$. So, ${\sf domC}(\simeq_1') = $ 
$(P_1 - \pi) \ \cup \ \bigcup_{i=i}^k \pi a_i$. This amounts to applying 
rule of type (\ref{congr_replacement}) to $\simeq_1$. 

Since \ $|{\sf Dom}(\simeq_1') - {\sf Dom}(\simeq_2)| \ < \ $
$|{\sf Dom}(\simeq_1) - {\sf Dom}(\simeq_2)|$, we conclude by induction that
$\simeq_2$ can be obtained from $\simeq_1'$ by rules of type
(\ref{congr_replacement}). Hence, the essentially equal restrictions from
$\simeq_1$ to $\simeq_1'$ and from there to $\simeq_2$ can be carried out by 
applying rules of type (\ref{congr_replacement}).  
 
In a similar way one proves that essentially equal extensions can be carried 
out by rules of type (\ref{congr_InvReplacement}). 
 \ \ \ $\Box$

\begin{defn} \label{ess_equ_congr} \
Two prefix code congruences $\simeq_1$ and $\simeq_2$ are {\bf essentially
equal} \ iff \ $\simeq_1$ and $\simeq_2$ can be obtained from each other by a
finite sequence of essentially equal extensions and essentially equal 
restrictions. We denote this by \ $\simeq_1 \ =_{\sf ess} \ \simeq_2$ .
\end{defn}
By Proposition \ref{congr_ext_restr_Rewriting}, this could be defined 
equivalently by: \ {\em $\simeq_1$ and $\simeq_2$ are essentially equal iff
$\simeq_1$ and $\simeq_2$ can be obtained from each other by a
finite sequence of replacement steps of form (\ref{congr_replacement}) and
(\ref{congr_InvReplacement}).} We also have the following characterization.

\begin{pro} \label{ess_equ_congrCharact} \   
Two prefix code congruences $\simeq_1$ and $\simeq_2$ are {\bf essentially
equal} \ iff \ 

\smallskip

\noindent {\bf (1)} \ \ \ ${\sf Dom}(\simeq_2) \ =_{\sf ess} \ {\sf
Dom}(\simeq_1)$,
 \ and

\smallskip

\noindent
{\bf (2)} \ \ \ $\simeq_2$ agrees with $\simeq_1$ on 
 \ ${\sf Dom}(\simeq_1) \cap {\sf Dom}(\simeq_2)$.
\end{pro}
{\bf Proof.} By Definition \ref{congr_ext_restr},  if $\simeq_1$ is an 
essentially equal extension of $\simeq_2$ then $\simeq_1$ and $\simeq_2$
agree on ${\sf Dom}(\simeq_1) \cap {\sf Dom}(\simeq_2)$. For finite 
sequence $\simeq_1$, $\ldots$, $\simeq_i$, $\ldots$, $\simeq_n$ of prefix 
code congruences, successively obtained from each other by essentially equal
extensions and restrictions, $\simeq_1$ and $\simeq_n$ will agree on
$\bigcap_{i=1}^n {\sf Dom}(\simeq_i)$. But then, if we extend (by applying 
replacement steps (\ref{congr_InvReplacement})), $\simeq_1$ and $\simeq_n$ 
will agree on ${\sf Dom}(\simeq_1) \cap {\sf Dom}(\simeq_2)$.

Conversely, suppose $\simeq_1$ and $\simeq_2$ agree on
${\sf Dom}(\simeq_1) \cap {\sf Dom}(\simeq_2)$. Then the restriction $\simeq$
of $\simeq_1$ to ${\sf Dom}(\simeq_1) \cap {\sf Dom}(\simeq_2)$ is an 
essentially equal restriction of $\simeq_1$. Moreover, the extension of 
$\simeq$ from ${\sf Dom}(\simeq_1) \cap {\sf Dom}(\simeq_2)$ to 
${\sf Dom}(\simeq_2)$ is an essentially equal extension of $\simeq$. So 
$\simeq_1$ and $\simeq_2$ are essentially equal.
 \ \ \ $\Box$

\medskip

The replacement system consisting of the rules of type
(\ref{congr_InvReplacement}) is {\em terminating} (since the number of
classes in the finite prefix code decreases at each step) and 
{\em confluent} (there are no overlaps between left-sides of rules). Hence, 
there is a unique result for the iterated replacement in the direction 
(\ref{congr_InvReplacement}). So we proved:

\begin{pro} \label{congr_Unique_max} \  
Every prefix code congruence has a {\bf unique} maximal essentially equal 
extension.  \ \  \ \ \ $\Box$
\end{pro}
This maximal essentially equal extension of $\simeq$ is denoted by 
${\sf max}(\simeq)$. 

\begin{pro} \label{congr_ess_simil_characteriz} \ 
For prefix code congruences $\simeq_1$ and $\simeq_2$ the following are 
equivalent:

\smallskip

\noindent {\bf (1)} \ \ \ $\simeq_1 \ =_{\sf ess} \ \simeq_2$ ;

\smallskip

\noindent
{\bf (2)} \ \ \ ${\sf max}(\simeq_1) \ = \ {\sf max}(\simeq_2)$ ;

\smallskip

\noindent
{\bf (3)} \ \ \ ${\sf ends}({\sf Dom}(\simeq_1)) \ = \ {\sf ends}({\sf
Dom}(\simeq_2))$, \ and
$\simeq_2$ agrees with $\simeq_1$ on \ ${\sf ends}({\sf Dom}(\simeq_1))$.
\end{pro}
{\bf Proof.} $[(1) \Rightarrow (2)]$ \ If $\simeq_1$ and $\simeq_2$ can be
obtained from each other by rewriting according to (\ref{congr_replacement})
and (\ref{congr_InvReplacement}), both can be rewritten to
${\sf max}(\simeq_1)$, as well as to ${\sf max}(\simeq_2)$. By uniqueness 
of the maximal essentially equal extension, we obtain (2).

\smallskip

\noindent
$[(2) \Leftarrow (1)]$ \ If $\simeq_1$ and $\simeq_2$ have the same maximal
essentially equal extension, we can rewrite $\simeq_1$ to $\simeq_2$ via 
this common maximal essentially equal extension. So, 
$\simeq_1 \ =_{\sf ess} \ \simeq_2$.

\smallskip

\noindent
$[(2) \Leftrightarrow (3)]$ \ The set of ends ${\sf ends}({\sf Dom}(\simeq_1))$
uniquely determines the prefix code ${\sf domC}({\sf max}(\simeq_1))$.
Namely, we take the set shortest prefixes of ends in 
${\sf ends}({\sf Dom}(\simeq_1))$ that are not prefixes of ends that are not
in ${\sf ends}({\sf Dom}(\simeq_1))$. Thus, we can write each end $w$ in
${\sf ends}({\sf Dom}(\simeq_1))$ uniquely as $w = pv$ with 
$p \in {\sf domC}({\sf max}(\simeq_1))$, and $v \in A^{\omega}$.
The partition $\simeq_1$ on ${\sf ends}({\sf Dom}(\simeq_1))$ then uniquely 
determines ${\sf max}(\simeq_1)$ on ${\sf domC}({\sf max}(\simeq_1))$.

Conversely, ${\sf max}(\simeq_1)$ on ${\sf domC}({\sf max}(\simeq_1))$
determines the partition $\simeq_1$ on ${\sf ends}({\sf Dom}(\simeq_1))$.
Since ${\sf max}(\simeq_1)$ determines $\simeq_1$ on 
${\sf ends}({\sf Dom}(\simeq_1))$, and vice versa, it follows that (2) is 
equivalent to (3).
 \ \ \ $\Box$


\subsection{The prefix code congruence of a right ideal homomorphism} 

It is well known that with any partial function $f: X \to Y$ one can associate 
an equivalence relation $\equiv_f$ on $X$, defined by $x_1 \equiv_f x_2$ iff 
$f(x_1) = f(x_2)$; the set of equivalence classes is
 \ $\{ f^{-1}(y) : y \in {\sf Im}(f)\}$. When $X$ and $Y$ have a structure 
and $f$ is a homomorphism for that structure, then $\equiv_f$ is a congruence 
for that structure. 

\begin{defn} \label{partition_riHom}
{\bf (The partition ${\sf part}(\varphi)$).} \
For $\varphi \in {\sf riHom}(A^*)$ we consider the right congruence on the
right ideal ${\sf Dom}(\varphi)$, defined by \ $x_1 \equiv x_2$ \ iff 
 \ $\varphi(x_1) = \varphi(x_2)$. 
This right congruence is called ${\sf part}(\varphi)$.
\end{defn}
As the following example shows, ${\sf part}(\varphi)$ is not always a 
{\it prefix code congruence} (according to 
Definition \ref{prefix_congruence}).
Let $A = \{a,b\}$, and let $\varphi$ be given by the table

\medskip

$\left[ \begin{array}{cc}
   a   & b   \\
   aa  & a
\end{array} \right] $, \ \  which has an essentially equal restriction 
$\varphi'$ with table \ 
$\left[ \begin{array}{ccc}
   a  & ba & bb   \\   
   aa & aa & ab
\end{array} \right] $ .

\medskip

\noindent The set of classes of ${\sf part}(\varphi)$ is 
 \ $\big\{ \{b\} \big\}$ $ \, \cup \, $
$\big\{ \{av, bav\} : v \in \{a,b\}^* \big\}$
$ \, \cup \, $  $\big\{ \{bbw\} : w \in \{a,b\}^*\big\}$. 
The presence of the class $\{b\}$ prevents ${\sf part}(\varphi)$ from 
being a prefix code congruence; indeed, for the class $\{b\}$ we see that 
 \ $\{b\} \, a$ \ is a strict subset of the class \ $\{a,ba\}$. 
On the other hand, ${\sf part}(\varphi')$ is a prefix code congruence. 

Interestingly, this is related to an issue that was mentioned when 
$M_{k,1}$ was first defined in \cite{BiThomMon}, namely the fact that 
$\varphi({\sf domC}(\varphi))$ {\it is not necessarily a prefix code}. 
In the above example, $\varphi({\sf domC}(\varphi)) = \{aa, a\}$.
The connection between these issues is given by the following.

\begin{pro} \label{partVS_imC1} \  
For any $\varphi \in {\sf riHom}(A^*)$ and its congruence 
${\sf part}(\varphi)$ we have: 

\smallskip

 \hspace{.4in} ${\sf part}(\varphi)$ is a prefix code congruence \ \ iff 
 \ \ $\varphi({\sf domC}(\varphi))$ is a prefix code.

\smallskip

\noindent When $\varphi({\sf domC}(\varphi))$ is a prefix code it is
also denoted by ${\sf imC}(\varphi)$ (called the {\it image code}).
\end{pro}
{\bf Proof.} \ $[ \Leftarrow ]$ \  We assume that 
$\varphi({\sf domC}(\varphi)) = {\sf imC}(\varphi)$ is a prefix code. 
Let $x_1, x_2 \in {\sf domC}(\varphi)$ and $u,v \in A^*$ 
be such that $\varphi(x_1u) = \varphi(x_2v)$.
Then $\varphi(x_1),\ \varphi(x_2) \in {\sf imC}(\varphi)$, and since this is 
a prefix code it follows that $\varphi(x_1)$ and $\varphi(x_2)$ are either 
prefix-incomparable or equal. Since $\varphi(x_1u) = \varphi(x_2v)$ it follows
that $\varphi(x_1) = \varphi(x_2)$ and that $u = v$. So,
$x_1u, x_2v \in \varphi^{-1}(yu) = \varphi^{-1}(y) \, u$ \ for some 
$y \in {\sf imC}(\varphi)$ (where $y = \varphi(x_1) = \varphi(x_2)$). 
In other words, every class of ${\sf part}(\varphi)$ is of the form 
 \ $\varphi^{-1}(y) \ u$ \ for some $y \in {\sf part}(\varphi)$, and 
$u \in A^*$. 

Moreover, $\varphi^{-1}(y) \subseteq {\sf domC}(\varphi)$ for
every $y \in {\sf imC}(\varphi) = \varphi({\sf domC}(\varphi))$.
Indeed, if $\varphi(xu) = y \in {\sf imC}(\varphi)$ for some 
$x \in {\sf domC}(\varphi)$ and $u \in A^*$, then $\varphi(x) \in $
${\sf imC}(\varphi) = \varphi({\sf domC}(\varphi))$; hence $u = \varepsilon$,
since ${\sf imC}(\varphi)$ is a prefix code, so 
$xu = x \in {\sf domC}(\varphi)$.
Hence, ${\sf part}(\varphi)$ is a prefix code congruence, with domain code
${\sf domC}(\varphi)$. 

\noindent $[ \Rightarrow ]$ \ We assume that 
$\varphi({\sf domC}(\varphi))$ is not a prefix code.
So, there exist $y, yu \in \varphi({\sf domC}(\varphi))$ with $u \in A^*$ 
and $u \neq \varepsilon$. Let $x_1, x_2 \in {\sf domC}(\varphi)$ be such 
that $\varphi(x_1) = y$ and $\varphi(x_2) = yu$. Since $y \neq yu$, 
$[x_1] \cap [x_2] = \varnothing$. On the other hand, $\varphi(x_1u) = yu$,
so $x_1 u \in [x_1] \, u \cap [x_2] \neq \varnothing$, and 
$[x_2] \neq [x_1] \, u$ (since $[x_1] \, u = [x_2] \subseteq {\sf
domC}(\varphi)$ would imply that $x_1$ and $x_1u$ both belong to the prefix
code ${\sf domC}(\varphi)$, which is impossible since $u \neq \varepsilon$).   
Hence, $[x_1] \, u$ is not a class of ${\sf part}(\varphi)$, so 
${\sf part}(\varphi)$ is not a prefix code congruence.
 \ \ \ $\Box$

\medskip

\noindent This motivates the following.

\begin{defn} \label{defn_riHom_pc} \   
Within the monoid ${\sf riHom}(A^*)$ we consider the submonoid 

\smallskip

 \ \ \ \ \  ${\sf riHom}_{\sf pc}(A^*) \ = \ $
$\{ \varphi \in {\sf riHom}(A^*) \ : \ \varphi({\sf domC}(\varphi)) \ $
${\rm is \ a \ prefix \ code} \}$.

\smallskip

\noindent The elements of ${\sf riHom}_{\sf pc}(A^*)$ are called 
{\bf prefix code preserving}.
\end{defn}
The subscript ``{\sf pc}'' stands for ``prefix code''.
It is easy to check that ${\sf riHom}_{\sf pc}(A^*)$ is indeed a monoid.
The reason for calling the elements of ${\sf riHom}_{\sf pc}(A^*)$ ``prefix 
code preserving'' is the following.

\begin{pro} \label{prefcodePres} \   
For all $\varphi \in {\sf riHom}(A^*)$ we have: 
 \ $\varphi({\sf domC}(\varphi))$ is a prefix code \ iff \ for every prefix 
code $P \subset A^*$, $\varphi(P)$ is a prefix code.
\end{pro}
{\bf Proof.} The right-to-left implication is trivial. We prove the 
left-to-right implication by contraposition. Let $x_1, x_2 \in A^*$ prefix 
incomparable, but assume by contradiction that 
$\varphi(x_2) = \varphi(x_1) \, w$, for some non-empty $w \in A^*$. 
Assume also that $x_1, x_2 \in {\sf Dom}(\varphi)$, so there are 
$p_1, p_2 \in {\sf domC}(\varphi)$ such that $x_1 = p_1u_1$, 
$x_2 = p_2u_2$ (for some $u_1, u_2 \in A^*$). Then 
$\varphi(x_2) = \varphi(x_1) \, w$ implies 
 \ $\varphi(p_2) = \varphi(p_1) \, u_1w$, which implies that $\varphi(p_2)$ 
and $\varphi(p_1)$ are prefix comparable.   \ \ \ $\Box$

\smallskip

\noindent The following further demonstrates the importance of the monoid 
${\sf riHom}_{\sf pc}(A^*)$.
\begin{pro} \
Every $\varphi \in {\sf riHom}(A^*)$ has an essential restriction to some 
element of ${\sf riHom}_{\sf pc}(A^*)$.
\end{pro}
{\bf Proof.} It is easy to restrict $\varphi$ (to some element $\Phi$) such
that the image code becomes ${\sf imC}(\Phi) = A^{\ell}$, where $\ell$ is
the length of a longest word in $\varphi({\sf domC}(\varphi))$. Obviously,
$A^{\ell}$ is a prefix code.
 \ \ \ $\Box$

\medskip

Henceforth, when we use ${\sf part}(\varphi)$ we will always assume that 
${\sf part}(\varphi)$ is a prefix code congruence; equivalently, we assume
that $\varphi({\sf domC}(\varphi)) = {\sf imC}(\varphi)$ is a prefix code.

\medskip

A related issue is the fact that essentially equal restrictions and
extensions of prefix code congruences are defined in a more limited way
than restrictions and extensions of elements of ${\sf riHom}(A^*)$: In an
essentially equal restriction of a prefix code congruence, an entire class
$C$ is replaced by the set of classes $\{Ca_1, \ldots, Ca_k\}$. On the other
hand, in an essentially equal restriction of $\varphi \in {\sf riHom}(A^*)$,
a single element $x \in {\sf domC}(\varphi)$ is replaced by
$\{xa_1, \ldots, xa_k\}$ (with accompanying replacements of the image
$\varphi(x)$ by $\{\varphi(x) \, a_1, \ \ldots, \ \varphi(x) \, a_k\}$).
So, besides the general essentially equal restrictions of
$\varphi \in {\sf riHom}(A^*)$ we will consider the following important
special case:

\begin{defn} \label{congruence_restr_ext} \
Let $\varphi, \Phi \in {\sf riHom}(A^*)$. Then $\Phi$ is an {\bf essentially
equal class-wise restriction} of $\varphi$ (and $\varphi$ is an
{\bf essentially equal class-wise extension} of $\Phi$) \ iff \ $\Phi$ is an
essentially equal restriction of $\varphi$ such that
${\sf Dom}(\varphi) - {\sf Dom}(\Phi)$ is a union of classes of
${\sf part}(\varphi)$.
\end{defn}
The best way to understand class-wise restrictions (or extensions) is to
think of them as {\em restrictions (or extensions) from the point of view 
of ${\sf imC}(\varphi)$}. More precisely, to create such a restriction we 
take $y \in {\sf imC}(\varphi)$, replace ${\sf imC}(\varphi)$ by 
$({\sf imC}(\varphi) - \{y\}) \cup \{ya_1, \ldots, ya_k\}$, and then
do the corresponding replacement in ${\sf domC}(\varphi)$ (as described in
Definition \ref{congruence_restr_ext}). 

It is easy to see that $\Phi$ is an {\it essentially equal class-wise
restriction} of $\varphi$ iff $\Phi$ can be obtained from $\varphi$ by a
finite number of the following type of replacement steps.
Below, $\chi \in {\sf riHom}(A^*)$ is any intermediate element obtained.

\begin{equat}\hspace{-.08in}{\bf )} \label{classwiseReplacement}
\hspace{.4in}
Replace \ $\{(x, \, \chi(C)) : \ x \in C\}$ \ in the table by
 \ \ $\{(xa_i, \, \chi(C) a_i) : \ x \in C, \ 1 \leq i \leq k \}$,
\end{equat}

\noindent where $C$ is a class of ${\sf part}(\chi)$ in ${\sf domC}(\chi)$.
An {\it essentially equal class-wise extension} is obtained by repeated
replacements in the opposite direction, i.e., of the form

\begin{equat}\hspace{-.08in}{\bf )} \label{classwiseReplacement_inv}
\hspace{.4in}
Replace \ $\{(xa_i, \, \chi(C) \, a_i) : \ x \in C, \ 1 \leq i \leq k\}$
 \ in the table by \ $\{(x, \, \chi(C)) : \ x \in C \}$,
\end{equat}

\noindent where $Ca_1, \ldots, Ca_k$ are classes of ${\sf part}(\chi)$ in
${\sf domC}(\chi)$.

\medskip

The significance of these replacement rules is demonstrated by the following.
\begin{pro} \label{partVS_imC2} \   
For any $\varphi \in {\sf riHom}_{\sf pc}(A^*)$ we have: \   
An essentially equal restriction or extension of $\varphi$ leads again to
an element of ${\sf riHom}_{\sf pc}(A^*)$
 \ iff \ this restriction or extension is an essentially equal
{\em class-wise} restriction or extension.
\end{pro}
{\bf Proof.} Consider an essentially equal restriction of $\varphi$
in which a class $C$ of ${\sf part}(\varphi)$ (contained in
${\sf domC}(\varphi)$) is replaced in part.
In other words, there are $x_1, x_2 \in C$ such that
$(x_1,y)$ is left unchanged (where  $\varphi(C) = y \in {\sf imC}(\varphi)$),
and $(x_2,y)$ is replaced by $\{ (x_2a_1, ya_1), \ldots, (x_2a_k, ya_k) \}$.
Then for the resulting element $\Phi \in {\sf riHom}(A^*)$ obtained from
$\varphi$ we have \ $y, y_a1, \ldots, ya_k \in \Phi({\sf domC}(\Phi))$,
hence $\Phi({\sf domC}(\Phi))$ is not a prefix code.

On the other hand, assume only entire classes are replaced; e.g., for a class
$C$ of ${\sf part}(\varphi)$ in ${\sf domC}(\varphi)$, we replace all of
$\{(x,y) : x \in C\}$ \ by
 \ $\{(xa_i, ya_i) : x \in C, \ 1 \leq i \leq k\}$, where $ y  = \varphi(C)$.
Then ${\sf imC}(\varphi)$ is replaced by
 \ $\big({\sf imC}(\varphi) - \{y\} \big) \cup yA$; this is a prefix code
if ${\sf imC}(\varphi)$ is a prefix code.

For essentially equal extensions, the proof is similar.
 \ \ \ $\Box$

\medskip

By Prop.\ \ref{partVS_imC2}, ${\sf riHom}_{\sf pc}(A^*)$ is closed under
essential class-wise restriction and extension, as introduced in Definition
\ref{congruence_restr_ext}.

\medskip

The replacement rules (\ref{classwiseReplacement}) and 
(\ref{classwiseReplacement_inv}). Are easily seen to form a confluent and 
terminating rewriting system, in the direction 
(\ref{classwiseReplacement_inv}).  Hence, for each 
$\varphi \in {\sf riHom}(A^*)$ there exists a unique maximal element
${\sf max}_{\sf pc}(\varphi)$ for the rules (\ref{classwiseReplacement_inv}).
 
\begin{pro} \label{equiv_vs_classrewrites} \hspace{-.06in}{\bf .}

\noindent {\bf (1)} \ For all 
$\varphi_1, \varphi_2 \in {\sf riHom}_{\sf pc}(A^*)$ we have:
 \ $\varphi_1 = \varphi_2$ in $M_{k,1}$ \ iff \ $\varphi_1$ and $\varphi_2$ 
can be obtained from each other by a finite number of applications of the 
replacement rules
(\ref{classwiseReplacement}) and (\ref{classwiseReplacement_inv}).

\noindent {\bf (2)} \ For all 
$\varphi_1, \varphi_2 \in {\sf riHom}_{\sf pc}(A^*)$ we have:
 \ $\varphi_1 = \varphi_2$ in $M_{k,1}$ \ iff 
 \ ${\sf max}_{\sf pc}(\varphi_1) = {\sf max}_{\sf pc}(\varphi_2)$.

\noindent {\bf (3)} \ For all $\varphi \in {\sf riHom}_{\sf pc}(A^*)$, 
 \ ${\sf max}_{\sf pc}(\varphi)$ is the maximum class-wise extension 
of $\varphi$.
\end{pro}
{\bf Proof.} For (1), the implication $[ \Leftarrow ]$ is obvious. For 
$[ \Rightarrow ]$ we consider a common essentially equal class-wise 
restriction $\varphi_0$ of both $\varphi_1$ and $\varphi_2$ (which exists 
since $\varphi_1$ and $\varphi_2$ are equal as elements of $M_{k,1}$). Next, 
we can extend $\varphi_0$ to $\varphi_1$ and to $\varphi_2$ by essentially 
equal class-wise extension steps. 

The proofs of (2) and (3) are straightforward.    \ \ \ $\Box$
 
\begin{pro} \label{part_invariant} \ 
If $\varphi_1, \varphi_2 \in {\sf riHom}_{\sf pc}(A^*)$ represent the same
element of $M_{k,1}$ then \ \
${\sf part}(\varphi_1) \ =_{\sf ess} \ {\sf part}(\varphi_2)$.
\end{pro}
{\bf Proof.}  By Propositions \ref{congr_ext_restr_Rewriting} and by Prop.\
\ref{congr_Unique_max} it is enough to prove that if $\varphi_2$ is obtained 
from $\varphi_1$ by one essential congruence extension (or restriction) step, 
then ${\sf part}(\varphi_2)$ is obtained from ${\sf part}(\varphi_1)$ by 
essential extension (or restriction) steps. We only consider the case of an
extension step, the case of a restriction step being similar.
Suppose $Ca_1, \ldots, Ca_k$ are classes of ${\sf part}(\varphi_1)$ in 
${\sf domC}(\varphi_1)$, and suppose 
 \ $\varphi_1(Ca_1) = ya_1, \ \ldots, \ \varphi_1(Ca_k) = ya_k$ \ for some 
$y \in A^*$. Let $\varphi_2$ be obtained by extending $\varphi_1$ by 
$\varphi_2(C) = y$. 
Then \ $\varphi_2^{-1}(y) \ a_i  \ = \ \varphi_2^{-1}(ya_i) \ = \ $
$\varphi_1^{-1}(ya_i) \in {\sf part}(\varphi_1)$, for $i = 1, \ldots, k$.
Rule (\ref{congr_InvReplacement}) can be applied to this situation; this 
leads to a new prefix code congruence, obtained by adding 
$\varphi_2^{-1}(y) \ (= C)$ to ${\sf part}(\varphi_1)$. 
But this new prefix code congruence is precisely ${\sf part}(\varphi_2)$,
since $\varphi_2$ is obtained from $\varphi_1$ by adding $\varphi_2(C) = y$.
So, ${\sf part}(\varphi_2)$ is obtained from ${\sf part}(\varphi_1)$ by one
extension step.
  \ \ \ $\Box$

\medskip

The converse of Proposition  \ref{part_invariant} is obviously 
not true. E.g., for every $\varphi \in G_{k,1}$ the essential congruence
${\sf max}({\sf part}(\varphi))$ is the same, namely the congruence given 
by the prefix code partition $\{ \{ \varepsilon \} \}$ (consisting of a 
single class, where $\varepsilon$ is the empty word); the prefix code 
congruence that 
corresponds to this is the discrete partition of $A^*$ (with singletons
as classes). This example also gives an instance where $\varphi$ is maximally 
extended, whereas ${\sf part}(\varphi)$ is not maximally extended (neither 
class-wise nor in the general sense). 

\medskip

We will show in Prop.\ \ref{func_vs_part} that every prefix code congruence 
is the prefix code congruence of some right-ideal homomorphism. 
First we will need a characterization of the idempotents of 
${\sf riHom}(A^*)$ and $M_{k,1}$.

\begin{lem} \label{idempot_characteriz_1} \hspace{-.08in}{\bf .}

\noindent 
{\bf (1)} \ An element $\eta \in {\sf riHom}(A^*)$ is an idempotent (for 
the operation of composition) \ iff \ for every $y \in {\sf Im}(\eta):$ 
 \ $y \in \eta^{-1}(y)$. \ The latter is equivalent to $\eta(y) = y$ for all
$y \in {\sf Im}(\eta)$.

\smallskip

\noindent {\bf (2)} \ If $\eta \in {\sf riHom}(A^*)$ is an idempotent (for 
the operation of composition) then all essentially equal extensions and 
restrictions of $\eta$ are also idempotents of ${\sf riHom}(A^*)$.

\smallskip

\noindent {\bf (3)} \ If $\eta \in {\sf riHom}(A^*)$ is an idempotent 
then $\eta \in {\sf riHom}_{\sf pc}(A^*)$.
\end{lem}
{\bf Proof.} Statement (1) is a basic fact about composition of partial 
functions.

Proof of (2): If $\eta(y) = y$ for some $y \in {\sf Im}(\eta)$ then 
$\eta(yw) = yw$ for all $w \in A^*$. Hence, essentially equal restrictions 
of $\eta$ are also idempotents.

If $ya_1, \ldots, ya_k \in {\sf Im}(\eta)$ and
$\eta(ya_1) = ya_1, \ldots, \eta(ya_k) = ya_k$ then an extension of $\eta$ 
will be a function $\eta'$ with the additional mapping $\eta'(y) = y$. This
preserves the condition for an idempotent.
 
Proof of (3): Since $\eta$ is an idempotent, it follows from (1) that 
${\sf Im}(\eta) \subseteq {\sf Dom}(\eta)$. Let 
$y_1, y_2 \in {\sf imC}(\eta) \subset {\sf Dom}(\eta)$. Then there exist
$p_1, p_2 \in {\sf domC}(\eta)$ and $u_1, u_2 \in A^*$ such that 
$y_1 = p_1 u_1$ and $y_2 = p_2 u_2$. If $y_2$ and $y_1$ were prefix-comparable
then $p_1$ and $p_2$ would also be prefix-comparable, contradicting the 
fact that ${\sf domC}(\eta)$ is a prefix code. 
 \ \ \ $\Box$

\begin{pro} \label{idempot_characteriz_2} \    
An element $\eta \in {\sf riHom}(A^*)$ represents an idempotent of $M_{k,1}$
 \ iff \ $\eta$ is an idempotent of ${\sf riHom}(A^*)$ (for the operation
of composition). 

An element $\varphi \in M_{k,1}$ is an idempotent iff at least one 
representative of $\varphi$ in ${\sf riHom}(A^*)$ is an idempotent, iff all 
representatives of $\varphi$ in ${\sf riHom}(A^*)$ are idempotents.
In other words, the inverse of the function
 \ $\varphi \in {\sf riHom}(A^*) \ \longmapsto \ {\sf max}(\varphi)$
$ \in M_{k,1}$ \ preserves idempotents.
\end{pro}
{\bf Proof.} We first prove the following. 

\smallskip
 
\noindent {\bf Claim.} \ 
{\it If for $\eta \in {\sf riHom}(A^*)$ we have $\eta = {\sf max}(\eta)$ 
and ${\sf max}(\eta \circ \eta) = \eta$, then
$\eta \circ \eta = \eta$. 
}

\smallskip
 
\noindent {\sf Proof of Claim.} For any $x_i \in {\sf domC}(\eta)$ and any 
$w \in A^*$ we have $\eta(x_i w) = y_i w$ for some $y_i \in {\sf imC}(\eta)$.
We also have $\eta \circ \eta(x_i w) = \eta(y_i w)$ if $w$ is long enough 
so that $y_i w \in {\sf Dom}(\eta)$. 
Since ${\sf max}(\eta \circ \eta) = \eta$, we then have $\eta(y_i w) = y_i w$.
Hence, for all $y_i \in {\sf imC}(\eta)$ and all long enough $w \in A^*$ we
have: \ $\eta(y_i w) = y_i \, w$. Since $\eta$ was assumed to be maximally
essentially extended it follows that for all $y_i \in {\sf imC}(\eta)$:  
 \ $\eta(y_i) = y_i$. Therefore, for all $x_i \in {\sf domC}(\eta)$: \   
$\eta \circ \eta(x_i) = \eta(y_i) = y_i$, so $\eta \circ \eta = \eta$.
This proves the Claim.

\smallskip
 
We complete the proof of the Proposition. 
If $\eta \in {\sf riHom}(A^*)$ is an idempotent then it represents an
idempotent of $M_{k,1}$, since $M_{k,1}$ is a homomorphic image of 
${\sf riHom}(A^*)$. 

If $\varphi \in M_{k,1}$ is an idempotent then $\varphi$ can be represented
by $\eta \in {\sf riHom}(A^*)$ such that $\eta = {\sf max}(\eta)$, and
${\sf max}(\eta \circ \eta) = \eta$. By the Claim, $\eta$ is an idempotent 
of ${\sf riHom}(A^*)$. 

Moreover, if $\eta = {\sf max}(\eta)$ is an idempotent of ${\sf riHom}(A^*)$
then by Lemma \ref{idempot_characteriz_1}, all its essentially equal 
restrictions (i.e., all representatives of $\varphi$) are idempotents of 
${\sf riHom}(A^*)$.
 \ \ \ $\Box$

\begin{defn} \label{func_congr} \  
With a prefix code congruence $\simeq$ we associate the following two 
right-ideal homomorphisms, ${\sf func}_0(\simeq)$, \ ${\sf func}_1(\simeq)$
$ \in {\sf riHom}(A^*)$. Both have domain ${\sf Dom}(\simeq)$, and they are
defined by

\medskip

${\sf func}_0(\simeq): \ \ $
$x \in {\sf Dom}(\simeq) \ \longmapsto \ {\sf min}_{\sf dict}([x]) \ \in [x]$ , 

\medskip

${\sf func}_1(\simeq): \ \ $
$ x \in {\sf Dom}(\simeq) \ \longmapsto \ {\sf max}_{\sf dict}([x]) \ \in [x]$ ,

\medskip
 
\noindent where $[x]$ denotes the $\simeq$-class of $x$, and 
${\sf min}_{\sf dict}([x])$ or ${\sf max}_{\sf dict}([x])$ denotes the 
minimum, respectively maximum, element in $[x]$ according to the 
dictionary order on $A^*$.
\end{defn}
It follows from Lemma \ref{idempot_characteriz_1} and Prop.\ 
\ref{idempot_characteriz_2} that 
${\sf func}_0(\simeq)$ and ${\sf func}_1(\simeq)$ are idempotents, both as
elements of ${\sf riHom}(A^*)$ and of $M_{k,1}$.

\begin{pro} \label{func_vs_part} \
The operations ${\sf part}$ and ${\sf func}$ are inverses of each other,
in the following sense: 

\smallskip

\noindent {\bf (1)} \ For any prefix code congruence $\simeq$ and for
$j = 0,1$, we have:

\smallskip

 \ \ \ \ \  ${\sf part}({\sf func}_j(\simeq)) \ = \ \ \simeq$ . 

\smallskip

\noindent {\bf (2)} \ For any right-ideal homomorphism $\varphi$ and for
$j = 0,1$, we have:

\smallskip

 \ \ \ \ \ ${\sf func}_j({\sf part}(\varphi)) \ \equiv_{\cal L} \ \varphi$ ,

\smallskip

 \ \ where $\equiv_{\cal L}$ is the $\cal L$-equivalence of 
  ${\sf riHom}(A^*)$.

 \ \ Hence (by Prop.\ \ref{idempot_characteriz_2}), we also have 
\ ${\sf func}_j({\sf part}(\varphi)) \equiv_{\cal L} \varphi$ 
 \ for the $\equiv_{\cal L}$-relation of $M_{k,1}$. 
\end{pro}
{\bf Proof.} \ Part (1) follows immediately from the definitions of 
${\sf part}$ and ${\sf func}_j$. 

For part (2) let $\varphi: PA^* \to QA^*$ be a right-ideal homomorphism 
where $P$ and $Q$ are finite prefix codes. 
Let ${\sf part}(\varphi) = \{P_1, \ldots, P_m\}$, where $m = |Q|$, and 
let $Q = \{q_1, \ldots, q_m\}$, where $\{q_i\} = \varphi(P_i)$.
We give the proof for ${\sf func}_1$; for ${\sf func}_0$ the proof is the 
same.  Let $f_1$ be a short-hand for ${\sf func}_1({\sf part}(\varphi))$. 
We want to show that $f_1 \geq_{\cal L} \varphi$ and 
$\varphi \geq_{\cal L} f_1$.

For all $p \in P_i$ ($i = 1, \ldots, m$) we have 
 \ $\varphi \circ f_1(p) = \varphi({\sf max}_{\sf dict}(P_i)) = \varphi(p)$, 
since $p$ and ${\sf max}_{\sf dict}(P_i)$ belong to $P_i$. So, 
$f_1 \geq_{\cal L} \varphi$.

Let $\psi \in {\sf riHom}(A^*)$ be defined by
$\psi(q_i) = {\sf max}_{\sf dict}(P_i)$ for $i = 1, \ldots, m$; so, 
${\sf domC}(\psi) = Q$.  Then for all $p \in P_i$ (for $i = 1, \ldots, m$) 
we have 
 \ $\psi \circ \varphi(p) = \psi(q_i) = {\sf max}_{\sf dict}(P_i) = f_1(p)$. 
Hence, $\varphi \geq_{\cal L} f_1$. 
  \ \ \ $\Box$

\medskip

\noindent 
It follows from Prop.\  \ref{func_vs_part}(1) that every prefix code 
congruence is the partition of some right-ideal homomorphism.

\begin{lem} \label{idempot_max} \ Let $\simeq$ be a prefix code congruence.
Then $\simeq$ is maximally extended \ iff \ ${\sf func}_j(\simeq)$ is 
maximally extended. It follows that (for $j = 0,1$), \

\smallskip
   
${\sf func}_j({\sf max}(\simeq)) \ = \ {\sf max}({\sf func}_j(\simeq))$
$ \ = \ {\sf max}({\sf func}_j({\sf max}(\simeq)))$.
\end{lem}
{\bf Proof.} An extension of $\simeq$ is possible iff $\simeq$ contains the
classes $Ca_1, \ldots, Ca_k$, but not $C$. This is equivalent to having
${\sf func}_1(\simeq)$ mapping $Ca_i$ to 
${\sf max}_{\sf dict}(Ca_i) = {\sf max}_{\sf dict}(C) \ a_i$ for 
$i = 1, \ldots, k$. 
So $\simeq$ is extendable iff ${\sf func}_1(\simeq)$ is extendable. 
The proof for ${\sf func}_0(\simeq)$ is the same.
 \ \ \ $\Box$

\begin{pro} \label{func_Mk1} \ For prefix code congruences $\simeq_1$ and 
$\simeq_2$ the following are equivalent: 

\smallskip

\noindent $\bullet$ \ \ \ $\simeq_1$ and $\simeq_2$ are essentially equal, 

\smallskip

\noindent $\bullet$ \ \ \ ${\sf func}_0(\simeq_1) = {\sf func}_0(\simeq_2)$ 
 \ in $M_{k,1}$ , 

\smallskip

\noindent $\bullet$ \ \ \ ${\sf func}_1(\simeq_1) = {\sf func}_1(\simeq_2)$ 
 \ in $M_{k,1}$. 
\end{pro}
{\bf Proof.} We prove this only for ${\sf func}_1$; for ${\sf func}_0$ the 
proof works in the same way. If $\simeq_1 \ =_{\sf ess} \ \simeq_2$ then 
${\sf max}(\simeq_1) = {\sf max}(\simeq_2)$, by Proposition
\ref{congr_Unique_max}(2).  Hence by Lemma \ref{idempot_max}, 
${\sf max}({\sf func}_1(\simeq_1)) = {\sf max}({\sf func}_1(\simeq_2))$, 
hence, 
${\sf func}_1(\simeq_1) = {\sf func}_1(\simeq_2)$ \ in $M_{k,1}$. 

Conversely, if ${\sf func}_1(\simeq_1) = {\sf func}_1(\simeq_2)$ in
$M_{k,1}$ then ${\sf part}({\sf func}_1(\simeq_1)) =_{\sf ess} $
${\sf part}({\sf func}_1(\simeq_2))$, by 
Prop.\ \ref{part_invariant}.
By Prop.\ \ref{func_vs_part}(1), 
${\sf part}({\sf func}_1(\simeq_1))  = \ \simeq_1$, and similarly for 
$\simeq_2$. Hence, $\simeq_1 \ =_{\sf ess} \ \simeq_2$. 
 \ \ \ $\Box$


\subsection{Refinements of prefix code congruences}

\begin{defn} \label{refine_ess_partition} \
Let $\simeq_1$ and $\simeq_2$ be prefix code congruences.
We say that $\simeq_1$ is an {\bf end refinement} of $\simeq_2$ \ iff 
 \ there exist essentially right congruences $\simeq_1'$ and $\simeq_2'$
such that:

\smallskip

\noindent $\bullet$ \ $\simeq_i'$ is an essentially equal restriction of 
$\simeq_i$ (for $i = 1, 2$), 

\smallskip

\noindent
$\bullet$ \ ${\sf domC}(\simeq_2') \ \subseteq \ {\sf domC}(\simeq_1')$ ,
 \ and

\smallskip

\noindent
$\bullet$ \ every class of $\simeq_2'$ is a union of classes of $\simeq_1'$.
\end{defn}
\noindent {\bf Notation:} \  If $\simeq_1$ is an end refinement of 
$\simeq_2$ we denote this by 
 \ \ \ $\simeq_2 \ \ \leq_{\sf end} \ \ \simeq_1$ .

\begin{lem} \label{refine_ess_partition_ext} \ Let $\simeq_2'$ and 
$\simeq_1'$ be prefix code congruences.  

\smallskip

\noindent   
{\bf (1)} \ If \ $\simeq_2' \ \leq_{\sf end} \ \simeq_1'$ \ then
 \ ${\sf Dom}(\simeq_2') \subseteq_{\sf end} \ {\sf Dom}(\simeq_1')$.  

\smallskip

\noindent {\bf (2)} \  
Assume that every class of $\simeq_2'$ is a union of classes of $\simeq_1'$,
and that \ ${\sf domC}(\simeq_2') \subseteq {\sf domC}(\simeq_1')$.  And
assume $\simeq_2'$ is essentially extendable, in one replacement step
(\ref{congr_InvReplacement}), to $\simeq_2$.
Then $\simeq_1'$ is essentially extendable to a prefix code congruence
$\simeq_1$ such that every $\simeq_2$-class is a union of $\simeq_1$-classes,
and \ ${\sf domC}(\simeq_2) \subseteq {\sf domC}(\simeq_1)$.  

\smallskip

\noindent {\bf (3)} \
Assume that every $\simeq_2'$-class is a union of $\simeq_1'$-classes, and 
that \ ${\sf domC}(\simeq_2') \subseteq {\sf domC}(\simeq_1')$. And
assume that $\simeq_2'$ can be essentially restricted, in one replacement step 
(\ref{congr_replacement}), to $\simeq_2''$. 
Then $\simeq_1'$ can be essentially restricted to a prefix code 
congruence $\simeq_1''$ such that every $\simeq_2''$-class is a union of 
$\simeq_1''$-classes and
 \ ${\sf domC}(\simeq_2'') \subseteq {\sf domC}(\simeq_1'')$. 
\end{lem}
{\bf Proof.} (1) Let $\simeq_1'$, $\simeq_2'$ be as in Definition
\ref{refine_ess_partition}. \\     
Then ${\sf domC}(\simeq_2') \subseteq {\sf domC}(\simeq_1')$ implies 
${\sf Dom}(\simeq_2') \subseteq {\sf Dom}(\simeq_1')$.
Therefore, ${\sf Dom}(\simeq_2) =_{\sf ess} {\sf Dom}(\simeq_2') \subseteq $
${\sf Dom}(\simeq_1') =_{\sf ess} {\sf Dom}(\simeq_1)$, since $\simeq_1'$ and 
$\simeq_2'$ are essentially equal restrictions of $\simeq_1$, respectively 
$\simeq_2$. Hence, 
${\sf Dom}(\simeq_2) \subseteq_{\sf end} {\sf Dom}(\simeq_1)$.

(2) Let $A = \{a_1, \ldots, a_k\}$.
If $\simeq_2'$ is extendable, it has a subset of classes of the form
$Ca_1$, $\ldots$, $Ca_k$, with 
 \ $Ca_1 \cup \ldots \cup Ca_k \subseteq {\sf domC}(\simeq_2')$. 
Since every class of $\simeq_2'$ is a union of $\simeq_1'$-classes, we have
for each $i = 1, \ldots, k:$ \ $Ca_i = \bigcup_{j} Q_{i,j}$, where each 
$Q_{i,j}$ is a $\simeq_1'$-class, and $Q_{i,j} \subset {\sf domC}(\simeq_1')$.
It follows that $Q_{i,j}$ has the form $Q_{i,j} = P_ja_i$ for all $i, j$, and
that $\bigcup_{j} P_j = C$.
Hence $\simeq_1'$ contains the classes $P_j a_1$, $\ldots$, $P_j a_k$.
So, $\simeq_1'$ can be essentially extended to a prefix code congruence 
$\simeq_1$ by adding the classes $P_j$ to $\simeq_1'$ for all $j$. 
The domain code of $\simeq_1$ is obtained from ${\sf domC}(\simeq_2')$ by
replacing the set \ $\bigcup_{j,i} P_j a_i$ \ by \ $\bigcup_{j} P_j$.
Since $\bigcup_{j} P_j = C$, it follows that every $\simeq_2$-class is a 
union of $\simeq_1$-classes, and that
 \ ${\sf domC}(\simeq_2') \subseteq {\sf domC}(\simeq_1')$.

(3) The proof is similar to the proof of (2). \ \ \ $\Box$


\bigskip

\noindent The following generalizes Definition 5.2 from \cite{BiCoNP}.
 
\begin{defn} \label{complement_pref_code} \ 
Let $Q, P \subset A^*$ be finite prefix codes such that
 $PA^* \subset QA^*$.
A {\bf complement of $P$ in $QA^*$} is any finite prefix code 
$C \subset A^*$ such that \ $CA^* \cap PA^* = \varnothing$ \ and 
 \ $CA^* \cup PA^* \ =_{\sf ess} \ QA^*$.
 
For the ends this means: 
 \ ${\sf end}(CA^*) \cap {\sf end}(PA^*) = \varnothing$, and 
 \ ${\sf end}(CA^*) \cup {\sf end}(PA^*) = {\sf end}(QA^*)$. 
\end{defn}

\begin{lem} \label{exist_complement_pref_code} \ 
Let $Q, P \subset A^*$ be finite prefix codes such that
 $PA^* \subset QA^*$.
Then there exists a complement of $P$ in $QA^*$. 
\end{lem}
{\bf Proof.} Let $\ell = {\sf max}\{ |p| : p \in P\}$, i.e., $\ell$ is
the length of the longest word in $P$.  We pick  

\smallskip 

 \ \ \ $C \ = \ \{ x \in QA^* - PA^* \ : \ |x| = \ell \}$ . 

\smallskip 

\noindent Since all elements of $C$ have the same length, $C$ is a 
finite prefix code. Also, the definition of $C$ immediately implies that 
$CA^* \cap PA^* = \varnothing$.

Let us prove that $CA^* \cup PA^* \ =_{\sf ess} \ QA^*$. It is enough to
show that every end $w \in {\sf ends}(QA^*)$ that does not pass trough $P$ 
passes trough $C$. The latter is true, since the prefix $x$ of $w$ of 
length $\ell$ belongs to $C$, by the definition of $C$.
 \ \ \ $\Box$

\begin{pro} {\bf (Characterizations of $\leq_{\sf end}$).} 
\label{refine_ess_partition_DomMax} \   
Let $\simeq_2$ and $\simeq_1$ be prefix code congruences.
The following are equivalent:

\smallskip

\noindent {\bf (1)} \ $\simeq_2 \ \leq_{\sf end} \ \simeq_1$ ;

\smallskip

\noindent {\bf (2)} \ there is an essentially equal restriction $\simeq_2'$ 
of $\simeq_2$ such that every $\simeq_2'$-class is a union of 
$\simeq_1$-classes;

\smallskip

\noindent {\bf (3)} \ there is an essentially equal extension 
$\simeq_1^{\sharp}$ of $\simeq_1$ such that every $\simeq_2$-class is a 
union of $\simeq_1^{\sharp}$-classes;

\smallskip

\noindent {\bf (4)} \ every $\simeq_2$-class is a union of 
${\sf max}(\simeq_1)$-classes; 

\smallskip

\noindent {\bf (5)} \ every ${\sf max}(\simeq_2)$-class is a union of
${\sf max}(\simeq_1)$-classes;

\smallskip

\noindent {\bf (6)} \ ${\sf ends}({\sf Dom}(\simeq_2)) \ \subseteq \ $
${\sf ends}({\sf Dom}(\simeq_1))$, \ and 

\hspace{.04in}  every $\simeq_2$-class of ${\sf ends}({\sf Dom}(\simeq_2))$ 
is a union of $\simeq_1$-classes of ${\sf ends}({\sf Dom}(\simeq_1))$. 
\end{pro}
{\bf Proof.}
[(1) $\Rightarrow$ (2)] \ Let $\simeq_2'$ and $\simeq_1'$ be as in Definition
\ref{refine_ess_partition}; so every $\simeq_2'$-class is a union of 
$\simeq_1'$-classes. Moreover, every $\simeq_1'$-class is also a
$\simeq_1$-classes, since $\simeq_1'$ is an essentially equal restriction of 
$\simeq_1$. Thus, every $\simeq_2'$-class is a union of classes of $\simeq_1$.  

\smallskip

\noindent
[(2) $\Rightarrow$ (3)] \ Let $\simeq_2'$, $\simeq_2$, and $\simeq_1$ be as
in (2). Now we apply Lemma \ref{refine_ess_partition_ext}(2) to $\simeq_2'$ 
and $\simeq_1$, repeatedly, until $\simeq_2'$ has been rewritten to its 
essentially equal extension $\simeq_2$. In this process, $\simeq_1$ is 
rewritten to some essentially equal extension $\simeq_1^{\sharp}$ such that 
(3) holds. 

\smallskip

\noindent
[(3) $\Rightarrow$ (4)] \ If every $\simeq_2$-class is a union of
$\simeq_1^{\sharp}$-classes, then $\simeq_2$-class is also a union of
${\sf max}(\simeq_1)$-classes, since ${\sf max}(\simeq_1)$ is an end
refinement of $\simeq_1^{\sharp}$.

\smallskip

\noindent
[(4) $\Rightarrow$ (5)] \ We repeatedly apply Lemma
\ref{refine_ess_partition_ext}(2) to $\simeq_2$ and ${\sf max}(\simeq_1)$,
until $\simeq_2$ has been rewritten to ${\sf max}(\simeq_2)$. In the process,
${\sf max}(\simeq_1)$ does not change (being already maximally extended).
As a result, every ${\sf max}(\simeq_2)$-class is a union of
${\sf max}(\simeq_1)$-classes.

\smallskip

\noindent
[(5) $\Rightarrow$ (6)] \ It follows immediately from (5) that every element
of a $\simeq_2$-class is also in a $\simeq_1$-class; thus, \     
${\sf Dom}({\sf max}(\simeq_2)) \subseteq {\sf Dom}({\sf max}(\simeq_1))$.
Hence, \ ${\sf ends}({\sf Dom}(\simeq_2)) \subseteq $
${\sf ends}({\sf Dom}(\simeq_1))$.
 
Let $w_1, w_2 \in {\sf ends}({\sf Dom}(\simeq_2))$ be such that 
$w_1 \simeq_1 w_2$; we want to show that $w_1 \simeq_2 w_2$.
It follows from $w_1 \simeq_1 w_2$ that $w_1 = q_1v$ and $w_2 = q_2 v$ with
$q_1 \simeq_1 q_2$, for some $q_1, q_2 \in {\sf Dom}(\simeq_1)$,
$v \in A^{\omega}$. Moreover, $q_1 \simeq_1 q_2$ implies that $q_1 = p_1x$
and $q_2 = p_2x$ with  $p_1 \ {\sf max}(\simeq_1) \ p_2$, for some 
$p_1, p_2 \in {\sf Dom}({\sf max}(\simeq_1))$, $x \in A^*$.
By (5), this implies $p_1 \ {\sf max}(\simeq_2) \ p_2$, and hence
$q_1 \ {\sf max}(\simeq_2) \ q_2$ (since this is a right congruence).
Hence, $q_1 \simeq_2 q_2$ (since $q_1, q_2 \in {\sf Dom}(\simeq_1)$).
Therefore, $w_1 \simeq_2 w_2$ (by the definition of $\simeq_2$ on ends). 

\smallskip

\noindent
[(6) $\Rightarrow$ (1)] \ Since ${\sf Dom}(\simeq_1)$ and 
${\sf Dom}(\simeq_2)$ are finitely generated right ideals, there intersection
is also a finitely generated right ideal (by Lemma \ref{idealInter}). So,
${\sf Dom}(\simeq_2) \cap {\sf Dom}(\simeq_1) = P_2A^*$ \ for some finite
prefix code $P_2$. By (6), \ ${\sf ends}({\sf Dom}(\simeq_2)) \subseteq $
${\sf ends}({\sf Dom}(\simeq_1))$; hence  

\smallskip

${\sf ends}(P_2A^*) = $
${\sf ends}({\sf Dom}(\simeq_2)) \cap {\sf ends}({\sf Dom}(\simeq_1)) = $
${\sf ends}({\sf Dom}(\simeq_2))$. 

\smallskip

\noindent In other words, \ $P_2A^* \ =_{\sf ess} \ {\sf Dom}(\simeq_2)$.
Let $\simeq_2'$ be the essentially equal restriction of $\simeq_2$ to 
$P_2A^*$

By Lemma \ref{exist_complement_pref_code} there exists a complementary 
prefix code (let's call it $Q_1$) of $P_2$ within ${\sf Dom}(\simeq_1)$.
Since $Q_1$ is a complementary prefix code, 
 \ $(Q_1 \cup P_2) \, A^* \ =_{\sf ess} \ {\sf Dom}(\simeq_1)$. 
Let $\simeq_1'$ be the essentially equal restriction of $\simeq_1$ to 
 \ $(Q_1 \cup P_2) \, A^*$.

Then $\simeq_2'$ and $\simeq_1'$ satisfy the conditions of Definition
\ref{refine_ess_partition}, so (1) holds.  \ \ \ $\Box$

\bigskip

\noindent For prefix codes congruences $\simeq_1$ and $\simeq_2$ we have

\smallskip

 \ \ \ \ \  $\simeq_1 \ =_{\sf ess} \ \simeq_2$ \ \ iff \ \    
$\simeq_1 \ \leq_{\sf end} \ \simeq_2$ \ and \   
$\simeq_1 \ \geq_{\sf end} \ \simeq_2$.

\smallskip

\noindent This follows immediately from Propositions 
\ref{congr_ess_simil_characteriz}(2) and \ref{refine_ess_partition_DomMax}(5). 
Recall that for prefix codes congruences, $=_{\sf ess}$ was defined in
Def.\ \ref{ess_equ_congr}, and $\leq_{\sf end}$ was defined in 
Def.\ \ref{refine_ess_partition}.

\medskip

 The relation $\leq_{\sf end}$ is a {\em lattice pre-order} on the set of 
all prefix code congruences of $A^*$. 
For two prefix code congruences $\simeq_1$ and $\simeq_2$ we consider the 
prefix code congruence \ $\simeq_1 \wedge \simeq_2$, called the 
{\em meet} or {\em wedge}. Its domain is 
 \ ${\sf Dom}(\simeq_1 \wedge \simeq_2) \ = \ $
${\sf Dom}(\simeq_1) \cap {\sf Dom}(\simeq_2)$,   
and for all $u, v \in {\sf Dom}(\simeq_1 \wedge \simeq_2)$ we have:
 \ \ $u \ (\simeq_1 \wedge \simeq_2) \ v$ \ \ iff \ \ $u \simeq_1 v$ \ and 
 \ $u \simeq_2 v$ .
Similarly, the {\em join} \ $\simeq_1 \vee \simeq_2$ \ has domain
 \ ${\sf Dom}(\simeq_1 \vee \simeq_2) \ = \ $
${\sf Dom}(\simeq_1) \cup {\sf Dom}(\simeq_2)$, 
and is defined by the transitive closure of the relation \    
$\simeq_1 \cup \simeq_2 \ = \ $
$\{ (u,v) \in A^* \times A^* : \ u \simeq_1 v \ {\sf or} \ u \simeq_2 v \}$.
Equivalently, we start with all the classes that belong to $\simeq_1$ or 
$\simeq_2$, and we iteratively replace classes that intersect by their 
union, until no two classes intersect.  


\subsection{Characterization of the $\cal L$-order of $M_{k,1}$ }

\begin{thm} \label{L_order} {\bf (${\cal L}$-order of $M_{k,1}$).} \ \
For any \ $\varphi, \psi \in M_{k,1} :$

\medskip

 \ \ \ \ \ $\psi(.) \ \leq_{\cal L} \ \varphi(.)$ \ \ \ iff \ \ \
${\sf part}(\psi) \ \leq_{\sf end} \ {\sf part}(\varphi)$
\end{thm}
{\bf Proof.}
[$\Rightarrow$] \ If $\psi \leq_{\cal L} \varphi$ then there exists 
$\alpha \in M_{k,1}$ such that $\psi$ and $\alpha \circ \varphi$ represent
the same element of $M_{k,1}$. Hence (by Lemma \ref{ess_equ_Im_Dom} and 
Prop.\ \ref{part_invariant}), 
 \ ${\sf Dom}(\psi) \ =_{\sf ess} \ {\sf Dom}(\alpha \circ \varphi)$, 
 \ ${\sf Im}(\psi) \ =_{\sf ess} \ {\sf Im}(\alpha \circ \varphi)$, \ and
 \ ${\sf part}(\psi) \ =_{\sf ess} \ {\sf part}(\alpha \circ \varphi)$. Also, 
when $\alpha \circ \varphi(x)$ is defined then $\varphi(x)$ must be defined,
so \ ${\sf Dom}(\alpha \circ \varphi) \subseteq {\sf Dom}(\varphi)$. 

Let $u, v \in {\sf Dom}(\alpha \circ \varphi)$. Then $u$ and $v$ are related
by ${\sf part}(\varphi)$ iff $\varphi(u) = \varphi(v)$, which implies 
$\alpha \circ \varphi(u) = \alpha \circ \varphi(v)$, hence $u$ and $v$ are
in the same ${\sf part}(\alpha \circ \varphi)$-class. It follows that
${\sf part}(\varphi)$ is a refinement of ${\sf part}(\alpha \circ \varphi)$. 
Since \ ${\sf part}(\psi) \ =_{\sf ess} \ {\sf part}(\alpha \circ \varphi)$, 
it follows that 
${\sf part}(\varphi)$ is an end refinement of ${\sf part}(\psi)$.

[$\Leftarrow$] \ If ${\sf part}(\varphi)$ is an end refinement of 
${\sf part}(\psi)$, then (by Definition \ref{refine_ess_partition}) there
exists an essentially equal restriction $\simeq_2'$ of ${\sf part}(\psi)$, 
and an essentially equal restriction $\simeq_1'$ of ${\sf part}(\varphi)$, 
such that every $\simeq_2'$-class is a union of $\simeq_1'$-classes. Let 
$P_2 = {\sf domC}(\simeq_2')$ and let $P_1 = {\sf domC}(\simeq_1')$. Let 
$\psi_0$ be the restriction of $\psi$ to $P_2A^*$, and let $\varphi_0$ be 
the restriction of $\varphi$ to $P_1A^*$. Since $\simeq_2'$ is an
essentially equal restriction, $\psi_0$ and $\psi$ represent the same 
element of $M_{k,1}$; similarly, $\varphi_0$ and $\varphi$ represent the 
same element of $M_{k,1}$. 

We define a right ideal homomorphism $\alpha$ with domain 
 \ ${\sf Dom}(\alpha) = {\sf Im}(\varphi_0)$ \ and image \ 
${\sf Im}(\alpha) = {\sf Im}(\psi_0)$, as follows:

\smallskip 

$\alpha: \ \ z \in {\sf Im}(\varphi_0) \ \longmapsto \ $
$\psi_0 ( \varphi_0^{-1}(z)) \ \in {\sf Im}(\psi_0)$ .

\smallskip 

\noindent Then $\alpha$ is a function. Indeed, $z \in {\sf Im}(\varphi_0)$ 
can be written as $z = \varphi_0(x)$ for any $x \in \varphi_0^{-1}(z)$.
Then $\alpha(z) = \psi_0 ( \varphi_0^{-1}(\varphi_0^{-1}(z)))$; and since
every $\simeq_2'$-class is a union of $\simeq_1'$-classes, the latter is a 
subset of $\psi_0 ( \psi_0^{-1}( \psi_0(z))) = \psi_0(z)$. Hence, 
$\alpha(z) = \psi_0 ( \varphi_0^{-1}(\varphi_0^{-1}(z))) = \psi_0(z)$, 
indepently of the choice of $x$.   

It follows also from the definition of $\alpha$ that for all 
$z \in {\sf Im}(\varphi_0) :$ \ $\alpha \circ \varphi_0(z) \ = \ $
$\psi_0 \circ \varphi_0^{-1} \circ \varphi_0(z)$, and since every
$\simeq_2'$-class is a union of $\simeq_1'$-classes, the latter is a 
subset of \ $\psi_0 \circ \psi_0^{-1} \circ \psi_0(z) \ = \ $ 
$\psi_0(z)$. So, $\psi_0 = \alpha \circ \varphi_0 \leq_{\cal L} \varphi_0$.
Since $\psi_0$ represents the same element as $\psi$ in $M_{k,1}$, and 
$\varphi_0$ represents the same element as $\varphi$ in $M_{k,1}$, we obtain 
$\psi \leq_{\cal L} \varphi$.   \ \ \ $\Box$


\section{Infinite $\cal L$- and $\cal R$-chains, and density }

\subsection{Embedding of the $\cal L$- and $\cal R$-orders into the 
idempotent order}

\begin{lem} \label{R_equiv_id} \ If $\varphi \in M_{k,1}$ is represented
by $\varphi: P_1A^* \to P_2A^*$ in ${\sf riHom}(A^*)$ where
$P_1 = {\sf domC}(\varphi)$ and $P_2 = {\sf imC}(\varphi)$
are prefix codes, then we have:
 \ \ \  $\varphi \ \equiv_{\cal R} \ {\sf id}_{P_2A^*}$.
\end{lem}
{\bf Proof.} [$\geq _{\cal R}$] \ \ Obviously,
 \ ${\sf id}_{P_2A^*} \circ \varphi(.) = \varphi(.)$, so
${\sf id}_{P_1A^*} \geq _{\cal R} \varphi$.

\smallskip

\noindent [$\leq_{\cal R}$] \ \ We want to define
$\alpha: P_1A^* \to P_2A^*$ so that
$\varphi \circ \alpha(.) = {\sf id}_{P_2A^*}$.
For every $y \in P_2$ we choose an element
$\overline{y} \in \varphi^{-1}(y)$ ($\subseteq P_1$);
next, for every $y \in P_2$ define \ $\alpha(y) = \overline{y}$.

Then, for every $y \in P_2$ we have: \ \
$\varphi \big( \alpha(y) \big) \ = \ \varphi \big( \overline{y} \big) \ $
$ \in \ \varphi \big( \varphi^{-1}(y) \big) \ = \ y$.
 \ \ \ $\Box$

\medskip

The following provides an embedding of the $\cal R$-order on the set
of $\cal R$-classes of $M_{k,1}$ into the idempotent order of $M_{k,1}$.

\begin{pro} \label{idempot_vs_Rorder} {\bf (embedding of the
$\cal R$-order into the idempotent order).}  \ 
Lef $X, Y \subseteq A^*$ be finitely generated right ideals. We have:

\smallskip

\noindent {\bf (1)} \ \ \
$Y \ =_{\sf ess} \ X$ \ iff \ ${\sf id}_Y \ = \ {\sf id}_X$ \ in $M_{k,1}$.

\smallskip

\noindent {\bf (2)} \ \ \ $Y \ \subseteq_{\sf end} \ X$ \ iff \
${\sf id}_Y \ \leq \ {\sf id}_X$ \ for the idempotent order of $M_{k,1}$.

\smallskip

\noindent
{\bf (3)} \ \ \ The $\cal R$-order (on $\cal R$-classes) of $M_{k,1}$ is
embedded in the idempotent-order order of partial

 \ \ \ \  identity elements of $M_{k,1}$.
\end{pro}
{\bf Proof.}  {\bf (1)} $[\Rightarrow]$ \ $Y \ =_{\sf ess} \ X$ \ implies that
$X$ and $Y$ have the same ends, hence $X \cap Y$ also has the same ends as
$X$ and $Y$. Hence, the restrictions of both ${\sf id}_Y$ and ${\sf id}_X$
to ${\sf id}_{X \cap Y}$ are essentially equal restrictions, hence \
${\sf id}_Y \ = \ {\sf id}_{X \cap Y} \ = \ {\sf id}_X$ \ in $M_{k,1}$.

\smallskip

\noindent
$[\Leftarrow]$ \ If ${\sf id}_Y = {\sf id}_X$ in $M_{k,1}$ then ${\sf id}_Y$
and ${\sf id}_X$ agree on some common essential subideal $Z$ of $X$ and $Y$.
Then both $X$ and $Y$ are essentially equal to $Z$, hence
$Y =_{\sf ess} X$.

\smallskip

\noindent
{\bf (2)} \ Since $\subseteq_{\sf end}$ is a lattice pre-order we have
 \ $Y \subseteq_{\sf end} X$ iff $Y =_{\sf ess} X \cap Y$. By {\bf (1)}
the latter is equivalent to ${\sf id}_Y = {\sf id}_{X \cap Y}$, and the latter
is equal to \ ${\sf id}_{X \cap Y} = {\sf id}_X \circ {\sf id}_Y = $
${\sf id}_Y \circ {\sf id}_X$. Moreover,
${\sf id}_Y = {\sf id}_Y {\sf id}_X \circ {\sf id}_Y = $
${\sf id}_Y \circ {\sf id}_X$ \ is equivalent to ${\sf id}_Y \leq {\sf id}_X$
for the idempotent order.

\smallskip

\noindent
{\bf (3)} \ This follows now from {\bf (2)} and Lemma \ref{R_equiv_id}.
 \ \ \ $\Box$

\begin{pro} \label{embed_L_in_idempotOrd} {\bf (embedding of the
$\cal L$-order into the idempotent order).}  \
The $\cal L$-order on the set of $\cal L$-classes of $M_{k,1}$ is embeddable
into the idempotent order of $M_{k,1}$ as follows.
The function \ ${\sf func}_j({\sf part}(.))$, defined on the set of 
$\cal L$-classes of $M_{k,1}$, is injective (for $j = 0,1$), and for any
$\psi, \varphi \in {\sf riHom}(A^*)$ and for $j = 0,1$, we have:

\smallskip

 \ \ \ \ \ $\psi \ \leq_{\cal L} \ \varphi$ \ \ \ iff \ \ \
${\sf func}_j({\sf part}(\psi)) \ \leq \ {\sf func}_j({\sf part}(\varphi))$

\smallskip

\noindent where $\leq$ is the idempotent order of $M_{k,1}$.
\end{pro}
{\bf Proof.}  The injectiveness of ${\sf func}_j(.)$, as a  from the set
of $\cal L$-classes of $M_{k,1}$ into the set of idempotents of $M_{k,1}$,
follows from Theorem \ref{L_order}. So we have an embedding, but we still
need to show that this is order-preserving.

\smallskip

\noindent $[ \Leftarrow ]$
 \ If ${\sf func}_j({\sf part}(\psi)) \leq {\sf func}_j({\sf part}(\varphi))$
then, since ${\sf func}_j({\sf part}(\chi) \equiv_{\cal L} \chi$ for all
$\chi$ (by Prop.\ \ref{func_vs_part}(2)), we obtain
 \ $\psi \leq_{\cal L} \varphi$.

\smallskip

\noindent $[ \Rightarrow ]$  \ If \ $\psi \leq_{\cal L} \varphi$ \ then (by
Prop.\ \ref{func_vs_part}(2)), \ ${\sf func}_j({\sf part}(\psi)) $
$ \leq_{\cal L} {\sf func}_j({\sf part}(\varphi))$. We also need to show
this order relation for $\leq_{\cal R}$, i.e., we need to show
 \ ${\sf Im}({\sf func}_j({\sf part}(\psi))) \subseteq_{\sf end}$
${\sf Im}({\sf func}_j({\sf part}(\varphi)))$. We prove the result only
for ${\sf func}_0$; for ${\sf func}_1$ the proof is similar.

Let the tables for $\psi$ and $\varphi$ be $\psi: S \to T$ and
$\varphi: P \to Q$, where $S, T, P, Q$ are finite prefix codes. After an
essential class-wise restriction (if necessary) we have $S \subseteq P$, and
every class $S_i$ of ${\sf part}(\psi)$ in $S$ is a union of classes
$P_1, \ldots, P_{n_i}$ of ${\sf part}(\varphi)$ in $S$ ($\subseteq P$);
this follows from $\psi \leq_{\cal L} \varphi$.
Then ${\sf func}_0({\sf part}(\psi))$ maps all of $S_i$ to
 \ ${\sf min}_{\sf dict}(S_i) \in S_i = \bigcup_{j=1}^{n_i} P_j$.
So there is $j$ with ${\sf min}_{\sf dict}(S_i) \in P_j$; hence,
${\sf min}_{\sf dict}(S_i) = {\sf min}_{\sf dict}(P_j)$. Therefore,
${\sf func}_0({\sf part}(\psi))$ maps all of $S_i$ to an element of
${\sf Im}({\sf func}_0({\sf part}(\varphi)))$. Therefore,
${\sf Im}({\sf func}_j({\sf part}(\psi))) $    $\subseteq $
${\sf Im}({\sf func}_j({\sf part}(\varphi)))$.
    \ \ \ $\Box$

\subsection{Density}

\begin{pro} \label{R_density} {\bf (Density of the $\cal R$-order).}
 \  For any two elements $\varphi, \psi \in M_{k,1}$ such that
$\varphi >_{\cal R} \psi$ there exists $\chi \in M_{k,1}$ such that
$\varphi >_{\cal R} \chi >_{\cal R} \psi$.
\end{pro}
{\bf Proof.}  Let $P = {\sf imC}(\varphi)$ and $Q = {\sf imC}(\psi)$. After
applying essentially equal reductions, if necessary, we can assume that $P$
and $Q$ are prefix codes. Since $\varphi >_{\cal R} \psi$ we have
$PA^* \neq_{\sf ess} QA^*$, by Theorem \ref{R_order}.  Hence, $PA^*$
contains some end $\eta$ which does not belong to $QA^*$. Let $w \in PA^*$
be a prefix of this end that is strictly longer
than any element in $P$ or $Q$,
and let $p$ be the prefix of $w$ that belongs to $P$.
Consider the right ideal \ $DA^* = (Q \cup \{w\})A^*$.
Then $QA^* \subset DA^* \subset \subset PA^*$.
Moreover, $DA^* \neq_{\sf ess} QA^*$ since $DA^*$ contains the end $\eta$.
And $DA^* \neq_{\sf ess} PA^*$ since $PA^*$ contains all ends with prefix
$p$, while $DA^*$ doesn't ($p$ being strictly shorter than $w$).
Hence (by Lemma \ref{R_equiv_id} and Theorem \ref{R_order}), we have
 \ $\varphi >_{\cal R} {\sf id}_{DA^*} >_{\cal R} \psi$.
 \ \ \ \ \ $\Box$

\begin{pro} \label{L_order_density} {\bf (Density of the $\cal L$-order).}
\  For any two elements $\varphi, \psi \in M_{k,1}$ such that
$\varphi >_{\cal L} \psi$ there exists $\chi \in M_{k,1}$ such that
$\varphi >_{\cal L} \chi >_{\cal L} \psi$.
\end{pro}
{\bf Proof.}  By Theorem \ref{L_order}, 
${\sf Dom}(\psi) \subseteq_{end} {\sf Dom}(\varphi)$ 
and ${\sf part}(\varphi)$ is an end-refinement of ${\sf part}(\psi)$.
To construct $\chi$ we distinguish two cases: 
Case (1), when ${\sf Dom}(\psi) \varsubsetneq_{\sf end} {\sf Dom}(\varphi)$; 
case (2), when ${\sf Dom}(\psi) =_{\sf ess} {\sf Dom}(\varphi)$, and 
${\sf part}(\varphi)$ is a strict end-refinement of ${\sf part}(\psi)$.
By Prop.\ \ref{refine_ess_partition_DomMax}, we can essentially restrict 
$\varphi$ and $\psi$ so that ${\sf domC}(\psi) \subseteq {\sf domC}(\varphi)$ 
and ${\sf part}(\varphi)$ is a finer partition of ${\sf domC}(\psi)$ than 
${\sf part}(\psi)$. So, after essentially equal restrictions, case (1) 
becomes ${\sf domC}(\psi) \varsubsetneq {\sf domC}(\varphi)$; case (2) 
becomes ${\sf domC}(\psi) = {\sf domC}(\varphi)$ and ${\sf part}(\varphi)$ 
is strictly finer than ${\sf part}(\psi)$. 

\medskip

\noindent {\sf Case (1):} 
 \ ${\sf domC}(\psi) \varsubsetneq {\sf domC}(\varphi)$.

\smallskip

\noindent Consider the prefix code \ $\{u_1, \ldots, u_n\} = $  
${\sf domC}(\varphi) - {\sf domC}(\psi)$. We take the essentially equal 
restriction of $\varphi$ to 
 \ $u_1A \cup \{u_2, \ldots, u_n\} \cup {\sf domC}(\psi)$.
Then we define $\chi$ to be the (non-essential) restriction of $\varphi$ to
 \ $u_1 (A - \{a_1\}) \cup \{u_2, \ldots, u_n\} \cup {\sf domC}(\psi)$.
This implies \ ${\sf Dom}(\psi) \varsubsetneq_{\sf end} {\sf Dom}(\chi) $
$\varsubsetneq_{\sf end} {\sf Dom}(\varphi)$.

\medskip

\noindent {\sf Case (2):} \ ${\sf domC}(\psi) = {\sf domC}(\varphi)$ and 
${\sf part}(\varphi)$ is strictly finer than ${\sf part}(\psi)$.

\smallskip

In this case will actually do not use the fact that 
${\sf domC}(\psi) = {\sf domC}(\varphi)$ in order to construct $\chi$. 

To simplify the notation we only prove case (2) when $k=2$ and 
$A = \{a,b\}$, but the same method works in general.

Let $\{Q_1, \ldots, Q_n\}$ be the classes of ${\sf part}(\psi)$ on 
${\sf domC}(\psi)$, and let $\{P_1, \ldots, P_m\}$ be the classes of
${\sf part}(\varphi)$ on ${\sf domC}(\varphi)$.
So, ${\sf imC}(\psi) = \{y_1, \ldots, y_n\}$, where $y_i = \psi(Q_i)$ for 
$1 \leq i \leq n$;
this is one word since $Q_i$ is a class of ${\sf part}(\psi)$.
By renaming the classes of ${\sf part}(\psi)$, if necessary, we can assume 
that $Q_1$ is the union of at least two classes of ${\sf part}(\varphi)$: 
 \ $Q_1 = P_1 \cup \ldots \cup P_s$, for some $s$ with $2 \leq s \leq m$.
We take an essentially equal restriction $\psi'$ of $\psi$ such that 
${\sf part}(\psi')$ on ${\sf domc}(\psi')$ is \   
$\{Q_1a, Q_1b, Q_2, \ldots, Q_n\}$. 
Similarly, we essentially restrict $\varphi$ to $\varphi'$ so that 
${\sf part}(\varphi')$ on ${\sf domc}(\varphi')$ is \  

\medskip

$\{P_1a, \ldots, P_s a\} \ \cup \ \{P_1b, \ldots, P_s b\}$  $ \ \cup \ $ 
$\{ P_{s+1}, \ldots , P_m\}$. 

\smallskip

\noindent So, ${\sf domc}(\psi') = {\sf domc}(\varphi')$, 
 \ $P_1a \cup \ldots \cup P_s a = Q_1a$, and 
 \ $P_1b \cup \ldots \cup P_s b = Q_1b$. 
We now define $\chi$ as follows, where $\{c_1, \ldots, c_s\}$ is any 
prefix code of cardinality $s$ :

\smallskip

\noindent $\bullet$ \ \ \ 
${\sf domC}(\chi) = {\sf domC}(\psi')$ , 

\smallskip

\noindent $\bullet$ \ \ \  
${\sf part}(\chi)$ on ${\sf domC}(\chi)$ is 
$ \ \{ P_1a, \ldots, P_s a\} \ \cup \ \{Q_1b, \ Q_2, \ldots, Q_n\}$, 

\smallskip

\noindent $\bullet$ \ \ \ the values of $\chi$ on ${\sf domC}(\chi)$ are

\smallskip

 \ \ $\left\{ \begin{array}{l} 
   \chi(P_j a) = y_1 a c_j \ \ {\rm for} \ j = 1, \ldots, s, \\  
   \chi(Q_1b) =  y_1 b , \\  
   \chi(Q_j) = y_j \ \ {\rm for} \ 2 \leq j \leq n.
           \end{array} \right. $

\medskip

\noindent Recall that $y_i = \psi(Q_i)$, $1 \leq i \leq n$.

Then on ${\sf domC}(\chi)$ ($= {\sf domC}(\psi') = $
${\sf domc}(\varphi')$) we have: \ ${\sf part}(\varphi')$ refines 
${\sf part}(\chi)$, and ${\sf part}(\chi)$ refines ${\sf part}(\psi')$.
Hence Theorem \ref{L_order} and Prop.\ \ref{refine_ess_partition_DomMax}
imply $\psi <_{\cal L} \chi <_{\cal L} \varphi$.
  \ \ \ $\Box$

\begin{pro} \label{L_order_density_Rclass} {\bf (Density of the 
$\cal L$-order within an $\cal R$-class).} 

\noindent {\bf (1)} \ For any two elements $\varphi, \psi \in M_{k,1}$ 
with $\varphi \equiv_{\cal R} \psi$ and $\varphi >_{\cal L} \psi$, there 
exists $\chi \in M_{k,1}$ such that 

 \ $\varphi \equiv_{\cal R} \chi \equiv_{\cal R} \psi$ \ and \  
$\varphi >_{\cal L} \chi >_{\cal L} \psi$.

\noindent {\bf (2)} \ For every $\varphi \in M_{k,1} - \{{\bf 0}\}$  there
exists $\chi \in M_{k,1}$ such that
 \ $\chi \equiv_{\cal R} \varphi$ \ and \ $\varphi >_{\cal L} \chi$. 
 
And for every $\psi \in M_{k,1}$ such that $\psi$ is not $\cal L$-maximal,
there exists $\chi \in M_{k,1}$ such that

 \ $\chi \equiv_{\cal R} \psi$ \ and \ $\chi >_{\cal L} \psi$. 
\end{pro}
{\bf Proof, (1).} \ As in the proof of Prop.\ \ref{L_order_density}, we 
consider two cases, and we take essentially equal class-wise restrictions. 
Then we have: 
Case (1), ${\sf domC}(\psi) \varsubsetneq {\sf domC}(\varphi)$; case (2),
${\sf domC}(\psi) = {\sf domC}(\varphi)$ and ${\sf part}(\varphi)$ is 
strictly finer than ${\sf part}(\psi)$.

To simplify the notation we only give the proof when $k=2$ and $A = \{a,b\}$,
but the same method works in general.

\medskip

\noindent {\sf Case (1):}
 \ ${\sf domC}(\psi) \varsubsetneq {\sf domC}(\varphi)$.

\smallskip

\noindent Consider the prefix code \ $\{u_1, \ldots, u_n\} = $
${\sf domC}(\varphi) - {\sf domC}(\psi)$. We take the essentially equal 
restriction $\varphi'$ of $\varphi$ to 
 \ $\{u_1a, u_1b, u_2, \ldots, u_n\} \cup {\sf domC}(\psi)$.
Then we define $\chi$ on ${\sf domC}(\chi) = {\sf domC}(\psi) \cup \{u_1a\}$
by 

\medskip

$\left\{ \begin{array}{l}
   \chi(x) = \psi(x) \ \ {\rm if} \ x \in {\sf domC}(\psi) , \\
   \chi(u_1a) = \varphi(u_1) \, a.
           \end{array} \right. $

\medskip
 
\noindent Then ${\sf Im}(\psi) \subseteq {\sf Im}(\chi) \subseteq $
${\sf Im}(\varphi')$, hence $\chi \equiv_{\cal R} \psi$ (since 
${\sf Im}(\psi) =_{\sf ess} {\sf Im}(\varphi')$).

Also, ${\sf Dom}(\psi) \varsubsetneq_{\sf end} {\sf Dom}(\chi)$ (since 
$u_1a$ and its ends are missing from ${\sf Dom}(\psi)$), and 
${\sf Dom}(\chi) \varsubsetneq_{\sf end} {\sf Dom}(\varphi)$ (since $u_1b$
and its ends are missing from ${\sf Dom}(\chi)$). So,
$\psi <_{\cal L} \chi <_{\cal L} \varphi$.

\medskip

\noindent {\sf Case (2):} \ ${\sf domC}(\psi) = {\sf domC}(\varphi)$ and
${\sf part}(\varphi)$ is strictly finer than ${\sf part}(\psi)$.

\smallskip

In this case will actually do not use the fact that
${\sf domC}(\psi) = {\sf domC}(\varphi)$ in order to construct $\chi$.

Let $\{Q_1, \ldots, Q_n\}$ be the classes of ${\sf part}(\psi)$ on
${\sf domC}(\psi)$, and let $\{P_1, \ldots, P_m\}$ be the classes of
${\sf part}(\varphi)$ on ${\sf domC}(\varphi)$.
So, ${\sf imC}(\psi) = \{y_1, \ldots, y_n\}$, where $y_i = \psi(Q_i)$ for
$1 \leq i \leq n$;
this is one word since $Q_i$ is a class of ${\sf part}(\psi)$.
By renaming the classes of ${\sf part}(\psi)$, if necessary, we can assume
that $Q_1$ is the union of at least two classes of ${\sf part}(\varphi)$:
 \ $Q_1 = P_1 \cup \ldots \cup P_s$, for some $s$ with $2 \leq s \leq m$.
We take the essentially equal class-wise restriction $\psi'$ of $\psi$ so 
that \ ${\sf part}(\psi') = \{Q_1a, Q_1b, Q_2, \ldots, Q_n\}$, 
and we define $\chi$ on ${\sf domC}(\chi) = {\sf domC}(\psi')$ by
 
\medskip

$\left\{ \begin{array}{l}
 \chi(Q_i) \ = \ y_i \ \ {\rm for} \ 2 \leq i \leq n , \\  
 \chi(Q_1b) \ = \ y_1b , \\  
 \chi(P_1a) \ = \ y_1 a a ,  \\   
 \chi(Q_1a - P_1a) \ = \ y_1 a b .
\end{array} \right. $

\medskip

\noindent Then ${\sf imC}(\chi) = \{y_1aa, y_1ab, y1b, y_2, \ldots, y_n\}$,
hence, since $\{aa, ab, b\}$ is a maximal prefix code, we have 
 \ ${\sf Im}(\chi) =_{\sf ess} {\sf Im}(\psi)$;  so 
$\chi \equiv_{\cal R} \psi$.

And ${\sf part}(\chi)$ is finer than ${\sf part}(\psi')$ since $Q_1a$ is
partitioned into $P_1a$ and $Q_1a - P_1a$; and ${\sf part}(\varphi')$ is
finer than ${\sf part}(\chi)$ since $Q_1b$ is partitioned into $P_1b$,
$P_2b$, $\ldots$, $P_sb$.
 
\medskip

\noindent 
{\bf (2)} \ The proof is just a simpler version of the proof of (1). 
 From $\varphi$ we can construct $\chi$ with an essentially smaller domain
or an essentially coarser partition, while leaving the image unchanged.
And from $\psi$ we can construct $\chi$ with an essentially larger domain
(if ${\sf domC}(\psi)$ is not a maximal prefix code) or with an essentially
finer partition (if $\psi$ is not injective), while leaving the image
unchanged.
 \ \ \ $\Box$

\medskip

\noindent 
We observe that the proof of Prop.\ \ref{L_order_density_Rclass} also shows
the following:

{\it If $\varphi \equiv_{\cal R} \psi$ and $\varphi >_{\cal L} \psi$, and in
addition, ${\sf Dom}(\psi) =_{\sf ess} {\sf Dom}(\varphi)$, then there
exists $\chi$ such that \ $\varphi \equiv_{\cal R} \chi \equiv_{\cal R} \psi$,
 \ $\psi <_{\cal L} \chi <_{\cal L} \varphi$, \ and in addition,
 \ ${\sf Dom}(\chi) =_{\sf ess} {\sf Dom}(\psi)$. }

The next proposition shows that for idempotents it is not possible to have
all these properties at the same time.

\begin{pro} \label{idempot_nonorder} \
If $\eta_0, \eta_1 \in M_{k,1}$ are idempotents such that
${\sf Dom}(\eta_0) =_{\sf ess} {\sf Dom}(\eta_1)$ and
${\sf Im}(\eta_0) =_{\sf ess} {\sf Im}(\eta_1)$, then either $\eta_0 = \eta_1$
or $\eta_0$ and $\eta_1$ are incomparable in the $\cal L$-order.
\end{pro}
{\bf Proof.} Since ${\sf Dom}(\eta_0) =_{\sf ess} {\sf Dom}(\eta_1)$ we can
take essentially equal restrictions so that 
${\sf domC}(\eta_0) = {\sf domC}(\eta_1)$.
By Lemma \ref{idempot_characteriz_1}(3), all idempotents of 
${\sf riHom}(A^*)$ belong to ${\sf riHom}_{\sf pc}(A^*)$; so all essentially 
equal restrictions of idempotents are essentially equal class-wise 
restrictions.

If ${\sf part}(\eta_0)$ and ${\sf part}(\eta_1)$ are not comparable for 
refinement then $\eta_0$ and $\eta_1$ are not $\cal L$-comparable; so, let us
assume now that ${\sf part}(\eta_0)$ refines ${\sf part}(\eta_1)$. If 
${\sf part}(\eta_0) = {\sf part}(\eta_1)$ then $\eta_0 \equiv_{\cal L} \eta_1$,
which in combination with ${\sf Im}(\eta_0) =_{\sf ess} {\sf Im}(\eta_1)$
implies that $\eta_0 = \eta_1$.
So we now assume that ${\sf part}(\eta_0)$ strictly refines 
${\sf part}(\eta_1)$. Let $P = {\sf domC}(\eta_0) = {\sf domC}(\eta_1)$.
Let $Q_1$ be a class of ${\sf part}(\eta_1)$ in $P$, which is the union of 
at least two classes $P_1, \ldots, P_n$ of ${\sf part}(\eta_0)$ in $P$.
Then $\eta_1(Q_1) = p_i \in P_i \subset Q_1$ for some $i = 1, \ldots, n$.
So all ends in $\eta_1(Q_1A^*)$ have $p_1$ as a prefix, whereas ends in
$\eta_0(Q_1A^*)$ can have prefixes in $P_j$ with $j \neq i$. This implies that
${\sf Im}(\eta_0) \neq_{\sf ess} {\sf Im}(\eta_1)$, contrary to the 
assumptions. So ${\sf part}(\eta_0)$ cannot strictly refine 
${\sf part}(\eta_1)$.
 \ \ \ $\Box$

\begin{pro} \label{R_order_density_Lclass} {\bf (Density of the $\cal R$-order
within an $\cal L$-class).}

\noindent {\bf (1)} \ For any two elements $\varphi, \psi \in M_{k,1}$ with
$\varphi \equiv_{\cal L} \psi$ and $\varphi >_{\cal R} \psi$, there
exists $\chi \in M_{k,1}$ such that 

 \ $\chi \equiv_{\cal L} \varphi \equiv_{\cal L} \psi$ \ and \
$\varphi >_{\cal R} \chi >_{\cal R} \psi$.

\noindent {\bf (2)} \ For every $\varphi \in M_{k,1} - \{{\bf 0}\}$  there
exists $\chi \in M_{k,1}$ such that
 \ $\chi \equiv_{\cal L} \varphi$ \ and \ $\varphi >_{\cal R} \chi$.

And for every $\psi \in M_{k,1}$ such that $\psi$ is not $\cal R$-maximal,
there exists $\chi \in M_{k,1}$ such that 
 
 \ $\chi \equiv_{\cal L} \psi$ \ and \ $\chi >_{\cal R}  \psi$.
\end{pro}
{\bf Proof, (1).} \ Since $\varphi \equiv_{\cal L} \psi$ we have, after an 
essentially equal restriction if necessary: 
 \ ${\sf domC}(\varphi) = {\sf domC}(\psi)$, \ and
 \ the ${\sf part}(\varphi)$ on ${\sf domC}(\varphi)$ is equal to 
${\sf part}(\psi)$ on this domain code; let us denote this partition by
$\{P_1, \ldots, P_m\}$. 
Let  $v_i = \psi(P_i)$ for $1 \leq i \leq m$. 

We essentially restrict $\psi$ to $\psi'$ such that ${\sf part}(\psi')$ 
on ${\sf domC}(\psi')$ is $\{P_1a, P_1b, P_2, \ldots, P_m\}$.
Since $\varphi >_{\cal R} \psi$, there are ends in ${\sf Im}(\varphi)$ that
do not belong to ${\sf Im}(\psi)$. Let $q \in A^*$ be a prefix of such an 
end; we can choose $q$ to be longer than any element of 
${\sf imC}(\varphi) \cup {\sf imC}(\psi')$.

We define $\chi$ such that ${\sf domC}(\chi) = {\sf domC}(\psi')$, and
${\sf part}(\chi)$ on ${\sf domC}(\chi)$ is equal to ${\sf part}(\psi')$,
i.e., 

\smallskip

${\sf domC}(\chi) \ = \ \{P_1a, P_1b, P_2, \ldots, P_m\}$. 

\smallskip

\noindent Hence $\chi \equiv_{\cal L} \psi'$.
The values of $\chi$ are defined by

\medskip

$\left\{ \begin{array}{l}
 \chi(P_i) \ = \ v_i \ \ {\rm for} \ 2 \leq i \leq m , \\   
 \chi(P_1b) \ = \ v_1 , \\  
 \chi(P_1a) \ = \ q .
\end{array} \right. $

\medskip

\noindent Then $\psi' <_{\cal R} \chi$ since ${\sf Im}(\psi)$ does not
have any end with prefix $q$. And $\chi <_{\cal R} \varphi'$, since
$q$ is longer than any element of ${\sf imC}(\varphi')$, hence 
${\sf Im}(\varphi')$ has ends that are missing in ${\sf imC}(\chi)$.

\medskip

{\bf (2)} \ The proof is just a simpler version of the proof of (1);
compare with \ref{L_order_density_Rclass}(2).
 \ \ \ $\Box$


\section{Complexity of the $\cal R$-order }

We are interested in the computational difficulty of deciding whether
$\psi \leq_{\cal R} \varphi$ or $\psi \equiv_{\cal R} \varphi$.
We assume at first that $\psi, \varphi \in M_{k,1}$ are given either by 
explicit tables, or by words over a finite generating set of $M_{k,1}$.
Then we consider {\it circuit-like generating sets} of $M_{k,1}$; they 
have the form $\Gamma \cup \tau$, where $\Gamma$ is any finite subset of 
$M_{k,1}$, and  $\tau = \{\tau_{i,i+1} : i \geq 1\}$. The permutation
$\tau_{i,i+1} \in G_{k,1}$ is the letter transposition which swaps positions
$i$ and $i+1$ in any word over $A$. More precisely, 
${\sf domC}(\tau_{i,i+1}) = {\sf imC}(\tau_{i,i+1}) = A^{i+1}$, and

\smallskip

 \ \ \ $\tau_{i,i+1} : $
  \ $u \ell_i \ell_{i+1} \longmapsto u \ell_{i+1} \ell_i $

\smallskip

\noindent for all $u \in A^{i-1}$ and $\ell_i, \ell_{i+1} \in A$.  
Including $\tau$ into the generating set makes word-length polynomially 
equivalent to circuit-size (see \cite{BiCoNP,BiDistor,BiThomMon,BiFact} 
for details, especially Prop.\ 2.4 and Theorem 2.9 in \cite{BiDistor}). 
The word problem of $M_{k,1}$ over $\Gamma$ is in {\sf P}, but over 
$\Gamma \cup \tau$ it is {\sf coNP}-complete \cite{BiThomMon}. 

The {\it word-length} of $\varphi \in M_{k,1}$ over $\Gamma$ or 
$\Gamma \cup \tau$ is the length of a shortest word over $\Gamma$, 
respectively $\Gamma \cup \tau$, that represents $\varphi$; for this,
the length of an element of $\Gamma$ is counted as $1$, and the length
of $\tau_{i,i+1} \in \tau$ is counted as $i+1$. 
We denote these word-lengths of $\varphi$ by $|\varphi|_{_{\Gamma}}$,
respectively $|\varphi|_{_{\Gamma \cup \tau}}$. 

General references on combinational circuits and circuit complexity are 
\cite{Wegener}, \cite{Savage}, and chapters 14 and 2 in \cite{Handb}.


\subsection{Deciding \, $\leq_{\cal R}$ \, over a finite generating set
        $\Gamma$ }

\noindent The first problems whose complexity we would like to know are 
as follows.

\smallskip

\noindent $\bullet$ Input: \ $\psi, \varphi \in M_{k,1}$, given by tables, 
or by words over a finite generating set $\Gamma$ of $M_{k,1}$.

\smallskip

\noindent $\bullet$ Question (the {\bf $\cal R$-order decision problem}): 
 \ Does $\psi \leq_{\cal R} \varphi$ hold?

\smallskip

\noindent $\bullet$ Search (the {\bf right-multiplier search problem}): 
 \ If $\psi \leq_{\cal R} \varphi$, find some $\alpha \in M_{k,1}$
such that $\psi (.) = \varphi \, \alpha(.)$; the multiplier $\alpha$ 
should be expressed by a word over $\Gamma$.

\begin{thm} \label{R_order_Decision_probl_fg} \   
The $\cal R$-order decision problem of $M_{k,1}$ is decidable in 
deterministic polynomial time, if inputs are given by tables or by words 
over a finite generating set.
\end{thm}
{\bf Proof.}
By Corollary 4.11 in \cite{BiThomps}, if $\psi$ and $\varphi$ are given by
words over a finite generating set $\Gamma$ of $M_{k,1}$ then
${\sf imC}(\psi)$ and ${\sf imC}(\varphi)$ can be computed (as explicit
lists of words over $A$) in deterministic polynomial time (where time is
taken as a function of the word-lengths of $\psi$ and $\varphi$ over 
$\Gamma$). If $\psi$ and $\varphi$ are given by tables, then
${\sf imC}(\psi)$ and ${\sf imC}(\varphi)$ can be immediately read from
the tables. By the characterization of $\leq_{\cal R}$ of $M_{k,1}$
(Theorem \ref{R_order}) it is now sufficient to solve the following 
problem in deterministic polynomial time.

\smallskip

\noindent The {\bf end inclusion problem for finite prefix codes}:

\smallskip

\noindent $\bullet$ Input: \ Two finite prefix codes $P, Q \subset A^*$, 
given explicitly as lists of words. 

\smallskip

\noindent $\bullet$ Question:
 \ ${\sf ends}(QA^*) \subseteq {\sf ends}(PA^*)$ ?

\medskip

By Lemma \ref{equiv_R_order}, this question has several equivalent 
formulations. Recall part (f) of Lemma \ref{equiv_R_order}, which says that
${\sf ends}(QA^*) \subseteq {\sf ends}(PA^*)$ \ iff \ for all $y \in Q$ we 
have: 

\noindent -- \ either $P$ contains a prefix of $y$, 

\noindent -- \ or the tree $T(y,P)$ is saturated (where $T(y,P)$ is the 
subtree of $A^*$ with root $y$ and leaf-set $yA^* \cap P$). 

\noindent This yields the following algorithm. 

\bigskip

\parbox{6.5in}{ 
{\it EndInclusion(Q, P)}  \hspace{1in} $||$ 
 \ {\sl Decide whether} \ ${\sf ends}(QA^*) \subseteq {\sf ends}(PA^*)$.

 \hspace{.1in} {\tt for} $y \in Q$
 
 \hspace{.3in}   {\tt if} $P$ does not contain a prefix of $y$

 \hspace{.5in}     {\tt then} 

 \hspace{.7in}        $L$ {\tt :=} ${\it Leaves}(y,P)$ 
 
 \hspace{.7in}        $N$ {\tt :=} ${\it NonLeaves}(y,P)$

 \hspace{.7in}        {\tt if not} ${\it saturated}(L,N)$ 

 \hspace{.9in}           {\tt then return false and halt}

 \hspace{.1in} {\tt return true and halt.} 
}

\bigskip

\noindent This algorithm uses the following three sub-routines.

\bigskip

\parbox{6.5in}{ 
${\it Leaves}(y,P)$ \hspace{1.5in} $||$ 
 \ {\sl Find the set of leaves of $T(y,P)$, i.e., the set $yA^* \cap P$.}

 \hspace{.1in} {\tt return} the set of all words in $P$ that have $y$ as
                a prefix.
}

\bigskip

\parbox{6.5in}{
${\it NonLeaves}(y,P)$ \hspace{1.3in} $||$ 
 \ {\sl Find the set of non-leaf vertices of $T(y,P)$.}

 \hspace{.1in} {\tt return} the set of all {\em strict} prefixes of words 
  in $P$ that have $y$ as a prefix.
}

\bigskip

\parbox{6.5in}{
${\it saturated}(L,N)$  \hspace{1.3in} $||$  
 \ {\sl Is the tree with leaf set $L$ and non-leaf set $N$ saturated? }

 \hspace{.1in} {\tt for} $x \in N$

 \hspace{.3in}    {\tt if} $(xa_1 \not\in N \cup L \ \ {\rm or} 
             \ \ \ \ldots \ \ \ $ ${\rm or} \ \ xa_k \not\in N \cup L)$
 
 \hspace{.6in}        {\tt then} {\tt return false and halt}

 \hspace{.1in} {\tt return true and halt.}
}

\bigskip

\noindent Since the prefix code $Q$ is given explicitly, the main loop 
``{\tt for} $y \in Q$'' repeats a linear number of times. Checking prefix 
relations between words takes linear or quadratic time (depending on the
model of computation).
Since $y$ and the prefix code $P$ is given explicitly, the functions 
${\it Leaves}(y,P)$, ${\it NonLeaves}(y,P)$, and ${\it saturated}(L,N)$
execute in polynomial time. Hence, each iteration of the main loop
``{\tt for} $y \in Q$'' takes polynomial time.
Note that since every vertex of $T(y,P)$ is a prefix of some element 
of $P$, the number of vertices of $T$ at most \ $\sum_{p \in P} |p|$. 
  \ \ \ $\Box$

\bigskip

Before we show that the right-multiplier search problem is also
solvable in polynomial time, we need more definitions and lemmas.

\medskip

By definition, if $\{(u_1,v_1), \ldots, (u_n,v_n)\}$  is a table for an
element of $M_{k,1}$, the {\it total table size} of that table is
 \ $\sum_{i=1}^n |u_i| + |v_i|$.

\begin{lem} \label{find_word_from_table} {\bf (From tables to generators).}
 \  Let $\Gamma$ be a finite generating set of $M_{k,1}$. The following
search problem is solvable in deterministic polynomial time.

\smallskip

 Input: \ $\varphi \in M_{k,1}$, described by a table of total size $n$.

\smallskip

 Output: \ A word over $\Gamma$ that represents $\varphi$.

\smallskip

\noindent It follows that $|\varphi|_{_{\Gamma}}$ has a polynomial upper
bound in terms of the total table size of $\varphi$.
\end{lem}
{\bf Proof outline.} \ First, if $\varphi \in G_{k,1}$ we apply the
constructions from \cite{BiThomps}: By Prop.\ 3.9 in \cite{BiThomps} we
can factor $\varphi$ as
$\varphi = \beta_{\varphi} \pi_{\varphi} \alpha_{\varphi}$ where
$\pi_{\varphi}$ is a permutation of a finite maximal prefix code, and
$\beta_{\varphi}, \alpha_{\varphi} \in F$; \ $F$ is the Thompson group
consisting of the elements of $G_{2,1}$ that preserve the dictionary order
of $A^*$. Then we factor elements of $F$ into generators (Prop.\ 3.10 in
\cite{BiThomps}, and \cite{CFP}). And we factor $\pi_{\varphi}$, first into
word transpositions, then into generators (Lemma 3.11 in \cite{BiThomps}).
The process takes polynomial time. In \cite{BiThomps} the construction was
done for $G_{2,1}$, but the same method applies to $G_{k,1}$ for all
$k \geq 2$.

Second, if $\varphi \in {\it Inv}_{k,1}$ (the monoid of injective elements
of $M_{k,1}$), we use the proof of Theorem 3.4 and Lemma 3.3 in
\cite{BiThomMon}, which is constructive and provides a polynomial-time
algorithm.

Finally, for a general element $\varphi \in M_{k,1}$ we use the proof of
Theorem 3.5 (especially the Claim) in \cite{BiThomMon}; again, the proof
is constructive, and provides a polynomial-time algorithm.
  \ \ \ $\Box$

\bigskip

\noindent {\bf Remark.}  The converse of Lemma \ref{find_word_from_table} 
is not true: By Prop.\ 4.1 in \cite{BiThomMon}, $M_{k,1}$ has
infinitely many elements whose total table size is exponentially larger
than their word-length.

\bigskip

When $\varphi \in M_{k,1}$ is by a word over a finite generating set
$\Gamma$ of $M_{k,1}$ we are interested in finding an {\it inverse} of
$\varphi$, i.e., an element $\chi \in M_{k,1}$ such that
$\varphi \chi \varphi = \varphi $ and $\chi \varphi \chi = \chi$.

\begin{lem} \label{inverse_search_problem} {\bf (Find an inverse).} \
Let $\Gamma$ be a finite generating set of $M_{k,1}$. For every
$\varphi \in M_{k,1}$ there exists an inverse $\chi \in M_{k,1}$ which
satisfies \

\smallskip

\noindent
{\bf (a)} \ $\varphi \, \chi(.) = {\sf id}_{{\sf Im}(\varphi)}$ \ \ and

\smallskip

  \ \ $\chi \, \varphi(.) = \eta_{\varphi}$ \ where $\eta_{\varphi}$ is
an idempotent such that
 \ ${\sf Dom}(\varphi) = {\sf Dom}(\eta_{\varphi})$,
 \ \ ${\sf part}(\eta_{\varphi}) = {\sf part}(\varphi)$;

\smallskip

\noindent {\bf (b)} \ the number of entries in the table of $\chi$ is
 \ $\|\chi\| = |{\sf imC}(\varphi)|$;

\smallskip

\noindent {\bf (c)} \ the length $\ell(\chi)$ of a longest word in the table
of $\chi$ satisfies \ $\ell(\chi) \leq O(|\varphi|_{_{\Gamma}})$.

\smallskip

\noindent Moreover {\em the inverse search problem}, specified as follows,
is solvable in deterministic polynomial time.

\smallskip

 Input: \ $\varphi \in M_{k,1}$, given by a word over $\Gamma$.

\smallskip

 Output: \ An inverse $\chi$ of $\varphi$ with properties {\bf (a),
(b), (c)}; $\chi$ should be described by a word over $\Gamma$.
\end{lem}
{\bf Proof.} \ By Corollary 4.11 in \cite{BiThomMon} we can compute
${\sf imC}(\varphi)$ as an explicit list of words, in polynomial time.
For each $y \in {\sf imC}(\varphi)$, to define $\chi(y)$ we choose one
element in $\varphi^{-1}(y)$. Any way of doing this will make $\chi$ satisfy
conditions (a), (b), and (c); for (b) and (c) this holds by Corollary
4.7 in \cite{BiThomMon}.

By Corollary 4.15 in \cite{BiThomMon}, we can construct a finite automaton
${\cal B}_y$ in polynomial time, such that ${\cal B}_y$ accepts
$\varphi^{-1}(y)$; the construnction takes as input
the given word over $\Gamma$ that represents $\varphi$, and an automaton
${\cal A}_y$ accepting just $\{y\}$; ${\cal A}_y$ can trivially constructed
from $y$; ${\cal A}_y$ has $|y| + 1$ states. The number of states of
${\cal B}_y$ is $\leq |y| + 1 + O(n)$, where $n$ is the length of the word
over $\Gamma$ that represents $\varphi$.

Then from ${\cal B}_y$ we can quickly pick a word $x$ accepted (i.e.,
$x \in \varphi^{-1}(y)$), e.g. by depth-first search in the transition
graph of ${\cal B}_y$. Then $|x| \leq |y| + 1 + O(n)$. Thus,
$\chi(y) = x \in \varphi^{-1}(y)$ can be chosen so that the inverse search
problem is solved in polynomial time.

It follows that the table of $\chi$ contains
$\leq O(|\varphi|_{_{\Gamma}})$ words of length
$\leq O(|\varphi|_{_{\Gamma}})$.
By Lemma \ref{find_word_from_table}, it follows that
$|\chi|_{_{\Gamma}} \leq O(|\varphi|_{_{\Gamma}})$.
   \ \ \ $\Box$

\medskip

\noindent As a consequence of Lemma \ref{inverse_search_problem}, the
elements of $M_{k,1}$, when described by words over a finite generating set,
are {\it not one-way functions} (for any reasonable definition of
``one-way function'').

\begin{pro} \label{R_order_Search_probl_fg} \
The right-multiplier search problem of $M_{k,1}$ is solvable in 
deterministic polynomial time, if inputs are given by tables or by words 
over a finite generating set.
\end{pro}
{\bf Proof.} \  Again, by Corollary 4.11 in \cite{BiThomps}, if $\psi$ and 
$\varphi$ are given by words over a finite generating set $\Gamma$ of 
$M_{k,1}$ then $Q = {\sf imC}(\psi)$ and $P = {\sf imC}(\varphi)$ can be 
computed (as explicit lists of words over $A$) in deterministic polynomial 
time. 

By Lemma \ref{R_equiv_id} above, we have   
 \ $\varphi \equiv_{\cal R} {\sf id}_{PA^*}$ \ and
 \ $\psi \equiv_{\cal R} {\sf id}_{QA^*}$.
The latter is equivalent to \ $\psi = {\sf id}_{QA^*} \circ \psi$.  
Moreover, $\psi \leq_{\cal R} \varphi$ is equivalent to 
 \ ${\sf id}_{QA^*} = {\sf id}_{PA^*} \circ {\sf id}_{QA^*}$.
Thus, after finding $\beta \in M_{k,1}$ such that 
 \ ${\sf id}_{PA^*} = \varphi \, \beta$ \ (by using Lemma
\ref{inverse_search_problem}), we will have
 \ $\psi = \varphi \circ \beta \circ {\sf id}_{QA^*} \circ \psi$.
Next, we express $\beta$ and ${\sf id}_{QA^*}$ as words over $\Gamma$ 
(by using Lemma \ref{find_word_from_table}). For $\varphi$ we already have a
word over $\Gamma$ as part of the input.  This yields a word for 
$\varphi \circ \beta \circ {\sf id}_{QA^*}$ in polynomial time, 
and this is a right-multiplier (since 
$\psi = \varphi \circ \beta \circ {\sf id}_{QA^*} \circ \psi$).
 \ \ \ $\Box$


\subsection{Deciding $\leq_{\cal R}$ over a circuit-like generating set 
  $\Gamma \cup \tau$ } 

Let $\Gamma$ be any finite generating set of $M_{k,1}$ and let 
$\tau = \{ \tau_{i,i+1} : i \geq 1\}$ be the set of transpositions of 
neighboring letters. We saw in \cite{BiCoNP, BiThomMon, BiDistor}  that
the set $\tau$ plays an important role in the representation of 
combinational circuits by elements of $M_{k,1}$ (in such a way that 
circuit size is polynomially related to word-length). For words over the
generating set $\Gamma \cup \tau$ we define the length by $|\gamma| = 1$ for 
$\gamma \in \Gamma$, and $|\tau_{i,i+1}| = i+1$. For 
$\varphi = \gamma_m \cdot \ldots \cdot \gamma_1$ we define
 \ $|\varphi| = \sum_{j=1}^m |\gamma_j|$.

We saw in \cite{BiThomMon} that the word problem of $M_{k,1}$ is 
{\sf coNP}-complete when the generating set $\Gamma \cup \tau$ is used to 
write words (and this is even true for $G_{k,1}$ over $\Gamma_G \cup \tau$, 
where $\Gamma_G$ is any finite generating set of $G_{k,1}$, \cite{BiCoNP}). 
At first impression one might think that the $\leq_{\cal R}$ decision 
problem of $M_{k,1}$ over $\Gamma \cup \tau$ might be in $\Sigma_2^{\sf P}$, 
since $\psi \leq_{\cal R} \varphi$ holds iff 
$\exists \alpha \forall x \, [\psi(x) = \varphi \alpha(x)]$; however, 
there is no guarantee that $\alpha$ has polynomial word-length. But, 
surprisingly at first, the problem turns out to be $\Pi_2^{\sf P}$-complete; 
the proof uses the
characterization of the $\cal R$-order (Theorem \ref{R_order}).   
The connection with $\Sigma_2^{\sf P}$ reappears in Subsection 5.3.

The complexity class $\Pi_2^{\sf P}$ consists of all decision problems 
that can be decided by alternating polynomial-time Turing machines of 
type $\forall \exists$, i.e., nondeterministic polynomial-time Turing 
machines whose computations first visit universal states and then 
existential states; see e.g.\ \cite{DuKo,HemaspOgi,Handb,Papadim} for 
details. The class $\Pi_2^{\sf P}$ contains interesting complete problems, 
e.g., the problem $\forall \exists${\sf -QBF} (also denoted by 
$\forall \exists${\sf Sat}), which asks whether a given 
$\forall \exists$-quantified boolean formula is true; 
see e.g. \cite{DuKo} pp.\ 84-89, \cite{HemaspOgi} pp.\ 270-274. More 
precisely, $\forall \exists${\sf -QBF} consists of all fully quantified 
boolean formulas of the form 
 \ $\forall x \exists y \, \beta(y, x)$, where $x$ and $y$ are finite 
sequences of boolean variables (each boolean variable ranging over 
$\{0,1\}$), and $\beta(y, x)$ is a boolean formula (in the usual sense).

We also consider two special versions of the $\cal R$-order decision problem,
called the ``lower-bound $\cal R$-order decision problem'' and the
``upper-bound $\cal R$-order decision problem'' for $M_{k,1}$ over
the generating set $\Gamma \cup \tau$. 
First we choose an element $\alpha \in M_{k,1}$, called the {\it 
bound-parameter} of the problem. The problems are then specified as follows.  

\smallskip

\noindent $\bullet$ Input: \ $\varphi \in M_{k,1}$, given by a word over 
$\Gamma \cup \tau$.

\smallskip

\noindent $\bullet$ Question ({\bf lower-bound $\cal R$-order 
decision problem for a fixed $\alpha$}): 
 \ \  $\varphi \geq_{\cal R} \alpha$ ? 

\smallskip

\noindent $\bullet$ Question ({\bf upper-bound $\cal R$-order 
decision problem for a fixed $\alpha$}): 
 \ \ $\varphi \leq_{\cal R} \alpha$ ? 

\medskip 

\noindent We consider also the following problems for $M_{k,1}$ over the
generating set $\Gamma \cup \tau$. 

\smallskip

\noindent $\bullet$ Input: \ $\varphi, \psi \in M_{k,1}$, given by words over
$\Gamma \cup \tau$. 

\smallskip

\noindent $\bullet$ Question ({\bf $\equiv_{\cal R}$-decision problem}):
 \ Does $\varphi  \equiv_{\cal R} \psi$ hold?

\smallskip

\noindent $\bullet$ Input: \ $\varphi \in M_{k,1}$, given by a word over
$\Gamma \cup \tau$. 

\smallskip

\noindent $\bullet$
Question ({\bf $\equiv_{\cal R} {\bf 1}$ decision problem}):
 \ Does $\varphi  \equiv_{\cal R} {\bf 1}$ hold?

\medskip

\noindent And we consider the membership problems of the domain, the 
domain code, and the image of an element $\varphi \in M_{k,1}$ 
(over $\Gamma \cup \tau$), specified as follows.

\smallskip

\noindent $\bullet$ Input: \ $\varphi \in M_{k,1}$, given by a word over
$\Gamma \cup \tau$, and a word $z \in A^*$.

\smallskip

\noindent $\bullet$ Question ({\bf domain membership problem}): 
 \ \ $z \in {\sf Dom}(\varphi)$ ?

\smallskip

\noindent $\bullet$ Question ({\bf domain code membership problem}):
 \ \ $z \in {\sf domC}(\varphi)$ ?

\smallskip

\noindent $\bullet$ Question ({\bf image membership problem}): 
 \ \ $z \in {\sf Im}(\varphi)$ ?

\medskip 

\noindent We also consider the special version of the image problem, 
namely the {\bf image membership problem for a fixed test string}.
We first choose a word $z_0 \in A^*$, called the {\it test string} or the
{\it test parameter} of the problem; for a fixed $z_0$ we consider: 

\smallskip

\noindent $\bullet$ Input: \ $\varphi \in M_{k,1}$, given by a word over
$\Gamma \cup \tau$.

\smallskip

\noindent $\bullet$ Question: \ \ $z_0 \in {\sf Im}(\varphi)$ ?

\bigskip

\noindent {\bf Notation:} \ For $x_2, x_1 \in A^*$ we write 
 \ $x_1 \ {\sf pref} \ x_2$ \ iff \ $x_1$ {\it is a prefix of} $x_2$.
We also write \ $x_1 \ {\sf spref} \ x_2$ \ iff \ $x_1$  is a prefix of 
$x_2$ and $x_1 \neq x_2$ (i.e., $x_1$ is a {\em strict prefix} of $x_2$).

\begin{pro} \label{complexity_domC_imC} \hspace{-.08in}{\bf .}

\noindent
{\bf (1)} \  The {\em domain} membership problem, and the {\em domain code} 
membership problem for $M_{k,1}$ over $\Gamma \cup \tau$ are in 
{\sf P}. \\      
{\bf (2)} \ The {\em image} membership problem for $M_{k,1}$ over 
$\Gamma \cup \tau$ is {\sf NP}-complete. The special image membership 
problem (with a fixed test string) is always in {\sf NP}, and it is 
{\sf NP}-complete for certain test strings.
\end{pro} 
{\bf Proof.} (1a) \ The domain membership problem:

Let $\gamma_N \cdot \ldots \cdot \gamma_1$ be the generator sequence that
represents $\varphi$, with $\gamma_j \in \Gamma \cup \tau$ 
($N \geq j \geq 1$). We simply apply this generator sequence to $z$, and 
check if the result is always defined.
By Theorem 4.5(2) in \cite{BiThomMon}, applied to 
$\gamma_n \cdot \ldots \cdot \gamma_1 \cdot {\sf id}_z$, we have for all
$r = 1, \ldots, n$:

\smallskip

$|\gamma_r \cdot \ldots \cdot \gamma_1 \cdot {\sf id}_z| \ \leq $
$ \ |z| + \sum_{j=1}^r \ell(\gamma_j)$,  

\smallskip

\noindent Here, $\ell(\gamma_j)$ is the length of a longest word in the 
table of $\gamma_j$. For every $\gamma_j \in \Gamma$,
$\ell(\gamma_j) \leq c$ for some constant $c$. 
For $\gamma_j = \tau_{i,i+1}$, $\ell(\gamma_j) = i+1$. It follows that 
$\sum_{j=1}^r \ell(\gamma_i) \leq c \ |\varphi|_{_{\Gamma \cup \tau}}$.
Hence, all the intermediate words (in $A^*$) obtained as all the generators
$\gamma_j$ are applied, have length 
$\leq O(|z| + |\varphi|_{_{\Gamma \cup \tau}})$, i.e., the length is linearly
bounded in terms of the size of the input $z, \varphi$.
Also, applying a generator to a word of linearly bounded length takes 
polynomial time.   

\smallskip

\noindent (1b) \ The domain code membership problem:
 
As in (1a) we first check whether $z \in {\sf Dom}(\varphi)$; if the answer
is ``no'' then, obviously, $z \not\in {\sf domC}(\varphi)$. If ``yes'', we 
check whether the prefix $s$ of $z$, obtained by removing the right-most 
letter, is in ${\sf Dom}(\varphi)$; if $s \in {\sf Dom}(\varphi)$ then 
$z \not\in {\sf domC}(\varphi)$, otherwise $z \in {\sf domC}(\varphi)$.
 
\smallskip

\noindent (2) \ The image membership problem, in general, and with fixed 
                test string:

The problems are in {\sf NP}. Indeed, to check whether 
$z \in {\sf Im}(\varphi)$  we can guess $x \in A^*$ and check whether 
$\varphi(x) = z$. By the reasoning in (1) above, the length of $x$ is
linearly bounded from above by the input length 
$|z| + |\varphi|_{_{\Gamma \cup \tau}}$;
moreover, $\varphi(x)$ can be computed in polynomial time by successive 
application of the generators $\gamma_j$.

To show {\sf NP}-hardness we reduce the satisfiability problem for boolean 
formulas to the image membership problem with fixed test string (namely the 
test string $1 \in \{0,1\}^*$). Here we identify $\{0,1\}$ with 
$\{a_1, a_2\} \subseteq \{a_1, \ldots, a_k\} = A$.
Let $\beta(x)$ be a boolean formula with list of boolean variables
$x = (x_1, \ldots, x_m)$. So, $\beta$ also represents a function
$\{0,1\}^m \to \{0,1\}$. Then $\beta(x)$ is satisfiable iff 
$1 \in {\sf Im}(\beta)$. A boolean formula can be viewed as an element
of $M_{k,1}$, with word length over $\Gamma \cup \tau$ linearly related
to the formula size (see \cite{BiDistor}, Proposition 2.4). 
Hence the satisfiability problem for boolean formulas reduces to the 
special image membership problem.
 \ \ \ $\Box$

\medskip

It is interesting that, complexity-wise, the membership problems of the 
domain and the image are quite different.
This is ultimately the cause of the discrepancy between the complexities
of the $\cal R$- and $\cal L$-orders.

\begin{lem} \label{R_order_Decision_infg_inPi2} \   
The $\cal R$-order decision problem of $M_{k,1}$ over $\Gamma \cup \tau$ 
is in $\Pi_2^{\sf P}$.
\end{lem}
{\bf Proof.} Let 
$\psi, \varphi \in M_{k,1}$ be given by words over $\Gamma \cup \tau$,
and let \ $n \ = \ $
${\sf max}\{|\psi|_{_{\Gamma \cup \tau}}, |\varphi|_{_{\Gamma \cup \tau}}\}$.
By Theorem 4.5(2) in \cite{BiThomMon}, the length of a longest word in
${\sf imC}(\psi)$ or ${\sf imC}(\varphi)$ satisfies 
 \ $\ell({\sf imC}(\psi)) \leq O(n)$, and 
 \ $\ell({\sf imC}(\varphi)) \leq O(n)$. Also,
 \ $\ell({\sf domC}(\psi)) \leq O(n)$ \ and 
 \ $\ell({\sf domC}(\varphi)) \leq O(n)$. 
Indeed, if $\varphi = \gamma_N \ldots \gamma_1$ (where 
$\gamma_i \in \Gamma \cup \tau$ for $i = 1, \ldots, N$), then 
$\ell(\varphi) \leq \sum_{i=1}^N \ell(\gamma_i)$ and this is 
$\leq O(|\varphi|_{_{\Gamma \cup \tau}})$; note that 
$\ell(\tau_{i-1,i}) = i$, and that we defined $|\tau_{i-1,i}| = i$ in
defining $|\varphi|_{_{\Gamma \cup \tau}}$.  

Recall Theorem \ref{R_order}(6) which says:  
 \ $\psi \leq_{\cal R} \varphi$ \ iff \    
for every $y \in {\sf imC}(\psi)$, either ${\sf imC}(\varphi)$ contains a 
prefix of $y$, or the subtree $T_{y,{\sf imC}(\varphi)}$ is saturated. Here, 
$T_{y,{\sf imC}(\varphi)}$ is the subtree of the tree of $A^*$ with root 
$y$ and leaf-set $yA^* \cap {\sf imC}(\varphi)$. 
So the set of leaves $L(T_{y,{\sf imC}(\varphi)})$ of 
$T_{y,{\sf imC}(\varphi)}$ is characterized by: 
 \ $z \in L(T_{y,{\sf imC}(\varphi)})$ \ iff \ $z \in {\sf imC}(\varphi)$ 
and $y$ is a prefix of $z$. Hence,

\smallskip

$z \in L(T_{y,{\sf imC}(\varphi)})$ \ \ iff \ \ 
$(\exists x \in {\sf domC}(\varphi))$
$[z = \varphi(x) \ \wedge \ y \ {\sf pref} \ \varphi(x)]$.

\smallskip 

\noindent The vertex set of the subtree $T_{y,{\sf imC}(\varphi)}$ consists 
of the words that are a prefix of a word in $L(T_{y,{\sf imC}(\varphi)})$ 
and that have $y$ as a prefix, i.e., 
$z \in V(T_{y,{\sf imC}(\varphi)})$ \ iff 
 \ $(\exists s \in {\sf imC}(\varphi))$
$[y \ {\sf pref} \ z \ {\sf pref} \ s]$.  Hence,

\smallskip

$z \in V(T_{y,{\sf imC}(\varphi)})$ \ \ iff 
 \ \ $(\exists x \in {\sf domC}(\varphi))$
$[ y \ {\sf pref} \ z \ {\sf pref} \ \varphi(x) ]$.

\smallskip 

\noindent Similarly, the set of non-leaf vertices of 
$T_{y,{\sf imC}(\varphi)}$ is characterized by

\smallskip

$z \in V(T_{y,{\sf imC}(\varphi)}) - L(T_{y,{\sf imC}(\varphi)})$ \ iff \  
$(\exists r \in {\sf domC}(\varphi))$ 
$[ y \ {\sf pref} \ z \ {\sf spref} \ \varphi(x) ]$.

\smallskip

\noindent The subtree $T_{y,{\sf imC}(\varphi)}$ is saturated \ iff \  
$(\forall z \in V(T_{y,{\sf imC}(\varphi)}) - L(T_{y,{\sf imC}(\varphi)}))$ 
$[ zA \subseteq V(T_{y,{\sf imC}(\varphi)}) ]$, since $zA$ is
the set of children of $z$ in the tree of $A^*$.
Thus $\psi \leq_{\cal R} \varphi$ is equivalent to 

\medskip

\parbox{6.5in}{
 $(\forall y \in {\sf imC}(\psi))$ 
 $\big[$
  $(\exists t \in {\sf imC}(\varphi)) [ t \ {\sf pref} \ y ]$
  \ \ $\vee$  
    
 \hspace{1.1in} 
 $(\forall z \in V(T_{y,{\sf imC}(\varphi)}) - L(T_{y,{\sf imC}(\varphi)}))$
  $ [ zA \subseteq V(T_{y,{\sf imC}(\varphi)}) ]$
 $\big]$.
}

\medskip

\noindent Since \ $y \in {\sf imC}(\psi)$ \ iff 
 \ $(\exists x \in {\sf domC}(\psi)) [y = \psi(x)]$, the above formula
is equivalent to

\medskip

\parbox{6.5in}{
 $(\forall x \in {\sf domC}(\psi))$
 $\big[$
   $(\exists t \in {\sf imC}(\varphi)) [ t \ {\sf pref} \ \psi(x) ]$
  \ \ $\vee$  
  
 \hspace{1.2in} $(\forall z \in V(T_{\psi(x),{\sf imC}(\varphi)}) $ $-$ 
     $ L(T_{\psi(x),{\sf imC}(\varphi)}) )$
    $[zA \subseteq V(T_{\psi(x),{\sf imC}(\varphi)}) ]$
 $\big]$.
}

\medskip

\noindent We saw that \ $z \in V(T_{\psi(x),{\sf imC}(\varphi)})$ $ - $
$L(T_{\psi(x),{\sf imC}(\varphi)})$ \ iff 
 \ $(\exists r \in {\sf domC}(\varphi))$
$[\psi(x) \ {\sf pref} \ z \ {\sf spref} \ \varphi(x) ]$.
Also, every $z \in V(T_{\psi(x),{\sf imC}(\varphi)}) $ $-$ 
$ L(T_{\psi(x),{\sf imC}(\varphi)})$ has length 
 \ $|z| < \ell({\sf imC}(\varphi)) \leq c \, \ell({\sf domC}(\varphi))$, 
where $c \geq 1$ is a constant (depending on the choice of the finite set 
$\Gamma$).   Thus, the above formula is equivalent to

\medskip

\parbox{6.5in}{
 $(\forall x \in {\sf domC}(\psi))$

 \ \ \ $\big[$ $(\exists s \in {\sf domC}(\varphi))$
       $ [ \varphi(s) \ {\sf pref} \ \psi(x) ]$     \ \ $\vee$ 

 \ \ \ \ \ \ $(\forall z \in A^{\leq c \, \ell({\sf domC}(\varphi))} ) $
   $(\exists r \in {\sf domC}(\varphi) ) $
   $[\psi(x) \ {\sf pref} \ z \ {\sf spref} \ \varphi(r) \ \Rightarrow \ $
   $zA \subseteq V(T_{\psi(x),{\sf imC}(\varphi)}) ]$
 $\big]$.
}

\medskip

\noindent Obviously, $X \Rightarrow Y$ iff ${\sf not}X \vee Y$. Also,
$zA \subseteq V(T_{\psi(x),{\sf imC}(\varphi)})$ \ iff 
 \ $\bigwedge_{i=1}^k (za_i \in V(T_{\psi(x),{\sf imC}(\varphi)}))$. Moreover, 
$za_i \in V(T_{\psi(x),{\sf imC}(\varphi)})$ \ iff 
 \ $(\exists x_i \in {\sf domC}(\varphi))$
$[ \psi(x) \ {\sf pref} \ za_i \ {\sf spref} \ \varphi(x) ]$.
Hence, the above formula is equivalent to

\medskip 

\noindent
\parbox{6.5in}{
 $(\forall x \in {\sf domC}(\psi))$ 

 \ \ \ \ {\bfseries \Large [} $(\exists s \in {\sf domC}(\varphi))$ 
    $[ \varphi(s) \ {\sf pref} \ \psi(x) \big)]$ \ \ \ \  $\vee$ 

 \ \ \ \ \ \ \  $(\forall z \in A^{\leq \ell({\sf domC}(\varphi))} ) $
  $(\exists r \in {\sf domC}(\varphi) ) $
  
  \ \ \ \   \ \ \ \     
   $[ \ {\sf not}[\psi(x) \ {\sf pref} \ z \ {\sf spref} \ \varphi(r)]$
   $ \ \vee \ $
   $(\exists x_1, \ldots, x_k \in {\sf domC}(\varphi)) $
        $\bigwedge_{i=1}^k [\psi(x) \ {\sf pref} \ za_i \ {\sf spref} \ $
               $\varphi(x_i)] \ ]$
 {\bfseries \Large ]}.
}
 
\medskip

\noindent We transform this to a $\forall \exists$-formula by using the
following facts, where $C$ is any formula that does not contain the variable 
$v$. First, we apply 
 \ $C \vee (\forall v \in S)B(v) \ \Leftrightarrow \ $
          $(\forall v \in S) [C \vee B(v)] $. 
Then, three times we apply 
 \ $C \vee (\exists v \in S)B(v) \ \Leftrightarrow \ $
          $(\exists v \in S)[C \vee B(v)] $.
We then obtain the following $\forall \exists$-formula that characterizes 
the $\cal R$-order of $M_{k,1}$:

\medskip

\parbox{6.5in}{
 $\psi \leq_{\cal R} \varphi$ \ \ iff 

\medskip

 $(\forall x \in {\sf domC}(\psi))$
 $(\forall z \in A^{\leq \ell({\sf domC}(\varphi))} ) $
 $(\exists s \in {\sf domC}(\varphi))$
 $(\exists r \in {\sf domC}(\varphi) ) $
 $(\exists x_1, \ldots, x_k \in {\sf domC}(\varphi)) $

  \ \ \  \ \ \ {\bfseries \Large [}
 $[ \varphi(s) \ {\sf pref} \ \psi(x)]$ \ \ $\vee$ \
 ${\sf not}[\psi(x) \ {\sf pref} \ z \ {\sf spref} \ \varphi(r)] $
   $ \ \vee \ $
 $\bigwedge_{i=1}^k $
   $[ \psi(x) \ {\sf pref} \ za_i \ {\sf spref} \ \varphi(x_i)] $
 {\bfseries \Large ]} .
}

\bigskip

\noindent Moreover, it is a $\Pi_2^{\sf P}$ formula when $\psi$ and
$\varphi$ are given by words over $\Gamma \cup \tau$.
Indeed, all the quantified variables are words of polynomially bounded 
lengths, since $\ell({\sf domC}(\psi)) \leq O(n)$ and 
$\ell({\sf domC}(\varphi)) \leq O(n)$, as we saw already.
And all the predicates  that appear in the formula are decidable in 
deterministic polynomial time; indeed (by Prop.\ \ref{complexity_domC_imC}), 
membership in ${\sf domC}(\psi)$ and in ${\sf domC}(\varphi)$ is in {\sf P},
and $\varphi(w)$ and $\psi(w)$ can be computed in polynomial time from $w$ 
(see the proof of Prop.\ \ref{complexity_domC_imC}(1)); moreover, 
obviously, prefix relations can be checked in polynomial time. 
 \ \ \ $\Box$

\bigskip

In order to prove $\Pi_2^{\sf P}$-hardness we will first consider another 
problem that turns out to be related to the $\cal R$-order of $M_{k,1}$, 
namely the {\bf surjectiveness problem} for $M_{k,1}$ over 
$\Gamma \cup \tau$. It is specified as follows:  

\smallskip

\noindent $\bullet$  
Input: \ $\varphi \in M_{k,1}$, given by a word over $\Gamma \cup \tau$. 

\smallskip

\noindent $\bullet$ 
Question: \ Is $\varphi$ surjective (on $A^{\omega}$) ? 

\smallskip

\noindent It is easy to see that in $M_{k,1}$, the surjective elements are
the same thing as the {\bf epimorphisms} of $M_{k,1}$, i.e., elements 
$\varphi \in M_{k,1}$ such that for all $\psi_1, \psi_2 \in M_{k,1}:$ 
 \ $\psi_1 \varphi(.) = \psi_2 \varphi(.) \Rightarrow \psi_1 = \psi_2$.  

\medskip

We will also use the {\bf surjectiveness problem for 
combinational circuits}, specified as follows:

\smallskip

\noindent $\bullet$ Input: \ A combinational circuit $C$.

\smallskip

\noindent $\bullet$ Question: \ Is the input-output function of $C$
surjective?

\bigskip

\noindent The following gives the connection between the surjectiveness 
problem for $M_{k,1}$ and the $\equiv_{\cal R}$-decision problem.

\begin{lem} \label{surj_equi_characteriz} \ 
The following are equivalent for any $\varphi \in M_{k,1}$: 

\smallskip

\noindent {\bf (1)} \ \ \ $\varphi$ is surjective (on $A^{\omega}$), or
equivalently, $\varphi$ is an epimorphism of $M_{k,1}$.

\smallskip

\noindent {\bf (2)} \ \ \ ${\sf imC}(\varphi)$ is a maximal prefix code.

\smallskip

\noindent {\bf (3)} \ \ \ $\varphi \equiv_{\cal R} {\bf 1}$ \ in $M_{k,1}$.

\smallskip

\noindent {\bf (4)} \ \ \ $(\forall y \in A^N)(\exists x \in A^{\leq N})$
 $[y$ {\rm is a prefix of} $\varphi(x)]$, 

\hspace{.2in}  where $N = c \cdot |\varphi|_{_{\Gamma \cup \tau}}$ and 
 $c$ is the constant $c = {\sf max}\{ \ell(\gamma) : \gamma \in \Gamma\}$.
\end{lem}
{\bf Proof.} The equivalence $(1) \Leftrightarrow (2)$ is straightforward 
(see the discussion of ends, ideals, and prefix codes in Subsections 1.2, 
1.3, and 2.1).

We have of course ${\sf Im}({\bf 1}) = A^*$. Moreover, $A^*$ is 
essentially equal to every essential right ideal, and a right ideal is 
essential iff its generating prefix code is maximal. By Theorem 
\ref{R_order}(1)(2), the equivalence $(2) \Leftrightarrow (3)$ then follows. 
 
Let us prove the equivalence $(4) \Leftrightarrow (2)$. By  Theorem 4.5(2) 
in \cite{BiThomMon}, all words in
${\sf imC}(\varphi) \cup {\sf domC}(\varphi)$ have lengths 
$\leq N = c \cdot |\varphi|_{_{\Gamma \cup \tau}}$, where
$c = {\sf max}\{ \ell(\gamma) : \gamma \in \Gamma\}$. 
Thus, $A^NA^* \subseteq {\sf Im}(\varphi)$.
Hence, ${\sf Im}(\varphi)$ is an essential right ideal iff every word in 
$A^N$ has a prefix in ${\sf imC}(\varphi)$.  This holds iff the formula  
 \ $(\forall y \in A^N)(\exists z \in {\sf imC}(\varphi))[y \ {\sf pref} \ z]$ 
 \ is true. Since all words in ${\sf domC}(\varphi)$ 
have lengths $\leq N$, the statement  \, 
$(\exists z \in {\sf imC}(\varphi))[y \ {\sf pref} \ z]$ \, is equivalent
to \ $(\exists x \in A^{\leq N}) [y \ {\sf pref} \ \varphi(x)]$.
Thus, ${\sf Im}(\varphi)$ is an essential right ideal iff
the formula
 \ $(\forall y \in A^N)(\exists x \in A^{\leq N})$
$[y \ {\sf pref} \ \varphi(x)]$ \ is true.  
 \ \ \ $\Box$

\bigskip

We mentioned already that the problem $\forall \exists{\sf -QBF}$ is 
$\Pi_2^{\sf P}$-complete. It will be useful to prove 
$\Pi_2^{\sf P}$-completeness for a slightly special form of 
$\forall \exists${\sf -QBF}, where the input consists of formulas of
the form $\forall y \exists x \beta_1(x,y)$, where $x \in \{0,1\}^M$ and 
$y \in \{0,1\}^N$ (for some $M,N$), and where $\beta_1(x,y)$ is a boolean 
formula that satisfies $\beta_1(1^{M+N}) = 1$. We call this problem 
$\forall \exists{\sf -QBF}_1$. 

\begin{lem} \label{special_forall_exists_QBF} \ 
The problem $\forall \exists{\sf -QBF}_1$ is $\Pi_2^{\sf P}$-complete.
\end{lem}
{\bf Proof.} The problem $\forall \exists{\sf -QBF}_1$ is obviously in
$\Pi_2^{\sf P}$ since $\forall \exists${\sf -QBF} is.
To prove hardness we map any formula $\forall x_2 \exists x_1 B(x_1,x_2)$
(where $x_1 \in \{0,1\}^m$ and $x_2 \in \{0,1\}^n$) to a formula
$\forall b \forall x_2 \exists x_1 \beta(x_1, x_2, b)$, where $b \in \{0,1\}$
and where $\beta$ is defined for all $x_1, x_2$ by 

\smallskip

 \ \ \ \ $\beta(x_1, x_2, 0) \ = \ B(x_1, x_2)$ ,

\smallskip

 \ \ \ \ $\beta(x_1, x_2, 1) \ = \ 1$ .

\smallskip

\noindent Then $\beta(1^{m+n+1}) = 1$.  
Moreover, $\forall b \forall x_2 \exists x_1 \beta(x_1, x_2, b)$ is 
true iff \ $\forall x_2 \exists x_1 B(x_1,x_2)$ is true. Indeed, 
$\forall b \forall x_2 \exists x_1 \beta(x_1, x_2, b)$ \ $\Longleftrightarrow$
 \ $\forall x_2 \exists x_1 \beta(x_1, x_2, 0)$  $ \ \wedge \ $
$\forall x_2 \exists x_1 \beta(x_1, x_2, 1)$ \ $\Longleftrightarrow$
 \ $\forall x_2 \exists x_1 B(x_1, x_2) \ \wedge \ 1$.  
 \ \ \ $\Box$

\begin{thm} \label{Surj_R_equiv_decision_Pi2hard} \hspace{-.13in} 
{\bf .}  
    
\noindent {\bf (1)} \ The surjectiveness problem for combinational circuits 
is $\Pi_2^{\sf P}$-complete. 

\smallskip

\noindent {\bf (2)} \ As a consequence, the following problem for $M_{2,1}$ 
over the generating set $\Gamma_{2,1} \cup \tau$ are $\Pi_2^{\sf P}$-hard 
(where $\Gamma_{2,1}$ is any finite generating set of $M_{2,1}$):
The surjectiveness problem for elements of $M_{2,1}$, the 
$\equiv_{\cal R} {\bf 1}$ decision problem, the $\equiv_{\cal R}$ decision
problem, and the $\leq_{\cal R}$ decision problem.
\end{thm}
{\bf Proof, (1).} It is easy to see that the surjectiveness problem for
combinational circuits is in $\Pi_2^{\sf P}$. Indeed, $C$ is surjective 
 \ iff \ $\forall y \exists x [C(x) = y]$; and for a given $(x,y)$, one 
can check in deterministic polynomial time whether $C(x) = y$. 

To prove hardness we will reduce $\forall \exists{\sf -QBF}_1$ to the the 
surjectiveness problem for combinational circuits. Let $\beta(x,y)$ be 
a boolean formula where $x$ is a sequence of $m$ boolean variables, and 
$y$ is a sequence of $n$ boolean variables.  We map the formula $\beta$ 
to the combinational circuit $C_{\beta,m,n}$ with input-output function 
defined by

\smallskip

 \ \ \ $(x, y) \ \longmapsto \ C_{\beta,m,n}(x, y) \ = \ $
  $\left\{ \begin{array}{ll}
           y  \  & \mbox{if \ $\beta(x, y) = 1$,} \\
           1^n \ & \mbox{if \ $\beta(x, y) = 0$.}
           \end{array} \right. $

\smallskip

\noindent From the formula $\beta(x, y)$ one can easily construct a 
combinational circuit or a word over $\Gamma \cup \tau$ for $C_{\beta,m,n}$. 
Moreover, ${\sf Im}(C_{\beta,m,n}) $
$=$ $ \{1^n\} \cup \{y : \exists x \beta(x,y)\}$ $ = $  
$\{y : \exists x \beta(x,y)\}$; the latter equality comes from the fact that
$\beta(1^{m+n}) = 1$.
Hence, $\forall y \exists x \beta(x,y)$ is true \ iff \    
${\sf Im}(C_{\beta,m,n}) = \{0,1\}^n$, i.e., iff $C_{\beta,m,n}$ is 
surjective.

\smallskip

\noindent {\bf (2)} By Lemma \ref{surj_equi_characteriz}, it follows that 
the $\equiv_{\cal R} {\bf 1}$ decision problem is 
$\Pi_2^{\sf P}$-hard. Since the latter is a special case of the 
$\cal R$-equivalence problem and of the $\leq_{\cal R}$ decision problem, 
these are also $\Pi_2^{\sf P}$-hard. 
 \ \ \ $\Box$

\begin{thm} \label{R_order_Decision_probl_Infg} \  
The $\equiv_{\cal R} {\bf 1}$ decision problem, the 
$\equiv_{\cal R}$ decision problem, and the $\leq_{\cal R}$ decision problem
of $M_{k,1}$ are $\Pi_2^{\sf P}$-complete, if inputs are given by words 
over $\Gamma \cup \tau$.

The lower-bound $\cal R$-order decision problem is always in
$\Pi_2^{\sf P}$, and it is $\Pi_2^{\sf P}$-complete for certain choices of
the lower-bound parameter (if inputs are given by words over 
$\Gamma \cup \tau$, where $\Gamma$ is any finite generating set of $M_{k,1}$).
\end{thm}
{\bf Proof.} By Lemma \ref{surj_equi_characteriz}(3), and by Lemma 
\ref{R_order_Decision_infg_inPi2} the problems are in $\Pi_2^{\sf P}$.
By Theorem \ref{Surj_R_equiv_decision_Pi2hard}, the $\equiv_{\cal R} {\bf 1}$,
$\equiv_{\cal R}$, and $\leq_{\cal R}$ decision problems are 
$\Pi_2^{\sf P}$-hard for $M_{k,1}$ over $\Gamma \cup \tau$.

Since $\equiv_{\cal R} {\bf 1}$ is equivalent to $\geq_{\cal R} {\bf 1}$,
the lower-bound $\cal R$-order decision problem is 
$\Pi_2^{\sf P}$-hard when {\bf 1} is chosen as the bound-parameter.

The lower-bound $\cal R$-order decision problem is also $\Pi_2^{\sf P}$-hard
for $M_{k,1}$ for each $k \geq 2$, by the same
reasoning as at the end of the proof of Proposition
\ref{complexity_domC_imC}.
 \ \ \ $\Box$

\bigskip

After seeing that the lower-bound decision problem
``$\varphi \geq_{\cal R} {\bf 1}$?'' is $\Pi_2^{\sf P}$-complete (when 
$\varphi$ is given by a word over $\Gamma \cup \tau$), we wonder what might 
be the complexity of upper-bound problems 
``$\varphi \leq_{\cal R} \alpha$?'' for a fixed $\alpha$. Surprisingly we 
have:

\begin{pro} \label{upper_bound_R_decision} \
The upper-bound $\cal R$-decision problem of $M_{k,1}$ over
$\Gamma \cup \tau$ is always in {\sf coNP}, and it is {\sf coNP}-complete
for certain choices of the bound parameter.
\end{pro}
{\bf Proof.} \ The problem is in {\sf coNP}:

Since {\sf coNP} is $\Pi_1^{\sf P}$ it is enough to find a 
$\Pi_1^{\sf P}$-formula that, for a fixed element $\alpha \in M_{k,1}$
and an input $\varphi \in M_{k,1}$, expresses that
$\varphi \leq_{\cal R} \alpha$. Here, $\varphi$ is given by a generator
sequence $\gamma_N \cdot \ldots \cdot \gamma_1$ with $\gamma_j \in $
$\Gamma \cup \tau$, $n \geq j \geq 1$. The length of the longest element 
in ${\sf imC}(\varphi)$ or in ${\sf domC}(\varphi)$
satisfies \ $\ell({\sf domC}(\varphi))$, $\ell({\sf imC}(\varphi)) $
$\leq c \cdot |\varphi|_{_{\Gamma \cup \tau}}$, as we saw in part (1a) of 
the Proof of Prop.\ \ref{complexity_domC_imC}.
By Theorem \ref{R_order}: \ $\varphi \leq_{\cal R} \alpha$ \ iff \
${\sf Im}(\varphi) \ \subseteq_{\sf ends} \ {\sf Im}(\alpha)$.
Also, ${\sf Im}(\alpha) = {\sf imC}(\alpha) \cdot A^*$, where 
${\sf imC}(\alpha)$ is a fixed finite set, and ${\sf imC}(\alpha) \cdot A^*$
is accepted by a finite automaton. We will use the following claims.

\smallskip

\noindent {\sf Claim 1:} \ If $x \in {\sf Dom}(\varphi)$ and
 \, $|x| \geq \ell({\sf domC}(\varphi)) + m$ \, then \,
$|\varphi(x)| \geq m$.

\noindent Indeed, when $\varphi$ is applied to $x$, at most the
$\ell({\sf domC}(\varphi))$ left-most letters of $x$ will be changed;
in particular, if there is a length decrease from $|x|$ to $|\varphi(x)|$,
the amount of decrease will be at most $\ell({\sf domC}(\varphi))$.

\smallskip

\noindent  {\sf Claim 2:} \ Let
 \ $N = \ell({\sf domC}(\varphi)) + \ell({\sf imC}(\alpha))$. Then,

\smallskip

\ \ \ \ ${\sf Im}(\varphi) \ \subseteq_{\sf ends} \ {\sf Im}(\alpha)$
 \ \ iff \ \ $(\forall x \in {\sf Dom}(\varphi) \cap A^{\geq N})$
   $[ \varphi(x) \in {\sf Im}(\alpha) ]$.

\smallskip

\noindent Indeed, if
${\sf Im}(\varphi) \subseteq_{\sf ends} {\sf Im}(\alpha)$ then every word 
$\varphi(x)$ of length $\geq \ell({\sf imC}(\alpha))$ is in
${\sf Im}(\alpha) = {\sf imC}(\alpha) \cdot A^*$. Also, if
$x \in {\sf Dom}(\varphi)$ and $|x| \geq N$ then (by Claim 1),
$|\varphi(x)| \geq \ell({\sf imC}(\alpha))$. Hence,
$\varphi(x) \in {\sf Im}(\alpha)$ for every
$x \in {\sf Dom}(\varphi) \cap A^{\geq N}$.

\noindent Conversely, if
$\varphi({\sf Dom}(\varphi) \cap A^{\geq N}) \subseteq {\sf Im}(\alpha)$,
then ${\sf Im}(\varphi) \ \subseteq_{\sf ends} \ {\sf Im}(\alpha)$, since
the right ideal ${\sf Dom}(\varphi) \cap A^{\geq N}$ is essential in
${\sf Dom}(\varphi)$.
This proves Claim 2.

\smallskip

The ideal ${\sf Dom}(\varphi) \cap A^{\geq N}$ is generated
by the finite prefix code ${\sf Dom}(\varphi) \cap A^N$.
Hence, by Claim 2, the upper-bound relation $\varphi \leq_{\cal R} \alpha$
is characterized by

\medskip

$\varphi \leq_{\cal R} \alpha$  \ \ \ iff
 \ \ \ $(\forall x \in A^N)$ $[ x \not\in {\sf Dom}(\varphi)$
 $ \ \vee \ $  $ \varphi(x) \in {\sf imC}(\alpha) \cdot A^* ]$.

\medskip

\noindent This is a $\Pi_1^{\sf P}$-formula, since the word-length of the
quantified variable $x$ is linearly bounded; indeed, we saw that
$N \leq c \, |\varphi|_{_{\Gamma \cup \tau}} + c_{\alpha}$, where
$c_{\alpha} = \ell({\sf imC}(\alpha))$, which is a constant. Moreover, the 
predicates in the formula can be decided in deterministic polynomial time.
Indeed, the membership problem of ${\sf Dom}(\varphi)$ is in {\sf P},
and ${\sf imC}(\alpha)$ is a fixed finite set so ${\sf imC}(\alpha) \, A^*$
is a fixed regular language (decided by a finite automaton). Moreover,
$\varphi(x)$ can be computed in deterministic polynomial time, as we saw
in part (1a) of the Proof of Prop.\ \ref{complexity_domC_imC}.

\medskip

\noindent Proof of {\sf coNP}-hardness:

We consider the tautology problem for boolean formulas, i.e., the question
whether $\forall x \, B(x)$ is true; here, $B(x)$ is any boolean formula 
with some list of boolean variable $x = (x_1, \ldots, x_m)$. The tautology
problem is a well known {\sf coNP}-complete problem. The formula $B(x)$ 
determines a function $\{0,1\}^m \to \{0,1\}$, which we also denote by 
$B(.)$, and this function determines an element of $M_{2,1}$, which we will 
denote by $\beta$.  By Prop.\ 2.4 in \cite{BiDistor}, the word-length of 
$\beta$ over $\Gamma \cup \tau$ is linearly bounded by the formula size of
$B(x)$, and an expression of $\beta$ as a word over $\Gamma \cup \tau$
can be found in polynomial time from the formula $B(x)$.

We have: \ $\forall x \, B(x)$ \ iff \ ${\sf Im}(B(.)) = \{1\}$  \ iff 
 \ ${\sf Im}(\beta) \subseteq 1 \, \{0,1\}^*$. Since $B(.)$ is a total
function (everywhere defined on $\{0,1\}^m$), the latter is equivalent to
 \ ${\sf Im}(\beta) \subseteq_{\sf ess} 1 \, \{0,1\}^*$. By the 
characterization of $\leq_{\cal R}$ (Theorem \ref{R_order}), this is 
equivalent to \ $\beta \leq_{\cal R} {\sf const}_{m,1}$, where 
${\sf const}_{m,1}$ is the function $xw \mapsto 1w$ (for all 
$x \in \{0,1\}^m$ and all $w \in \{0,1\}^*$).
So the upper-bound $\cal R$-order decision problem with bound parameter
${\sf const}_{m,1}$ is {\sf coNP}-hard for $M_{2,1}$.
The upper-bound problem is also {\sf coNP}-complete for $M_{k,1}$ for each 
$k \geq 2$, by the same reasoning as at the end of the proof of Proposition 
\ref{complexity_domC_imC}.
 \ \ \ $\Box$

\subsection{The right-multiplier search problem for $M_{k,1}$ over
$\Gamma \cup \tau$ }

\begin{thm} \label{RmultiplierLengthCirc} \hspace{-.11in} {\bf .} \\     
{\bf (1)} \ The lengths of right-multipliers for $\leq_{\cal R}$ are not
polynomially bounded, unless the polynomial hierarchy {\sf PH} collapses.
  
 More precisely, suppose there is a polynomial $p(.)$ such that for all
$\psi, \varphi \in M_{k,1}$ we have the following: 
 \ $\psi \leq_{\cal R} \varphi$ implies that there is a right-multiplier 
$\alpha \in M_{k,1}$ such that $\psi = \varphi \alpha$, and such that 
 \ $|\alpha|_{_{\Gamma \cup \tau}} \ \leq \ $
$p(|\psi|_{_{\Gamma \cup \tau}} + |\varphi|_{_{\Gamma \cup \tau}})$.
Then $\Pi_2^{\sf P} = \Sigma_2^{\sf P} = {\sf PH}$.

\smallskip

\noindent {\bf (2)} \ The lengths of right-inverses of surjective 
elements of $M_{k,1}$ are not polynomially bounded, unless the polynomial 
hierarchy {\sf PH} collapses.

More precisely, suppose there is a polynomial $p(.)$ such that we have:
For every $\varphi \in M_{k,1}$ that is surjective (on $A^{\omega}$) 
there exists a right-inverse $\alpha \in M_{k,1}$ (i.e., 
$\varphi \, \alpha(.) = {\bf 1}$) such that 
 \ $|\alpha|_{_{\Gamma \cup \tau}}$ $ \ \leq \ $
$p(|\varphi|_{_{\Gamma \cup \tau}})$. 
Then $\Pi_2^{\sf P} = \Sigma_2^{\sf P} = {\sf PH}$.

\smallskip

\noindent {\bf (3)} \ The circuit sizes of right-inverses of surjective 
boolean functions are not polynomially bounded, unless 
$\Pi_2^{\sf P} = \Sigma_2^{\sf P} = {\sf PH}$.
\end{thm}
{\bf Proof.}  (1) follows from (2), since by Lemma 
\ref{surj_equi_characteriz}, $\varphi$ is surjective iff 
${\bf 1} \leq_{\cal R} \varphi$.

Proof of (2): Recall that $\varphi$ is surjective iff 
${\bf 1} \equiv_{\cal R} \varphi$, iff $\varphi$ has a right-inverse.
Let us assume for a contradiction that for all 
$\varphi \in M_{k,1}$ with ${\bf 1} \equiv_{\cal R} \varphi$ there exists
a right-inverse $\alpha \in M_{k,1}$ with 
$|\alpha|_{_{\Gamma \cup \tau}} \leq p(|\varphi|_{_{\Gamma \cup \tau}})$.
This would imply that the $\equiv_{\cal R} {\bf 1}$ decision problem is 
in $\Sigma_2^{\sf P}$. Indeed, ${\bf 1} \leq_{\cal R} \varphi$ is 
equivalent to 
 \ $(\exists \alpha) (\forall x)$
$[\varphi \alpha(x) = x \ {\sf or} \ \varphi \alpha(x) = \varnothing]$. 
Here, the lengths of the quantified variables are polynomially bounded 
(as a function of $n = |\varphi|_{_{\Gamma \cup \tau}}$); indeed,
$|\alpha|_{_{\Gamma \cup \tau}} \leq p(n)$ by assumption, and 
$x \in {\sf domC}(\varphi \alpha) \subseteq A^{\leq N}$ for some
$N \leq c \cdot (p(n) + n)$ by Theorem 4.5 in \cite{BiThomMon}. 
And $\varphi \alpha(x)$ can be computed in polynomial time when $x$, 
$\varphi$, and $\alpha$ are given, by part (1a) of the Proof of Prop.\
\ref{complexity_domC_imC}.

Since we also saw in Theorem \ref{R_order_Decision_probl_Infg} that the 
question ``${\bf 1} \equiv_{\cal R} \varphi$?'' is an
$\Pi_2^{\sf P}$-complete problem, it follows that 
$\Pi_2^{\sf P} = \Sigma_2^{\sf P}$. 

The proof of (3) is similar to the proof of (2). 
For circuits, as for elements of $M_{2,1}$ in general, we have: 
 \ $C$ is surjective iff there exists a circuit $\alpha$ such that 
$C \circ \alpha(.) = {\sf id}$. Indeed, if the circuit $C$ is surjective 
then $C$ is also surjective as an element of $M_{2,1}$, hence (by Lemma 
\ref{surj_equi_characteriz}), there exists $\alpha \in M_{2,1}$ such that 
$C \circ \alpha(.) = {\sf id}$. 
Let $C: \{0,1\}^m \to \{0,1\}^n$ and ${\sf id}: \{0,1\}^n \to \{0,1\}^n$,
for some $n,m$. Then we have $\alpha: \{0,1\}^n \to \{0,1\}^m$. 
Since both $C$ and ${\sf id}$ are total (i.e., defined on all inputs in
$\{0,1\}^m$, respectively $\{0,1\}^n$), it follows that $\alpha$ is also
total, so $\alpha$ belongs to the monoid ${\it lep}M_{2,1}$ (studied in 
\cite{BiDistor}), i.e., $\alpha$ is the input-output function of a
combinational circuit. By Prop.\ 2.4 and Theorem 2.9 in \cite{BiDistor},
the circuit-size of $\alpha$ and its word-length over $\Gamma \cup \tau$
(as an element of $M_{2,1}$) are polynomially related.
Hence, if $\alpha$ always had polynomially bounded circuit-size 
then the surjectiveness problem for combinational circuits would be in
$\Sigma_2^{\sf P}$.
We saw that the surjectiveness problem for combinational circuits is 
$\Pi_2^{\sf P}$-complete (Theorem \ref{Surj_R_equiv_decision_Pi2hard}),  
hence $\Pi_2^{\sf P} = \Sigma_2^{\sf P}$.
 \ \ \ $\Box$


\section{Complexity of the $\cal L$-order}

\subsection{Deciding \, $\leq_{\cal L}$ \, over a finite generating set}

\noindent We address the same complexity problems as in the previous 
subsections, but for the $\cal L$-order. The main problems are specified by:

\smallskip

\noindent $\bullet$ Input: \ $\psi, \varphi \in M_{k,1}$, given by tables,
or by words over a finite generating set $\Gamma$ of $M_{k,1}$.

\smallskip

\noindent $\bullet$ Question ({\bf $\cal L$-order decision problem}):
 \ Does $\psi \leq_{\cal L} \varphi$ hold?

\smallskip

\noindent $\bullet$ Search ({\bf left-multiplier search problem}): 
 \ If $\psi \leq_{\cal L} \varphi$, find some $\alpha \in M_{k,1}$
such that $\psi (.) = \alpha \, \varphi(.)$; $\alpha$ should be expressed
by a word over $\Gamma$.

\begin{thm} \label{L_order_Decision_probl_fg} \hspace{-.08in}{\bf .} 

\noindent
{\bf (1)} \ The $\cal L$-order decision problem of $M_{k,1}$ is decidable 
in deterministic polynomial time, if inputs are given by tables or by 
words over a finite generating set. \\   
{\bf (2)} \ The left-multiplier search problem of $M_{k,1}$ is
solvable in deterministic polynomial time, if inputs are given by tables or 
by words over a finite generating set.
\end{thm}
{\bf Proof.}  {\bf (1)} Given $\varphi$, we can find an inverse $\varphi'$ 
in deterministic polynomial time, satisfying
$\varphi' \, \varphi(.) = \eta_{\varphi}$, with $\eta_{\varphi}$ as in
Lemma \ref{inverse_search_problem}. In particular, $\eta_{\varphi}$ is
an idempotent; and from $\varphi' \, \varphi(.) = \eta_{\varphi}$ it 
follows that $\varphi \equiv_{\cal L} \eta_{\varphi}$. Also, by Coroll.\ 
\ref{inverse_search_problem}, 
$|\varphi'|_{_{\Gamma}} \leq O(|\varphi|_{_{\Gamma}})$,
hence (since $ \eta_{\varphi} = \varphi' \, \varphi$) we also have
$|\eta_{\varphi}|_{_{\Gamma}} \leq O(|\varphi|_{_{\Gamma}})$.
Similarly, for $\psi$ we have an inverse $\psi'$ with all the properties 
of Lemma \ref{inverse_search_problem}, so  $\psi' \, \psi =  \eta_{\psi}$,
$\eta_{\psi}$ is an idempotent, $\psi \equiv_{\cal L} \eta_{\psi}$,
$|\eta_{\psi}|_{_{\Gamma}} \leq O(|\psi|_{_{\Gamma}})$.

So, $\psi \leq_{\cal L} \varphi$ iff 
$\eta_{\psi} \leq_{\cal L} \eta_{\varphi}$, and the latter holds iff
$\eta_{\psi} = \eta_{\psi}  \eta_{\varphi}$ (since they are idempotents). 
The question ``$\eta_{\psi} = \eta_{\psi}  \eta_{\varphi}$?'' is an 
instance of the word problem of $M_{k,1}$ over $\Gamma$.
Since the word problem of $M_{k,1}$ over $\Gamma$ is in {\sf P}, as was 
proved in \cite{BiThomMon}, the $\leq_{\cal L}$-decision problem is in
{\sf P}. 

\smallskip

{\bf (2)} By Lemma \ref{inverse_search_problem} and the proof of (1) 
above, we have 
 \ $\psi = \psi \,\eta_{\psi} = \psi \, \eta_{\psi} \, \eta_{\varphi} = $
$\psi \, \eta_{\psi} \, \varphi' \, \varphi = \psi \, \varphi' \, \varphi$. 
So, $\psi \, \varphi'$ serves as a left multiplier, and by Lemma
\ref{inverse_search_problem}, $\psi \, \varphi'$  can be found (as a word 
over $\Gamma$) in deterministic polynomial time.
 \ \ \ $\Box$

\subsection{Deciding $\leq_{\cal L}$ over a circuit-like generating set
  $\Gamma \cup \tau$ }

We saw that the $\leq_{\cal R}$-decision problem is $\Pi_2^{\sf P}$-complete.
The characterization of $\leq_{\cal L}$ (Theorem \ref{L_order}) is more 
complicated than the characterization of $\leq_{\cal R}$. Nevertheless we
will see that the $\leq_{\cal L}$-decision problem is easier than the the 
$\leq_{\cal R}$-decision problem over $\Gamma \cup \tau$ (assuming that
{\sf NP} is not equal to {\sf coNP}).
Before we deal with the $\leq_{\cal L}$-decision problem we consider related 
problems that are of independent interest.

\bigskip

Let {\bf 0} denote the zero element of $M_{k,1}$ (represented by the empty
function). We will consider a special case of the word problem of 
$M_{k,1}$, called the {\bf 0 word problem}.

\smallskip

\noindent $\bullet$ Input: \ $\varphi \in M_{k,1}$, given by a word over
the generating set $\Gamma \cup \tau$,

\smallskip

\noindent $\bullet$ Question:
 \ Is \ $\varphi = {\bf 0}$ \ in $M_{k,1}$ ?

\begin{pro} \label{0_word_probl_infgen} \
The {\bf 0} word problem of $M_{k,1}$ over $\Gamma \cup \tau$ is
{\sf coNP}-complete.
\end{pro}
{\bf Proof.} We reduce the tautology problem for boolean formulas to the
{\bf 0} word problem. Let $B$ be any boolean formula, with corresponding
boolean function $\{0,1\}^m \to \{0,1\}$. We identify $\{0,1\}$ with 
$\{a_1,a_2\} \subseteq \{a_1, \ldots, a_k\} = A$.
The function $B$ can be viewed as an element $\beta \in M_{k,1}$, 
represented by a word over $\Gamma \cup \tau$. The length of that word is
linearly bounded by the size of the formula $B$ (by Prop.\ 2.4 in 
\cite{BiDistor}). In $M_{k,1}$ we consider the element ${\sf id}_{0 A^*}$
(i.e., the identity function restricted to $0A^*$), and we assume that 
some fixed representation of ${\sf id}_{0 A^*}$ by a word over $\Gamma$ 
has been chosen. We have:

\smallskip

 \ \ \ ${\sf id}_{0 A^*} \circ \beta(.) = {\bf 0}$ 
 \ \ iff \ \ ${\sf Im}(\beta) \subseteq  1 \, A^*$.

\smallskip

\noindent The latter holds iff $B$ is a tautology. Thus we reduced the
tautology problem for $B$ to the special word problem
 \ ${\sf id}_{0A^*} \, \beta = {\bf 0}$. Note that ${\sf id}_{0A^*}$ is
fixed, and independent of $B$. 

It follows that the {\bf 0} word problem of $M_{k,1}$ over
$\Gamma \cup \tau$ is {\sf coNP}-hard for all $k \geq 2$. Moreover,
since the word problem of $M_{k,1}$ over $\Gamma \cup \tau$ is in
{\sf coNP} (by \cite{BiThomMon}), it follows that the {\bf 0} word problem 
is {\sf coNP}-complete.
 \ \ \ $\Box$

\begin{lem} \label{L_upperbound_decision} \ 
The $\equiv_{\cal L} {\bf 0}$ decision problem for $M_{k,1}$ over 
$\Gamma \cup \tau$ is {\sf coNP}-hard. Hence, the $\cal L$ upper-bound 
decision problem is {\sf coNP}-hard for certain choices of the upper-bound.
\end{lem}
{\bf Proof.} This follows from the {\sf coNP}-hardness of the ${\bf 0}$ 
word problem of $M_{k,1}$ over $\Gamma \cup \tau$ 
(Prop.\ \ref{0_word_probl_infgen}), since $\varphi \equiv_{\cal L} {\bf 0}$ 
iff $\varphi = {\bf 0}$ in $M_{k,1}$. Moreover, this is an $\cal L$
upper-bound decison problem since $\varphi = {\bf 0}$  iff  
$\varphi \leq_{\cal L} {\bf 0}$.
 \ \ \ $\Box$


\begin{lem} \label{L_vs_injective} \
An element $\varphi \in M_{k,1}$ is an {\em injective total} function 
(on the ends space $A^{\omega}$) \ iff \ $\varphi \equiv_{\cal L} {\bf 1}$.  
\end{lem}
{\bf Proof.} Suppose $\varphi \equiv_{\cal L} {\bf 1}$, and let 
$\alpha \, \varphi = {\bf 1}$. If  $\varphi(w)$ were not defined (for some 
$w \in A^{\omega}$), then $\alpha \, \varphi(w)$ would also be undefined;
but $\alpha \, \varphi(w) = {\bf 1}(w) = w$, which is defined. So $\varphi$ 
is total. If $\varphi(w_1) = \varphi(w_2)$ for some $w_1, w_2 \in A^{\omega}$
then $w_1 = w_2$ (by application of $\alpha$ on the left). So $\varphi$ is
injective. 

Conversely, if $\varphi$ is total and injective on $A^{\omega}$ then 
${\sf Dom}(\varphi)$ is an essential right ideal, and the partition of 
$\varphi$ is the identity. Thus, $\varphi$ has the same partition as {\bf 1} 
so, by Theorem \ref{L_order}, $\varphi \equiv_{\cal L} {\bf 1}$.
 \ \ \ $\Box$.

\medskip

It is easy to see the total injective elements of $M_{k,1}$ are the same 
thing as the {\bf monomorphisms} of $M_{k,1}$, i.e., the elements 
$\varphi \in M_{k,1}$ such that for all $\psi_1, \psi_2 \in M_{k,1}:$ 
 \ $\varphi \psi_1(.) = \varphi \psi_2(.) \Rightarrow \psi_1 = \psi_2$. 
By Lemma \ref{L_vs_injective}, the $\cal L$-class of {\bf 1} is exactly
the set of monomorphisms of $M_{k,1}$.

\bigskip

As a consequence of Lemma \ref{L_vs_injective}, the following problems about 
acyclic circuits are relevant for the complexity of the $\cal L$-order.

\smallskip

\noindent $\bullet$ Input: \ An acyclic boolean circuit $B$.

\smallskip

\noindent $\bullet$
Question ({\bf injectiveness problem}): \ Is the input-output
function of $B$ injective?

\smallskip

\noindent $\bullet$ Question ({\bf identity problem}): \ Is the input-output
function of $B$ the identity function?

\begin{pro} \label{injectivenessProbl} \ 
The injectiveness problem and the identity problem for acyclic boolean 
circuits are {\sf coNP}-complete.
\end{pro}
{\bf Proof.}  It is easy to see that the injectiveness problem and the 
identity problem are in {\sf coNP}. To show hardness we reduce the 
tautology problem for boolean formulas to the injectiveness problem (and to
the identity problem) for acyclic circuits. Let $B$ be any boolean 
formula with $n$ variables; the formula defines a function 
$\{0,1\}^n \to \{0,1\}$, which we also call $B$. From $B$ we define a new 
boolean function $F: \{0,1\}^{n+1} \to \{0,1\}^{n+1}$ by 

\smallskip

$F(x_1, \ldots, x_n, x_{n+1}) \ = \ $
$ \left\{ \begin{array}{ll}
(x_1, \ldots, x_n, x_{n+1}) \ \ \ \ \ & 
  \mbox{if \ $B(x_1, \ldots, x_n) = 1$ \ or \ $x_{n+1} = 0$} \\   
( \ 0, \ \ldots, \ 0, \ 1) \ \ \ \ \ & \mbox{otherwise.}   \end{array}  
  \right.  $          

\smallskip

\noindent Let us check that $B$ is a tautology iff $F$ is injective, iff
$F$ is the identity function on $\{0,1\}^{n+1}$.
  
When $B(x_1, \ldots, x_n) = 1$ then 
$F(x_1, \ldots, x_n, x_{n+1}) = (x_1, \ldots, x_n, x_{n+1})$.  So, if 
$B$ is a tautology, i.e., $B(x_1, \ldots, x_n) = 1$ always holds, then $F$ 
is the identity function (which is injective) on $\{0,1\}^{n+1}$.

If $B$ is a not a tautology then $B(c_1, \ldots, c_n) = 0$ for some 
$(c_1, \ldots, c_n) \in \{0,1\}^n$, 
hence $F(c_1, \ldots, c_n,1) = (0,\ldots,0,1)$. But we also have 
$F(0,\ldots,0,0) = (0,\ldots,0,1)$. Hence, $F$ is not injective (and hence 
not the identity function).
 \ \ \ $\Box$.
 
\begin{lem} \label{L_lowerbound_decision} \
The $\equiv_{\cal L} {\bf 1}$ decision problem for $M_{k,1}$ over
$\Gamma \cup \tau$ is {\sf coNP}-hard. Hence, the $\cal L$ lower-bound
decision problem is {\sf coNP}-hard for certain choices of the lower-bound.
\end{lem}
{\bf Proof.} By Lemma \ref{L_vs_injective}, for $\varphi \in M_{k,1}$ we
have $\varphi \equiv_{\cal L} {\bf 1}$ iff $\varphi$ is total and injective
on $A^{\omega}$, or equivalently, iff $\varphi$ is total and injective on 
$A^n$ for some $n$. Combinational circuits are a special case of elements
of $M_{k,1}$ given over a generating set $\Gamma \cup \tau$, and 
circuit-size is polynomially related to word-length in $M_{k,1}$ over 
$\Gamma \cup \tau$ (by Prop.\ 2.4 and Theorem 2.9 in \cite{BiDistor}). The
input-output functions of combinational circuits are total. Hence there is
a reduction from the injectiveness problem of combinational circuits to the
$\equiv_{\cal L} {\bf 1}$ decision problem for $M_{k,1}$ over
$\Gamma \cup \tau$. 
Since the injectiveness problem for circuits is {\sf coNP}-complete, the 
Lemma follows. 
 \ \ \ $\Box$

\begin{thm} \label{L_order_Decision_probl_Infg} \
The $\leq_{\cal L}$, the $\equiv_{\cal L}$, and the 
$\equiv_{\cal L} {\bf 1}$ (i.e., the monomorphism) decision problems of 
$M_{k,1}$ are {\sf coNP}-complete, if inputs are given 
by words over $\Gamma \cup \tau$.

The lower-bound and upper-bound $\cal L$-order decision problem are
always in {\sf coNP}, and they are {\sf coNP}-complete for certain choices of
the bound-parameters (if inputs are given by words over $\Gamma \cup \tau$).
\end{thm}
{\bf Proof.} (1) Hardness of the problems follows immediately from Lemmas
\ref{L_upperbound_decision} and \ref{L_lowerbound_decision}.

\smallskip

\noindent
(2) \ To show that the $\cal L$-order decision problem is in {\sf coNP} we 
use the characterization of the $\cal L$-order from Theorem \ref{L_order},
for all $\varphi, \psi \in M_{k,1} :$

\smallskip

$\psi(.) \ \leq_{\cal L} \ \varphi(.)$ \ \ \ iff \ \ \
${\sf part}(\psi) \ \leq_{\sf end} \ {\sf part}(\varphi)$ .
 
\smallskip

\noindent Let 
 \ $N = {\sf max}\{ \ell({\sf domC}(\psi), \ \ell({\sf domC}(\varphi)\}$.
By the Prop.\ \ref{complexity_domC_imC}(1) (based on Theorem 4.5(2) in 
\cite{BiThomMon}), we have \ $N \leq $
$O(|\psi|_{_{\Gamma \cup \tau}} + |\varphi|_{_{\Gamma \cup \tau}})$. 
Let $\Psi$ and $\Phi$ be the restrictions of $\psi$, respectively $\varphi$,
to $A^N A^*$. Then we have

\smallskip

${\sf part}(\psi) \ \leq_{\sf end} \ {\sf part}(\varphi)$ \ 

\smallskip

\noindent iff 

\smallskip

${\sf domC}(\Psi) \subseteq {\sf domC}(\Phi)$ \ ($\subseteq A^N$), \ \ and  
 
${\sf part}(\Psi) \leq {\sf part}(\Phi)$ 
 \ \ \ (where ``$\leq$'' means ``is a coarser partition than'')

\smallskip

\noindent iff  

\smallskip

$(\forall x \in {\sf domC}(\Psi)) [ x \in {\sf domC}(\Phi)]$ \ \ $\wedge$

$(\forall x_1, x_2 \in {\sf domC}(\Psi))$ $[ \Phi(x_1) = \Phi(x_2) \ $
$ \ \Longrightarrow \ $  $ \ \Psi(x_1) = \Psi(x_2)]$.

\smallskip

\noindent This is a $\forall$-formula. Moreover, all arguments have linearly
bounded length $N \leq $
$O(|\psi|_{_{\Gamma \cup \tau}} + |\varphi|_{_{\Gamma \cup \tau}})$.
Membership in ${\sf domC}(\Psi)$ or ${\sf domC}(\Phi)$ can be tested in
deterministic polynomial time, by Prop.\ \ref{complexity_domC_imC}(1).
Also, on input $\varphi$, $\psi$ (as words over $\Gamma \cup \tau$), and 
$x_1, x_2 \in A^N$ we can compute $\Phi(x_i)$ ($ = \varphi(x_i)$) and
$\Psi(x_i)$ ($ = \psi(x_i)$) ($i = 1,2$) in deterministic polynomial time 
(by Prop.\ \ref{complexity_domC_imC}(1)). So, the above is a 
$\Pi_1^{\sf P}$-formula, i.e., the problem is in {\sf coNP}.
 \ \ \ $\Box$



\bigskip

\bigskip




\bigskip

\bigskip

\noindent {\bf Jean-Camille Birget} \\
Dept.\ of Computer Science \\
Rutgers University at Camden \\
Camden, NJ 08102, USA \\
{\tt birget@camden.rutgers.edu}


{\small

\begin{thebibliography}{99}


\bibitem{BiDrelJord} J.C.~Birget, ``The Thompson-Higman monoids $M_{k,1}$:
  the $\cal J$-order, the $\cal D$-relation, and their complexity''. 
  (In preparation.) 

\bibitem{BiThomMon} J.C.~Birget, ``Monoid generalizations of the Richard
  Thompson groups'', {\it J.~of Pure and Applied Algebra}, 213(2) (Feb.\
  2009) 264-278. (Online pre-publication DOI:
  http://dx.doi.org/10.1016/j.jpaa.2008.06.012 )

  (Preprint: \ ArXiv
  http://arXiv.org/abs/math.GR/0704.0189 , April 2007.)

\bibitem{BiDistor} J.C.\ Birget, ``One-way permutations, computational
  asymmetry and distortion'', {\it J.~of Algebra}, 320(11) (Dec.\ 2008)
4030-4062. \ (Online pre-publication DOI:
  http://dx.doi.org/10.1016/j.jalgebra.2008.05.035 )

  (Preprint: \ ArXiv
  http://arxiv.org/abs/0704.1569,  April 2007).

\bibitem{BiFact} J.C.\ Birget, ``Factorizations of the Thompson-Higman
  groups, and circuit complexity'', {\it International J.~of Algebra and
  Computation}, 18.2 (March 2008) 285-320.
  (Preprint: \  ArXiv
  http://arXiv.org/abs/math.GR/0607349, July 2006.)

\bibitem{BiCoNP} J.C.~Birget, ``Circuits, coNP-completeness, and the groups
  of Richard Thompson'', {\it International J.~of Algebra and Computation}
  16(1) (Feb.\ 2006) 35-90.

  (Preprint: \ ArXiv http://arXiv.org/abs/math.GR/0310335, Oct.\ 2003).

\bibitem{BiThomps} J.C.\ Birget, ``The groups of Richard Thompson and
  complexity'', {\it International J. of Algebra and Computation} 14(5,6)
  (Dec.\ 2004) 569-626.

 (Preprint: \ ArXiv http://arXiv.org/abs/math.GR/0204292, Apr.\ 2002).

\bibitem{CFP} J.\ W.\ Cannon, W.\ J.\ Floyd, W.\ R.\ Parry,
 ``Introductory notes on Richard Thompson's groups'',
 {\it L'Enseignement Math\'ematique} 42 (1996) 215-256.

\bibitem{CliffPres} A.H.\ Clifford, G.B.\ Preston, {\it The Algebraic
  Theory of Semigroups}, Vol.\ 1 (Mathematical Survey, No 7 (I))  American
  Mathematical Society, Providence (1961).

\bibitem{Grillet} P.A.\ Grillet, {\it Semigroups, An Introduction to the
  Structure Theory}, Marcel Dekker, New York (1995).

\bibitem{DuKo} D.Z.\ Du, K.I.\ Ko, {\it Theory of computational complexity},
  Wiley (2000). 

\bibitem{HemaspOgi} L.\ Hemaspaandra, M.\ Ogihara, {\it The complexity 
  theory companion}, Springer Verlag (2002).

\bibitem{Hig74} G.\ Higman, ``Finitely presented infinite simple groups'',
  Notes on Pure Mathematics 8, The Australian National University,
  Canberra (1974).

\bibitem{McKTh} R.\ McKenzie, R.\ J.\ Thompson, ``An elementary 
  construction of unsolvable word problems in group theory'', in 
  {\it Word Problems}, (W.\ W.\ Boone, F.\ B.\ Cannonito, R.\ C.\ Lyndon,
  editors), North-Holland (1973) pp.\ 457-478.

\bibitem{Papadim} Ch.\ Papadimitriou, {\it Computational complexity},
  Addison-Wesley (1994).

\bibitem{PerrinPin} D.\ Perrin, J.E.\ Pin, {\it Infinite words}, 
  Elsevier (2004).

\bibitem{Savage} John E.\ Savage, {\it Models of computation}, 
  Addison-Wesley (1998).

\bibitem{Scott} Elizabeth A. Scott, ``A construction which can be used
  to produce finitely presented infinite simple groups'',
  {\it J. of Algebra} 90 (1984) 294-322.

\bibitem{Th0} Richard J. Thompson, Manuscript (1960s).

\bibitem{Th} Richard J. Thompson, ``Embeddings into finitely generated
  simple groups which preserve the word problem'',
  in {\it Word Problems II}, (S.\ Adian, W.\ Boone, G.\ Higman, editors),
  North-Holland (1980) pp.\ 401-441.

\bibitem{Handb} J.\ van Leeuwen (editor), {\it Handbook of Theoretical
  Computer Science}, volume {\bf A}, MIT Press and Elsevier (1990).

\bibitem{Wegener} I.\ Wegener, {\it The complexity of boolean functions},
  Wiley/Teubner (1987).




\end{thebibliography}

}
\end{document}